\documentclass{amsart}

\usepackage{amssymb}
\usepackage[margin=1.15in]{geometry}
\usepackage[dvipsnames]{xcolor}
\usepackage{amsmath}
\usepackage{amssymb}
\usepackage{amsthm}
\usepackage[shortlabels]{enumitem}
\usepackage{mathdots}
\usepackage{tikz-cd}
\usepackage{subcaption}
\usepackage{mathrsfs}
\usepackage{booktabs}

\usepackage{todonotes}
\usepackage{wrapfig}

\usepackage{array}
\usepackage{tabularx}
\usepackage{url} 

\usepackage{pgf,tikz,pgfplots}

\usepackage{longtable}

\usepackage{setspace}

\makeatletter
\def\paragraph{\@startsection{paragraph}{4}%
  \z@\z@{-\fontdimen2\font}%
  {\normalfont\itshape}}
\makeatother

\usetikzlibrary{shapes}
\usepackage{blkarray}
\usepackage{multirow}
\usepackage{mathtools}
\usepackage{graphicx}

\usepackage{hyperref}
\hypersetup{hyperfootnotes=false,bookmarksnumbered,colorlinks=true, allcolors={RoyalBlue}}

\usepackage{theoremref}

\allowdisplaybreaks

\definecolor{light-gray}{gray}{0.95}

\let\amsamp=&

\newcommand{\pcoor}[1]{%
  \begingroup\lccode`~=`: \lowercase{\endgroup
  \edef~}{\mathbin{\mathchar\the\mathcode`:}\nobreak}%
  [
  \begingroup
  \mathcode`:=\string"8000
  #1%
  \endgroup 
  ]
}

\usepackage{float} 

\makeatletter
\newcommand{\proofstep}[1]{%
  \par
  \addvspace{\medskipamount}
  \textit{#1\@addpunct{.}}\enspace\ignorespaces
}
\makeatother

\theoremstyle{definition}
\newtheorem{defn}{Definition}[section]
\newtheorem{theorem}[defn]{Theorem}

\newtheorem{prop}[defn]{Proposition}
\newtheorem{lemma}[defn]{Lemma}
\newtheorem{cor}[defn]{Corollary}
\newtheorem{rk}[defn]{Remark}

\newtheorem*{stepk+1}{Step $\boldsymbol{k+1}$}

\numberwithin{subcase}{case}
 
\newtheorem{notation}[defn]{Notation}

\newtheorem{example}[defn]{Example}

\theoremstyle{definition}
\newtheorem{innercustomthm}{Theorem}
\newenvironment{customthm}[1]
  {\renewcommand\theinnercustomthm{#1}\innercustomthm}
  {\endinnercustomthm}

\def\l{\lambda}
\def\b{\beta}
\newcommand{\GG}{{\mathbb G }}
 
\newcommand{\zmin}{Z_{\min}}

\newcommand{\End}{\operatorname{End}}

\newcommand{\GL}{\operatorname{GL}}
\newcommand{\SL}{\operatorname{SL}}
\newcommand{\CC}{\mathbb{C}}
\newcommand{\Spec}{\operatorname{Spec}}

\usetikzlibrary{decorations.pathreplacing,arrows}
\tikzset{
  commutative diagrams/.cd,
  arrow style=tikz
}

\tikzset{
  symbol/.style={
    draw=none,
    every to/.append style={
      edge node={node [sloped, allow upside down, auto=false]{$#1$}}}
  }
}

\tikzset{>=to}

\title[Three instability stratifications of the stack of Higgs bundles]{Three instability stratifications of the stack of Higgs bundles on a smooth projective curve}
\author{Eloise Hamilton}

\begin{document}

\begin{abstract}
We study three instability stratifications of the stack of twisted Higgs bundle of a fixed rank and degree on a smooth complex projective curve. The first is the Harder-Narasimhan (HN) stratification, defined by the instability type of the Higgs bundle. The second is the bundle Harder-Narasimhan (bHN) stratification, defined by the instability type of the underlying bundle. While an unstable HN stratum fibres over the stack parametrising Higgs bundles which are isomorphic to their graded, this is not true for Higgs bundles of unstable bHN type. Obtaining such a fibration requires refining the bHN stratification; this is the third instability stratification. After introducing these three stratifications, we establish comparison results. In particular we obtain explicit criteria for determining semistability of a Higgs bundle of low rank with unstable underlying bundle. Then we show how the HN and bHN stratifications can be used to filter the stack of Higgs bundles by global quotient stacks in two different ways. Finally we use these filtrations to relate the HN and bHN stratifications to GIT instability stratifications, and the refined bHN stratification to a Bialynicki-Birula stratification. 
\end{abstract}

\maketitle

\vspace{-20pt}

\section{Introduction} 

\subsection*{The stack of Higgs bundles} 
The classification problem for Higgs bundles, first considered by Hitchin in 1987 \cite{Hitchin1987}, is encoded in the stack $\mathscr{H}_{r,d}(\Sigma)$ of Higgs bundles of rank $r$ and degree $d$ over a compact Riemann surface $\Sigma$. The notion of slope semistability defines an open stratum $\mathscr{H}_{r,d}^{ss}(\Sigma)$ parametrising semistable Higgs bundles \cite{Casalaina-Martin2018}. A quasi-projective moduli space for semistable Higgs bundles can be constructed analytically using gauge theory \cite{Hitchin1987,Fan2020,Fan2022}, or algebraically using Geometric Invariant Theory (GIT) \cite{Nitsure1991,Simpson1994a}. The richness of the geometry of this moduli space has contributed to Higgs bundles becoming a central object of study in mathematics since their introduction; see for example \cite{Garcia-Prada2008,Schaposnik2020, Swoboda2021} for some surveys on the developments in the subject over the last thirty years. 

While the semistable stratum and its associated moduli space are well-studied, less is known about the structure of the stack as a whole and in particular whether moduli spaces with an explicit modular interpretation can be constructed for Higgs bundles which aren't necessarily semistable. 
The aim of this paper is to study the stack of Higgs bundles beyond the semistable stratum, and by doing so to open the door to the construction of moduli spaces for possibly unstable Higgs bundles. 

Stratifications represent a useful tool for studying moduli stacks, one of the reasons being that they break up the stack into pieces whose geometry is typically better behaved than that of the whole stack and hence more amenable to the construction of coarse moduli spaces \cite{Berczi2018}. Stratifications are also at the heart of Halpern-Leistner's `Beyond GIT' programme, according to which a solution to a classification problem should consist in an algebraic stack together with a certain stratification of the stack, called a $\Theta$-stratification \cite{Halpern-Leistner2015}. The results proved in the current paper relate to stratifications of the stack of Higgs bundles, and can be split into three classes: 
\begin{enumerate}[(A)]
\item the definition of three instability stratifications on the stack of Higgs bundles with an explicit moduli-theoretic interpretation (see Theorem \ref{thmA} below);
\item a comparison of the three stratifications, both in arbitrary and low rank, with sharper results obtained in low rank (see \thref{thmB} below);
\item an identification of the strata of each stratification as quotient stacks, and proofs that each stratification can be recovered asymptotically from an algebraic stratification of the parameter space involved in the definition of the quotient stacks (see Theorem \ref{thmC} below).  
\end{enumerate}

The results proved in this paper are all valid in the more general context of $L$-twisted Higgs bundles on $\Sigma$, first considered in \cite{Ramanan1989,Nitsure1991}. For this reason we work throughout with the stack $\mathscr{H}_{r,d}(\Sigma,L)$ of $L$-twisted Higgs bundles instead. Given a line bundle $L$ on $\Sigma$, the Higgs field $\phi$ of an $L$-twisted Higgs bundle $(E,\phi)$ is a map $\phi: E \to E \otimes L$. Classical Higgs bundles are recovered if $L$ is the cotangent bundle. 

\subsection*{Three stratifications}

The stack of Higgs bundles on a smooth projective curve has the particularity that it admits (at least) two instability stratifications with an explicit moduli-theoretic interpretation: a stratification based on the instability type of the Higgs bundle itself (i.e.\ its Harder-Narasimhan type), and a stratification based on the instability type of the underlying bundle (i.e.\ the Harder-Narasimhan type of the underlying bundle). We call the former the Harder-Narasimhan (HN) stratification, and the latter the bundle Harder-Narasimhan (bHN) stratification. 

The HN stratification of the stack of Higgs bundles has appeared in work of Halpern-Leistner \cite{Halpern-Leistner2016}\footnote{In \cite{Halpern-Leistner2016} the HN stratification is considered in the more general case of $L$-twisted $G$-Higgs bundles on a smooth projective curve for $G$ a semisimple group.} and of Gurjar and Nitsure \cite{Gurjar2014}\footnote{In \cite{Gurjar2014} the HN stratification is considered for $\Lambda$-modules on an arbitrary dimensional smooth scheme over a fixed base scheme, and also over a field of possibly mixed characteristic (such objects correspond to classical Higgs bundles if the underlying scheme is a smooth projective curve $\Sigma$ and if $\Lambda = \operatorname{Sym}^{\bullet} (T \Sigma)$).}. By contrast the bHN stratification does not seem to have been studied at the level of the stack of Higgs bundles, although it has been studied at the level of the moduli space of semistable Higgs bundles with the purpose of comparing it with the Bialynicki-Birula stratification of the moduli space induced by the Higgs field scaling action (see \cite{Hausel1998b,Gothen2015,Huang2020,Sancho2022}). 

The HN stratification has the property that each unstable stratum fibres over a product of stacks of semistable Higgs bundles of smaller rank under the map taking a Higgs bundle to its associated graded, defined to be the direct sum of the successive semistable quotients in the HN filtration of the Higgs bundle \cite{Halpern-Leistner2016}. Obtaining the analogous structure of a fibration over moduli stacks of semistable objects for the bHN strata is more subtle. It requires refining the bHN stratification, and doing so gives a third stratification of the stack of Higgs bundles with an explicit moduli-theoretic interpretation. The results we obtain concerning this refinement can be summarised as follows (see Section \ref{sec:thirdstrat} for precise definitions and statements).

\begin{customthm}{A}[The refined bundle Harder-Narasimhan stratification] \label{thmA} 
Let $\tau$ denote a HN type of rank $r$, degree $d$ and length $s$, and $\mathscr{H}_{\tau}^{r,d}(\Sigma,L)$ the stack of Higgs bundles with underlying bundle of HN type $\tau$. Then there is a stratification of $\mathscr{H}_{\tau}^{r,d}(\Sigma,L)$, indexed by suitably defined equivalence classes of pairs $[i,j]$ with $i,j \in \{1,\hdots,s\}$, such that each stratum $\mathscr{H}_{\tau,[i,j]}^{r,d}(\Sigma,L)$ fibres over a product of moduli stacks of semistable vector bundles and holomorphic chains of smaller rank. 
\end{customthm} 

The notion of a `bHN graded' for a Higgs bundle which arises from Theorem \ref{thmA} (see \thref{defn:bHNgraded} for its explicit definition) does not seem to have been explored in the literature, even for semistable Higgs bundles. Nevertheless, the notion we define relates closely to holomorphic chains, which represent an important tool in the study of Higgs bundles through the Higgs field scaling action (see for example \cite{GarciaPrada2011}). The relationship between the bHN graded of a semistable Higgs bundle, as defined in this paper, and its limit under the Higgs field scaling action merits further exploration (see \thref{linkwithcstar}).

\subsection*{Comparison of the stratifications}

The relationship between the HN and bHN stratifications, as well as between the HN and refined bHN stratifications, is difficult to pin down in general: it is not obvious whether the HN filtration of the Higgs bundle can be recovered from the HN filtration of the underlying vector bundle. The comparison results which we establish can be summarised as follows.

\begin{customthm}{B}[Comparison results for the three stratifications] \thlabel{thmB} 
Given HN types $\mu$ and $\tau$, let $\mathscr{H}_{r,d}^{\mu}(\Sigma,L)$ and $\mathscr{H}_{\tau}^{r,d}(\Sigma,L)$ denote the substacks of Higgs bundles of HN type $\mu$ and with underlying bundle of HN type $\tau$ respectively. Then we have: 
\begin{enumerate}[(i)]
\item \label{Bi}  a given stratum $\mathscr{H}_{r,d}^{\mu}(\Sigma,L)$ intersects only a finite number of strata $\mathscr{H}_{\tau}^{r,d}(\Sigma,L)$, and similarly a given stratum $\mathscr{H}_{\tau}^{r,d}(\Sigma,L)$ intersects only a finite number of strata $\mathscr{H}_{r,d}^{\mu}(\Sigma,L)$ (see \thref{intersection}); 
\item \label{Bii} if $r=2$ or $\operatorname{deg} L =0$, then the HN filtration of an unstable Higgs bundle coincides with the HN filtration of its underlying bundle (see \thref{rk2strats}). Moreover, when $r=2$ semistability of a  Higgs bundle with unstable underlying bundle can be established simply by checking the Higgs field invariance of the unique destabilising subbundle (\thref{sscritrk2});
\item if $r=3$, then the HN filtration of a Higgs bundle can be recovered from the HN filtration of its underlying bundle and from the Higgs field (see \thref{rk3strats}). In particular, semistability of a Higgs bundle with unstable underlying bundle can be established simply by checking the Higgs field invariance and slopes of two explicitly defined subbundles (see \thref{sscriterion}); 
\item the refined bHN stratification of $\mathscr{H}_{r,d}(\Sigma,L)$ does not in general coincide with its HN stratification, unless $r=2$. In particular the refined bHN stratum in which a Higgs bundle lies does not in general determine the HN stratum in which it lies (see Section \ref{subsec:intersection}). 
\end{enumerate}
\end{customthm}

Theorem \ref{thmB} \ref{Bi} is a direct generalisation of a result obtained by Nitsure in the case where $\mu$ is the trivial HN type, namely that of a semistable Higgs bundle (see \cite[Prop 3.2]{Nitsure1991}).  Theorem \ref{thmB} \ref{Bii} in the case where $\operatorname{deg} L = 0$ can be viewed as a generalisation of a result by Franco given in \cite{Franco2014} for classical Higgs bundles on an elliptic curve, stating that if a Higgs bundle in this setting is semistable, its underlying bundle is also semistable.

\subsection*{Recovering the stratifications from algebraic stratifications} 

It follows from our comparison results that none of the three stratifications coincide in general. The third class of results in this paper shows that each of the stratifications can be recovered algebraically, either from a GIT instability stratification or from a Bialynicki-Birula stratification. The first step to obtaining these results is an identification of the HN and bHN strata as quotient stacks, and this can be achieved using the set-up from Nitsure's GIT construction of the moduli space of semistable Higgs bundles in \cite{Nitsure1991}.

Vector bundles of rank $r$ and degree $d$ on $\Sigma$ can be parametrised  for $d$ sufficiently large by a subvariety $R_{r,d}$ of the Quot scheme parametrising quotient sheaves of $\mathcal{O}_{\Sigma}^{\oplus m}$ where $m = d+ r(1-g)$ \cite{Newstead1978}. Nitsure constructs a parameter space $F_{r,d}$ for Higgs bundles for $d$ sufficiently large, and which admits a map $f: F_{r,d} \to R_{r,d}$ given by forgetting the Higgs field \cite{Nitsure1991}. The natural action of $\GL(m)$ on $R_{r,d}$ extends to an action on $F_{r,d}$ in such a way that $\GL(m)$-orbits in $F_{r,d}$ correspond to isomorphism classes of Higgs bundles. 

This parameter space can be used to filter the stack of Higgs bundles by quotients stacks in two different ways, via the HN or bHN stratifications respectively. Moreover, the HN and bHN stratifications of the quotient stacks can be recovered from GIT instability stratifications of $F_{r,d}$ and $R_{r,d}$ respectively, while the bHN stratification can be recovered from a Bialynicki-Birula stratification. These results can be summarised as follows. 

\begin{customthm}{C}[Filtrations of $\mathscr{H}_{r,d}(\Sigma,L)$ by quotients stacks and link with algebraic stratifications] \label{thmC} 
Let $\mu$ and $\tau$ denote HN types of rank $r$ and degree $d$, and let $m= d+ r(1-g)$. Let $\mathscr{H}_{r,d}^{\leq \mu}$ and $\mathscr{H}_{\leq \tau}^{r,d}(\Sigma,L)$ denote the substacks of Higgs bundles of HN type less than or equal to $\mu$ and with underlying bundle of HN type less than or equal to $\tau$ respectively. Similarly let $F^{\leq \mu}_{r,d}$ and $F_{\leq \tau}^{r,d} = f^{-1}(R_{r,d}^{\leq \tau})$ denote the corresponding subvarieties of the parameter space $F_{r,d}$. Then for $d$ sufficiently large there are isomorphisms of stacks $$ \mathscr{H}_{r,d}^{\leq \mu}(\Sigma,L) \cong \left[ F_{r,d}^{\leq \mu}/\GL(m) \right] \text{ and } \mathscr{H}^{r,d}_{\leq \tau}(\Sigma,L) \cong \left[ F_{\leq \tau}^{r,d}/\GL(m) \right]. $$ Moreover, the varieties $R_{r,d}^{\leq \tau}$ and $F_{r,d}^{\leq \mu}$ admit $\GL(m)$-equivariant projective completions $\overline{R_{r,d}^{\leq \tau}}$ and $\overline{F_{r,d}^{\leq \mu}}$ with linear $\GL(m)$-actions such that for  $d$ sufficiently large: \begin{enumerate}[(i)]
\item \label{Ci} the HN stratification of $\mathscr{H}_{r,d}^{\leq \mu}(\Sigma,L)$ coincides with the stratification of the quotient stack $[F_{r,d}^{\leq \mu}/\GL(m)]$ induced by the GIT instability stratification associated to the linear action of $\GL(m)$ on $\overline{F_{r,d}^{\leq \mu}}$ (see \thref{HNasGIT});
\item \label{Cii} the bHN stratification of $\mathscr{H}_{\leq \tau}^{r,d}(\Sigma,L)$ coincides with the stratification of the quotient stack $[F^{r,d}_{\leq \tau}/\GL(m)]$ induced by the GIT instability stratification associated to the linear action of $\GL(m)$ on $\overline{R_{r,d}^{\leq \tau}}$ (see \thref{corbHN});
\item the refined bHN stratification of $\mathscr{H}_{\tau}^{r,d}(\Sigma,L)$ coincides with the stratification of $[F_{\tau}^{r,d} /\GL(m)]$ induced by the Bialynicki-Birula stratification associated to the action of a one-parameter subgroup of $\GL(m)$ on $\overline{F_{\leq \tau}^{r,d}} \cap F_{\tau}^{r,d}$ (see \thref{higgstrat}).
\end{enumerate} 
\end{customthm}

Theorem \ref{thmC} \ref{Ci} and \ref{Cii} are generalisations of existing results for vector and Higgs bundles on a smooth projective curve. The first to mention is a result by Seshadri from \cite{Seshadri1967}, proved again by Newstead in \cite{Newstead1978}, which shows that semistability for vector bundles coincides with GIT semitability for the corresponding points in a suitable parameter space. This result is extended to vector bundles of any fixed HN type (i.e.\ not just the trivial HN type associated to semistable bundles) by Kirwan in \cite{Kirw1985}: a correspondence is established between HN types and GIT instability types which satisfies the property that a vector bundle has a given HN type if and only if it lies in the corresponding GIT unstable stratum in a suitable parameter space.\footnote{This result was later generalised to sheaves on a smooth projective variety of arbitrary dimension by Hoskins and Kirwan in \cite{Hoskins2012}.} Theorem \ref{thmC} \ref{Cii} is proved directly from this result. 

Seshadri's result is extended to Higgs bundles by Nitsure in \cite{Nitsure1991}, which shows that semistability for Higgs bundles can be made to coincide with GIT semistability for the associated points in a suitable parameter space. The main result from which Theorem \ref{thmC} \ref{Ci} follows is \thref{bigprop}, which shows that under Kirwan's correspondence between HN types and GIT instability types, a Higgs bundle has a given HN type if and only if it lies in the corresponding GIT unstable stratum in Nitsure's parameter space. This result should be viewed as the Higgs bundle analogue of Kirwan's generalisation to unstable vector bundles of Seshadri's result.

Theorem \ref{thmC} \ref{Ci} and \ref{Cii} show that there is an asymptotic agreement between the HN and bHN stratifications of the stack of Higgs bundles and algebraic stratifications coming from GIT instability stratifications, since the open substacks $\mathscr{H}_{r,d}^{\leq \mu}(\Sigma,L)$ and $\mathscr{H}_{\leq \tau}^{r,d}(\Sigma,L)$ can be made as large as necessary by choosing $\mu$ and $\tau$ suitably large. These comparison results are not the first of their kind: in \cite{Hoskins2014} Hoskins compares HN stratifications to GIT instability stratifications on two different moduli stacks, the stack of quiver representations and the stack of coherent sheaves on a smooth projective scheme, and shows that they coincide when defined in a suitable way. Theorem \ref{thmC} can therefore be viewed as an analogue of these results for the stack of $L$-twisted Higgs bundles on a smooth projective curve.

\subsection*{Links with past and future work} 

The HN stratification for Higgs bundles appears to have been  first considered on an infinite-dimensional space parametrising these objects, rather than on their moduli stack: in \cite{Hausel2004} Hausel and Thaddeus use the HN stratification to study generators for the cohomology ring of the moduli space of semistable Higgs bundles. In \cite{Wilkin2006}, the same stratification is considered by Wilkin who shows that it coincides with the analytic stratification  of the parameter space by the gradient flow of the Yang-Mills-Higgs functional, as a stepping stone to establishing a Morse theory for the Yang-Mills-Higgs functional on the space of Higgs bundles over a compact Riemann surface.

At the level of the stack of Higgs bundles, the HN stratification was first studied by Gurjar and Nitsure in \cite{Gurjar2014}, which shows that Higgs bundles of a fixed HN type form an algebraic stack. This stratification is shown to be a $\Theta$-stratification by Halpern-Leistner in \cite{Halpern-Leistner2016}, and this property along with the theory of derived $\Theta$-stratifications is used to establish the equivariant Verlinde formula for the stack of $L$-twisted $G$-Higgs bundles.

By contrast the bHN stratification of the stack of Higgs bundles, and by extension its refined bHN stratification, do not appear to have been studied at the level of the stack of Higgs bundles, nor on a parameter space for all Higgs bundles. It would be interesting to explore whether these stratifications share similar properties to the HN stratification, in particular in relation to $\Theta$-stratifications, and whether they can be used to shed further light on the moduli space of semsistable Higgs bundles. 

In a different direction, the study of the stack and its stratifications achieved in this paper can also be viewed as a stepping stone to the construction of moduli spaces for Higgs bundles which are not necessarily semistable, thanks to tools from Non-Reductive GIT. This has already been achieved by Jackson in \cite{Jackson2021} in the case of sheaves with a length two HN type on a smooth projective variety. The case of sheaves with arbitrary length HN type is adressed by Hoskins and Jackson in \cite{Hoskins2021}, with a construction given under the assumption of `$\tau$-stability' of the HN type. In future work we will address the construction of moduli spaces for Higgs sheaves of fixed HN type, as well as of fixed bHN type, building on the set-up and results of this paper. We note that in the more general case of Higgs sheaves on an arbitrary dimensional smooth projective variety, a correspondence between the HN stratification and a GIT instability stratification can be established by appealing to the spectral correspondence and the existing correspondence for sheaves given in \cite{Hoskins2012}.

\subsection*{Structure of the paper}

Section \ref{subsec:background} provides some background. Section \ref{sec:stratifications} introduces the three stratifications considered in this paper: the HN, bHN and refined bHN stratifications. Section \ref{sec:comparison} establishes comparison results for the three stratifications. Section \ref{sec:idasqs} describes two filtrations of the stack of Higgs bundles by quotient stacks, using the HN and bHN stratifications respectively. Section \ref{sec:HNGIT} shows how the HN and bHN stratifications of these quotients stacks can be obtained from GIT instability stratifications on the parameter spaces involved in the definition of the quotient stacks, and how the refined bHN stratification can be obtained from a Bialynicki-Birula stratification of a suitable parameter space. 

\subsection*{Conventions}
We work throughout over $\CC$, although the results of this paper remain valid if we replace $\CC$ by an algebraically closed field $k$ of characteristic zero. By variety we mean an integral scheme over $\CC$, separated and of finite type. 

\subsection*{Acknowledgements} This paper builds on work completed during my DPhil under the supervision of Frances Kirwan. I am indebted to her for her help and support throughout. I am grateful also to Gergely B\'erczi, Victoria Hoskins, Joshua Jackson, Yuuji Tanaka and Steven Rayan for many useful discussions.

\section{Background}  \label{subsec:background}

In Section \ref{subsubsec:strat} we define stratifications of an algebraic stack. In Section \ref{subsec:examples} we give two prototypical examples of such stratifications: the Harder-Narasimhan (HN) stratification of the stack of vector bundles and the GIT instability stratification associated to the linear action of a reductive group on a projective variety. In Section \ref{subsubsec:stackofhbs} we define the stack of $L$-twisted Higgs bundles on a smooth projective curve. 

\subsection{Stratification of an algebraic stack}  \label{subsubsec:strat}

Stratifications represent a useful tool for studying the geometry of algebraic stacks; they can be used to break up the stack into pieces whose geometry is more tractable than that of the entire stack. We recall below the definition of a stratification.

\begin{defn}[Stratification of an algebraic stack] \thlabel{defn:stratification}
A \emph{stratification} of an algebraic stack $\mathscr{X}$ is a decomposition $$\mathscr{X} = \bigsqcup_{\lambda \in \Lambda} \mathscr{X}_{\lambda}$$ into substacks $\mathscr{X}_{\lambda}$ indexed by a discrete set $\Lambda$, such that the following properties hold: 
\begin{enumerate}[(i)]
\item the set $\Lambda$ is partially ordered and for each $\lambda \in \Lambda$ the set $\Lambda_{\lambda} : = \{ \lambda' \in \Lambda  \ | \ \lambda' \leq \lambda\}$ is finite; \label{poset}
\item for each $\lambda \in \Lambda$ there is a containment\footnote{In some settings the definition of a stratification involves a stronger condition than that given at \eqref{inclusion}, namely that equality rather than inclusion holds. However as the stratifications studied in this paper do not satisfy this stronger condition, the inclusion is the more natural condition for us.}  \begin{equation} \overline{\mathscr{X}_{\lambda}} \subseteq \mathscr{X}_{\lambda} \sqcup \bigsqcup_{\lambda' > \lambda} \mathscr{X}_{\lambda'},\label{inclusion} 
\end{equation} so that each stratum $\mathscr{X}_{\lambda}$ is a closed substack in $$\mathscr{X}_{\leq \lambda}: = \mathscr{X}_{\lambda} \sqcup \bigsqcup_{\lambda' > \lambda} \mathscr{X}_{\lambda'},$$ which is itself open in $\mathscr{X}$.
\end{enumerate}
\end{defn}

The above definition is the stack-theoretic analogue of the Morse-theoretic stratification introduced by Atiyah and Bott in \cite{Atiyah1983} on the way to studying the moduli space of semistable vector bundles on a smooth projective curve. The stratification they use is the HN stratification, which is the canonical example of a stratification as defined above. We recall it in Section \ref{HNfiltrationvb} below.

\subsection{Examples}  \label{subsec:examples}

\subsubsection{Harder-Narasimhan stratification of the stack of vector bundles} 

\label{HNfiltrationvb} Let $\mathscr{V}_{r,d}(\Sigma)$ denote the stack of vector bundles of rank $r$ and degree $d$ on $\Sigma$; it is an algebraic stack of finite type over $\CC$. It is a classical result that any vector bundle $E$ on $\Sigma$ admits a unique filtration $$0 = E^0 \subseteq E^1 \subseteq \cdots \subseteq E^s = E $$ such that each $E^i/E^{i-1}$ is semistable and such that the slopes $\mu(E^i/E^{i-1})$ are strictly decreasing as $i$ ranges from $1$ to $s$. This canonical filtration is called the \emph{Harder-Narasimhan (HN) filtration} of $E$, and we call the vector $\tau=(d_1/r_1,\hdots,d_1/r_1,d_2/r_2,\hdots,d_s/r_s)$ the \emph{HN type} of $E$, where $d_i/r_i$ is the slope of $E^i/E^{i-1}$ and each $d_i/r_i$ is repeated $r_i$ times. Note that if $E$ has rank $r$ and degree $d$ then $d = d_1 + \cdots + d_s$ and $r = r_1 + \cdots + r_s$, and we say that $\tau$ has rank $r$ and degree $d$, and we call $s$ the \emph{length} of $\tau$. The \emph{HN graded of $E$} is the vector bundle $\operatorname{gr} E = E_1 \oplus \cdots \oplus E_s$ where $E_i = E^i / E^{i-1}$ for each $i = 1,\hdots, s$. 

The HN type induces a stratification of the moduli stack $\mathscr{V}_{r,d}(\Sigma)$ of rank $r$ and degree $d$ vector bundles on $\Sigma$.

\begin{defn}[Harder-Narasimhan stratification of the stack of vector bundles] 
The \emph{HN stratification} of the moduli stack $\mathscr{V}_{r,d}(\Sigma)$ of vector bundles is given by the decomposition \begin{equation} \mathscr{V}_{r,d}(\Sigma) = \bigsqcup_{\tau} \mathscr{V}_{r,d}^{\tau}(\Sigma),
\label{eq:HNstrat}
\end{equation}  where the disjoint union is over HN types $\tau$ of rank $r$ and degree $d$ and $\mathscr{V}_{r,d}^{\tau}(\Sigma) \subseteq \mathscr{V}_{r,d}(\Sigma)$ denotes the substack of vector bundles with Harder-Narasimhan type $\tau$.
\end{defn} 

Note that if $\tau_0$ denotes the HN type $(d/r,\hdots, d/r)$, then $\mathscr{V}_{r,d}^{\tau_0}(\Sigma)$ coincides with the open stratum $\mathscr{V}_{r,d}^{ss}(\Sigma)$ of semistable vector bundles. The HN stratification is a stratification in the sense of \thref{defn:stratification}, with partial order on the set of HN types defined as follows. Each HN type $\tau$ has an associated convex polygon $P_{\tau}$ in the plane, given by connecting the vertices $(d_1/r_1), \hdots, (d_s,r_s)$, and we say that $\tau > \tau'$ if $P_{\tau}$ lies strictly above $P_{\tau'}$. The upper semi-continuity of the HN type with respect to this partial ordering (see \cite[Thm 3, Prop 10]{Shatz1977}) ensures that the HN stratification satisfies the properties of \thref{defn:stratification}. The structure of $\mathscr{V}_{r,d}^{\tau}(\Sigma)$ as an algebraic stack that is a locally closed substack $\mathscr{V}_{r,d}(\Sigma)$ follows from \cite{Nitsure2011}. 

For each HN type $\tau=(d_1/r_1,\hdots, d_s/r_s)$ of rank $r$ and degree $d$, there is a surjective algebraic map $$ \operatorname{gr}: \mathscr{V}_{r,d}^{\tau}(\Sigma)   \to \prod_{i=1}^s \mathscr{V}_{r_i,d_i}^{ss}(\Sigma)$$  
given by $E  \mapsto (E_1,\hdots,E_s)$ where $\operatorname{gr} E: = E_1 \oplus \cdots \oplus E_s$. 
Thus a vector bundle of HN type $\tau$ can be thought of as an $s$-tuple of semistable vector bundles, together with some extension data indicating how the vector bundle can be reconstructed from its graded pieces. This feature is also captured in the fact that the HN stratification of $\mathscr{V}_{r,d}(\Sigma)$ is a $\Theta$-stratification in the sense of \cite{Halpern-Leistner2015,Halpern-Leistner2018}.

\subsubsection{GIT instability stratifications}  \label{subsubsec:GITinstab}
Another prototypical example of a stratification is a GIT instability stratification.\footnote{Such stratifications are sometimes called HKKN stratifications due to Hesselink \cite{Hesselink1978}, Kempf \cite{Kempf1978}, Kirwan \cite{Kirwan1984} and Ness \cite{Ness1984}.} We assume that the algebraic stack is a global quotient stack of the form $[X/G]$ where $G$ is a reductive group acting linearly on a projective variety $X$ with respect to an ample line bundle.

A $G$-invariant stratification of $X$ (which induces a stratification of the quotient stack) can be obtained as follows (see \cite{Kirwan1984} for more detail). By taking a tensor power of the linearisation, we can assume that $G$ acts on $X \subseteq \mathbb{P}(V)$ via a representation $G \to \GL(V)$. We fix a maximal torus $T \subseteq G$ and an invariant inner product on the Lie algebra $\mathfrak{t}$ of $T$, which we use to identify $\mathfrak{t}$ with its dual $\mathfrak{t}^{\vee}$. 

The $V$-representation of $T$ can be diagonalised so that $T$ acts with weight $\alpha_i \in \mathfrak{t}^{\vee}$ on each projective coordinate $x_i$ of $\mathbb{P}(V)$. Given a point $x = [x_0: \cdots : x_n] \in \mathbb{P}(V)$, we let $\operatorname{conv} (x)$ denote the convex hull of the set of weights $\alpha_i$ such that $x_i \neq 0$, viewed as a subset of $\mathfrak{t}$ under the isomorphism $\mathfrak{t}^{\vee} \cong \mathfrak{t}$.

By the Hilbert-Mumford criterion, a point $x = [x_0: \cdots: x_n] \in X$ is semistable for the linear action of $G$ if and only if $\operatorname{conv}(g \cdot x)$ contains the origin for every $g \in G$. If $x$ is unstable, then its `degree of instability' can be measured by how far away from the origin the closest point to the convex hull $\operatorname{conv}(g \cdot x)$ of $g \cdot x$ can be as $g$ ranges over $G$. The GIT instability stratification makes this idea precise. 

Fix a positive Weyl chamber $\mathfrak{t}_{+} \subseteq \mathfrak{t}$ and let $\mathcal{B}$ denote the set of all $\beta \in \mathfrak{t}_{+}$ such that $\beta$ is the closest point to the origin of $\operatorname{conv}(x)$ for some $x \in X$. 

Fix some $\beta \in \mathcal{B} \setminus \{0\}$ and define $$Z_{\beta} : = \{ [x_0:\cdots:x_n] \in X \ | \ x_i =0 \text{ if } \alpha_i \cdot \beta \neq || \beta||^2\}$$ and $$Y_{\beta} : = \left\{ [x_0:\cdots:x_n] \in X \ \biggm| \ \begin{array}{l}
 x_i =0 \text{ if } \alpha_i \cdot \beta < || \beta||^2 \\
 x_i \neq 0 \text{ for some } i \text{ such that } \alpha_i \cdot \beta = || \beta ||^2 \end{array}  \right\}.$$  
There is a retraction $p_{\beta}: Y_{\beta} \to Z_{\beta}$ given by $x_i \mapsto x_i$ if $\alpha_i \cdot \beta = ||\beta||^2$ and $x_i \mapsto 0$ otherwise.

Now let $\operatorname{Stab} \beta \subseteq G$ denote the stabiliser of $\beta \in \mathfrak{t}^{\vee}$ under the co-adjoint action of $G$ on the dual $\mathfrak{g}^{\vee}$ of its Lie algebra $\mathfrak{g}$. Then $Z_{\beta}$ is $\operatorname{Stab} \beta$-invariant and the linear action of $G$ on $X$ induces a linear action of $\operatorname{Stab} \beta$ on $Z_{\beta}$. Identifying $\beta$ as a character of $\operatorname{Stab} \beta$, we let $Z_{\beta}^{ss}$ denote the semistable locus for the above linear action of $\operatorname{Stab} \beta$ on $Z_{\beta}$ twisted by the character $-\beta$, and define $Y_{\beta}^{ss} = p_{\beta}^{-1}(Z_{\beta}^{ss})$. Finally, we set $S_{\beta} := G Y_{\beta}^{ss}.$
The subvarieties $S_{\beta}$ form the strata of the GIT instability stratification of $X$, which is given by \begin{equation} X = \bigsqcup_{\beta \in \mathcal{B}} S_{\beta}, \label{instabilitystrat}
\end{equation} where we set $S_0 = X^{ss}$. Each stratum $S_{\beta}$ for $\beta \neq 0$ satisfies the property that $$\overline{S_{\beta}} \subseteq S_{\beta} \sqcup \bigsqcup_{||\beta'|| >  ||\beta||} S_{\beta'}.$$

We note that each $\beta \neq 0$ has a corresponding one-parameter subgroup $\lambda_{\beta}: \mathbb{G}_m \to T$, determining the following parabolic subgroup of $G$ (see \cite[Prop 2.6]{Mumford1994}): $$P_{\beta}= \left\{ g \in G \ \left| \ \lim_{t \rightarrow  0} \lambda_{\beta}(t)g \lambda_{\beta}(t^{-1}) \in G \right\} \right..$$ By \cite[Thm 13.5]{Kirwan1984}, for each $\beta \neq 0$ there is an isomorphism \begin{equation} 
S_{\beta} \cong G \times_{P_{\beta}} Y_{\beta}^{ss}. \label{unstablestratumisom}
\end{equation} 
The above property of the unstable strata will be used in Section \ref{subsec:refbHNasBB}.

\subsection{The stack of $L$-twisted Higgs bundles} 
\label{subsubsec:stackofhbs}

The stack studied in this paper is the moduli stack of $L$-twisted Higgs bundles of rank $r$ and degree $d$ on a smooth projective curve $\Sigma$. We define it below. 

\begin{defn}[Moduli stack of Higgs bundles]
Fix a smooth projective curve $\Sigma$, a line bundle $L \to \Sigma$ with $h^0(\Sigma,L) > 0$, and integers $r \in \mathbb{Z}_{>0}$ and $d \in \mathbb{Z}$. The definitions of $L$-twisted Higgs bundles of rank $r$ and degree $d$ on $\Sigma$, morphisms of such objects, families of such objects and equivalences of such families are as follows:  
\begin{enumerate}[(i)]
\item an \emph{$L$-twisted Higgs bundle of rank $r$ and degree $d$ on $\Sigma$} is a pair $(E,\phi)$ where $E$ is a vector bundle\footnote{By vector bundle we mean a locally free sheaf of finite rank.} of rank $r$ and degree $d$ on $\Sigma$ and $\phi: E \to E \otimes L$ is a map of vector bundles;
\item a \emph{morphism from a Higgs bundle $(E,\phi)$ to another $(E',\phi')$} is a morphism $\psi: E \to E'$ of the underlying bundles such that the diagram $$ \begin{tikzcd}
E \ar{r}{\psi} \ar{d}{\phi} & E' \ar{d}{\phi'}\\
E \otimes L \ar{r}{\psi \otimes \operatorname{id}_L} & E' \otimes L,
\end{tikzcd} $$ commutes; it is an isomorphism if $\psi$ is an isomorphism of vector bundles;
\item a \emph{family of Higgs bundles parametrised by a scheme $B$} is a pair $(E_B,\phi_B)$ where $E_B \to \Sigma \times B$ is a vector bundle over $\Sigma \times B$ of rank $r$ and fibrewise degree $d$ and $\phi_B: E_B \to E_B \otimes \pi_{\Sigma}^{\ast} L$ is a morphism of bundles where $\pi_{\Sigma}: \Sigma \times B \to \Sigma$ is the natural projection; 
\item two families $(E_B,\phi_B)$ and $(E_{B'},\phi_{B'})$ parametrised by schemes $B$ and $B'$ respectively are \emph{equivalent} if there exists an isomorphism $\psi_B: E_B \to E_{B'}$ such that the following diagram commutes: $$ \begin{tikzcd}
E_B \ar{r}{\psi_B} \ar{d}{\phi_B} & E' \ar{d}{\phi_{B'}}\\
E_B \otimes \pi^{\ast}_{\Sigma} L \ar{r}{\psi_B \otimes \operatorname{id}_{\pi^{\ast}_{\Sigma}L}} & E_{B'} \otimes \pi^{\ast}_{\Sigma} L.
\end{tikzcd} $$
\end{enumerate}
The \emph{stack of $L$-twisted Higgs bundles of rank $r$ and degree $d$ on $\Sigma$}, denoted $\mathscr{H}_{r,d}(\Sigma,L)$, is the moduli stack associated to the above data. That is, as a category fibred in groupoids it associates to a scheme $B$ the category whose objects are families of $L$-twisted Higgs bundles of rank $r$ and degree $d$ parametrised by $B$ and whose morphisms are given by equivalences of families. 
\end{defn}

The stack $\mathscr{H}_{r,d}(\Sigma,L)$ is an algebraic stack, locally of finite type over $\Spec \CC$ (see \cite[Thm 7.15]{Casalaina-Martin2018}). Our aim in this paper is to study stratifications of this stack.

\section{Three stratifications}
 \label{sec:stratifications}

In this section we describe three stratifications on the stack of Higgs bundles, each admitting an explicit moduli-theoretic interpreation. The first is the HN stratification, determined by the instability type of the Higgs bundle; we introduce it in Section \ref{subsec:HNstrat}. The second is the sheaf HN stratification, determined by the instability type of the underlying bundle; we introduce it in Section \ref{subsec:bHNstrat}. The third stratification arises from trying to obtain for the bHN strata the structure of a fibration over a moduli stack parametrising `associated graded' objects, as is the case for the HN strata. This is the refined bHN stratification, which we define in Section \ref{sec:thirdstrat}.

\subsection{The Harder-Narasimhan stratification}  \label{subsec:HNstrat}
A natural instability stratification of the stack of Higgs bundles can be obtained by considering the analogue for Higgs bundles of the HN type for vector bundles. Indeed, the notion of HN filtrations for vector bundles reviewed in Section \ref{HNfiltrationvb} can be extended to Higgs bundles (see for example \cite[\S 3]{Simpson1994}): given any Higgs bundle $(E,\phi)$, there exists a unique filtration $$0 = E^0 \subseteq E^1 \subseteq \cdots \subseteq E^s = E$$ satisfying the following conditions: \begin{enumerate}[(i)]
\item $E^i$ is a $\phi$-invariant subbundle of $E$ for each $i = 1, \hdots, s$; \label{1}
\item $(E_i,\phi_i)$ is semistable for each $i$, where $E_i = E^i/E^{i-1}$ and $\phi_i$ is the map $E_i \rightarrow E_i \otimes L $ induced by $\phi$; \label{2} 
\item $\mu(E_1) > \mu(E_2) > \cdots > \mu(E_s)$ where $E_i := E^i / E^{i-1}$ for $i = 1,\hdots, s$. \label{3} 
\end{enumerate} 
This filtration is the \emph{HN filtration} of $(E,\phi)$. The vector $\mu= (\mu(E_1),\hdots, \mu(E_1),\mu(E_2) , \hdots, \mu(E_s))$, where each $\mu(E_i)$, is repeated $\operatorname{rk}(E_i)$ times is the \emph{HN type} of $(E,\phi)$, which we say has \emph{rank $r$} and \emph{degree $d$}. The \emph{HN graded} of $(E,\phi)$ is the Higgs bundle $\operatorname{gr}(E,\phi) = (E_1,\phi_1) \oplus \cdots \oplus (E_s,\phi_s)$ where $\phi_i:E_i \to E_i \otimes L$ is induced from $\phi$. Just as for vector bundles, the HN type induces a stratification of the stack of Higgs bundles. 

\begin{defn}[HN stratification of $\mathscr{H}_{r,d}(\Sigma,L)$]
The \emph{HN stratification} of the moduli stack $\mathscr{H}_{r,d}(\Sigma,L)$ of Higgs bundles is the decomposition \begin{equation} 
\mathscr{H}_{r,d}(\Sigma,L) = \bigsqcup_{\mu} \mathscr{H}_{r,d}^{\mu} (\Sigma,L), \label{eqn:HNstrat}
\end{equation} where the disjoint union is over HN types of rank $r$ and degree $d$ and $\mathscr{H}_{r,d}^{\mu}(\Sigma,L) \subseteq \mathscr{H}_{r,d}(\Sigma,L)$ denotes the substack of Higgs bundles with HN type $\mu$. 
\end{defn} 

Just like the HN stratification of $\mathscr{V}_{r,d}(\Sigma)$, the HN stratification of $\mathscr{H}_{r,d}(\Sigma,L)$ satisfies the properties of \thref{defn:stratification} (see for example \cite[\S 2.1]{Wentworth2016}). The structure of $\mathscr{H}_{r,d}^{\mu}(\Sigma,L)$ as an algebraic stack that is a locally closed substack of $\mathscr{H}_{r,d}(\Sigma,L)$ follows from \cite[Cor 4.9]{Gurjar2014}. Note that as in the vector bundle case, if we let $\mu_0$ denote the trivial HN type $(d/r,\hdots, d/r)$, then $\mathscr{H}_{r,d}^{\mu_0}(\Sigma,L)$ coincides with the open stratum $\mathscr{H}_{r,d}^{ss}(\Sigma,L)$ of semistable Higgs bundles. 

Again as in the case of vector bundles, given a HN type $\mu = (d_1/r_1,\hdots, d_s/r_s)$ there is an algebraic map $$
\operatorname{gr}: \mathscr{H}_{r,d}^{\mu}(\Sigma,L)  \to \prod_{i=1}^s \mathscr{H}_{r_i,d_i}^{ss}(\Sigma,L) $$
given by $(E,\phi)  \mapsto ((E_1,\phi_1),\hdots, (E_s,\phi_s))$  where $\operatorname{gr}(E,\phi) = (E_1,\phi_1) \oplus \cdots \oplus (E_s,\phi_s)$. A Higgs bundle of HN type $\mu$ is therefore determined by an $s$-tuple of semistable Higgs bundles together with some extension data indicating how the Higgs bundles can be reconstructed from its graded pieces. This feature is also captured by the fact that the HN stratification of $\mathscr{H}_{r,d}(\Sigma,L)$ is a $\Theta$-stratification in the sense of \cite{Halpern-Leistner2015} (see \cite[\S 1.0.6]{Halpern-Leistner2016}). 

\begin{rk}[Non-separatedness of the HN strata]  \thlabel{HNnotsep}
It is an important fact that the stable stratum (or the semistable stratum if $r$ and $d$ are coprime) of $\mathscr{H}_{r,d}^{s}(\Sigma,L)$ is separated; this is what makes the construction of an associated quasi-projective coarse moduli space possible. By contrast the HN strata may not be separated. Indeed, given an unstable Higgs bundle $(E,\phi)$, one can find a family of Higgs bundles parametrised by $\mathbb{A}^1$ such that for each non-zero $t \in \mathbb{A}^1$, the
fibre over $t$ is isomorphic to $(E,\phi)$, while the fibre over $0$ is isomorphic to $\operatorname{gr}(E,\phi)$ (see \cite[Rk 4.5]{Nitsure1991}). Since $(E,\phi)$ may not be isomorphic to $\operatorname{gr}(E,\phi)$, it follows that the HN strata may not be separated. This precludes the possibility of obtaining a separated coarse moduli space for a whole HN stratum in general, as any such coarse moduli space would have to identify a Higgs bundle with its HN graded. 
\end{rk}

\subsection{The bundle Harder-Narasimhan stratification}  \label{subsec:bHNstrat}

A second instability stratification of the stack of Higgs bundles can be obtained by pulling back the HN stratification of $\mathscr{V}_{r,d}(\Sigma)$, the stack of vector bundles, via the natural forgetful map $F: \mathscr{H}_{r,d}(\Sigma,L) \to \mathscr{V}_{r,d}(\Sigma)$ which maps a Higgs bundle to its underlying vector bundle. The map $F$ is a representable morphism of stacks (see \cite[Cor 7.16]{Casalaina-Martin2018}). 

\begin{defn}[Bundle HN stratification of $\mathscr{H}_{r,d}(\Sigma,L)$] \thlabel{defn:bHNstrat}
The \emph{bundle HN stratification} of the moduli stack $\mathscr{H}_{r,d}(\Sigma,L)$ of Higgs bundles is the decomposition \begin{equation} \mathscr{H}_{r,d}(\Sigma,L) = \bigsqcup_{\tau} \mathscr{H}^{r,d}_{\tau}(\Sigma,L) \label{eqn:bHNstrat}
\end{equation} where the disjoint union is taken over HN types of rank $r$ and degree $d$ and $\mathscr{H}^{r,d}_{\tau}(\Sigma,L) : = F^{-1}(\mathscr{V}_{r,d}^{\tau}(\Sigma))$ for each such HN type $\tau$. 
\end{defn}  

When $\tau$ is the trivial bundle HN (bHN) type $\tau_0 = (d/r,\hdots,d/r)$, we denote $\mathscr{H}^{r,d}_{\tau}(\Sigma,L)$ by $\mathscr{H}_{ss}^{r,d}(\Sigma,L)$ since it parametrises Higgs bundles with underlying semistable vector bundle. 

\begin{rk}[Comparison with the Bialynicki-Birula stratification of the moduli space of semistable Higgs bundles]
While the bHN stratification does not appear to have been studied in the literature at the level of the stack of Higgs bundles, it has been considered at the level of the moduli space of semistable Higgs bundles for the purpose of comparing it to the Bialynicki-Birula stratification of the moduli space induced by the Higgs field scaling action. Indeed, in \cite{Hausel1998b} Hausel proves that the two stratifications coincide in rank two, while in \cite{Gothen2015} Gothen and Zu\~niga-Rojas show that this is no longer the case in rank three (the rank four case is dealt with by Ant\'on Sancho in \cite{Sancho2022}). In \cite{Huang2020} Huang compares the bHN stratification of the moduli space of rank three semistable Higgs bundles to its partial oper stratification; the latter is compared to the Bialynicki-Birula stratification of the moduli space in arbitrary rank in \cite{Collier2019}, using the theory of conformal limits.  
\end{rk}

Given a Higgs bundle $(E,\phi)$ of bHN type $\tau$, we would like to define its associated `bHN graded' $\operatorname{sgr}(E,\phi)$, just as we can define the HN graded $\operatorname{gr}(E)$ and $\operatorname{gr}(E,\phi)$ of a vector bundle $E$ and Higgs bundle $(E,\phi)$ respectively. This is because the existence of an associated graded object endows the unstable strata with the structure of a fibration over products of moduli stacks of semistable objects, a useful tool in the study of the unstable strata. 

Since the HN graded of a vector or Higgs bundle has the same HN type as the vector or Higgs bundle itself, we would expect the bHN graded of a Higgs bundle of bHN type $\tau$ to have bHN type $\tau$. It is reasonable therefore to expect the underlying bundle of the bHN graded of $(E,\phi)$ to be $\operatorname{gr} E$.  However the Higgs field of a Higgs bundle $(E,\phi)$ does not necessarily preserve its bHN filtration and so each graded piece $E_i$ need not have an associated Higgs field. For this reason there is no obvious way to define its associated `bHN graded' as a direct sum of Higgs bundles of smaller rank and degree. The situation is therefore not as simple as in the case of the HN strata for vector bundles or the HN strata for Higgs bundles. This issue is the motivation behind the definition of the refined bHN stratification, introduced in Section \ref{sec:thirdstrat} below.

\subsection{The refined bHN stratification} \label{sec:thirdstrat} 

In this section we introduce a third instability stratification of the stack of Higgs bundles: a refinement of the bHN stratification. This refinement arises from trying to define a `bHN graded' for a Higgs bundle with a given bHN type $\tau$. Section \ref{subsec:bHNgraded} is devoted to proposing a definition of the `bHN-graded' of a Higgs bundle of bHN type $\tau$. In Section \ref{subsec:refofbHN} we define the refined bHN stratification and construct the analogue of the surjective algebraic maps $\operatorname{gr}: \mathscr{V}_{r,d}^{\tau}(\Sigma) \to \prod_{i=1}^s \mathscr{V}_{r_i,d_i}(\Sigma)$ and $\operatorname{gr}: \mathscr{H}_{r,d}^{\mu}(\Sigma,L) \to \prod_{i=1}^s \mathscr{H}^{ss}_{r_i,d_i}(\Sigma,L)$ for the refined bHN strata.

 \subsubsection{bHN graded of a Higgs bundle} \label{subsec:bHNgraded} 

The definition of the bHN graded of a Higgs bundle $(E,\phi)$ which we give satisfies the property that its underlying bundle is the HN graded of $E$, namely $\operatorname{gr} E= E_1 \oplus \cdots \oplus E_s$. The Higgs field for the bHN graded of $(E,\phi)$ then needs to be a well-defined Higgs field for $\operatorname{gr} E$, determined roughly by how much of the HN filtration of $E$ is preserved by $\phi$. If the filtration is preserved by the Higgs field, then the underlying bundle of the HN graded $\operatorname{gr}(E,\phi)$ of $(E,\phi)$ coincides with $\operatorname{gr} (E)$ and therefore its Higgs field is a well-defined Higgs field for $\operatorname{gr} (E)$. In this case, the bHN graded of $(E,\phi)$ should coincide with its HN graded. If the HN filtration of $E$ is not preserved by $\phi$, then we would like to extract from $\phi$ a 		`component' which, in a well-defined sense, is `most responsible' for not preserving the HN filtration of $E$, and moreover which gives a well-defined Higgs field for the HN graded $\operatorname{gr} E$ of $E$.  This can be achieved as follows. 

Let $(E,\phi)$ have bHN type $\tau$ and let $0=E^0 \subseteq E^1 \subseteq \cdots \subseteq E^s = E$ denote its bHN filtration. For each $i$ let $\pi_i: E \to E/E^i$ denote the quotient map. By taking a tensor product with $\operatorname{id}_{L} : L \to L$, we obtain a map $E \otimes L \to E/E^i \otimes L$ which we also denote by $\pi_i$ for simplicity. Given $i$ and $j$, consider the map $\pi_{i-1} \circ \phi|_{E^j}: E^j \to E/E^{i-1} \otimes L.$ Then if $\pi_{i-1} \circ \phi|_{E^j}$ restricts to the zero map on $E^{j-1} \subseteq E^j$, the map descends to a map $E^j/E^{j-1} \to E^i/E^{i-1} \otimes L$. If moreover we have that $\pi_{i} \circ \phi|_{E^j}: E^j \to E/E^i \otimes L$ is the zero map, then the image of $\pi_{i-1} \circ \phi|_{E^j}$ is contained in $E^i/E^{i-1} \otimes L$. Under these two assumptions, we obtain a well-defined map $ \phi_{ij}: E^j/E^{j-1} \to E^i/E^{i-1} \otimes L$, which makes the following diagram commute: 
\begin{center} \begin{tikzcd}
E^j \ar["\phi|_{E^j}",r] \ar["\pi_j",d] & E \otimes L \ar["\pi_{i-1}",r] & E/E^{i-1} \otimes L \ar["\supseteq",d] \\
E^j/E^{j-1} \cong E_j \ar["\phi_{ij}",rr] & & E^i / E^{i-1} \otimes L \cong E_i \otimes L.
\end{tikzcd}
\end{center}

The choice of $\tau=(d_1/r_1,\hdots, d_s/r_s)$ induces an ordering on the pairs $(i,j)$ for $i,j \in \{1,\hdots,s\}$ given by $(i,j) \leq (i',j')$ if $d_i/r_i - d_j/r_j \leq d_{i'}/r_{i'} - d_{j'}/{r_{j'}}.$  As a consequence of the strict inequalities $d_1/r_1 > \cdots > d_s/r_s$, we have that $(i,j) > (i',j)$ for any $i'  > i$ and that $(i,j) < (i,j')$ for any $j' > j$. Note that $(s,1)$ is the smallest pair under this ordering, and $(1,s)$ the largest. Given $\tau$, we let $[i,j]_{\tau}$ denote the equivalence class of $(i,j)$ under the equivalence class $(i,j) \sim_{\tau} (i',j')$ if $d_i/r_i - d_j/r_j = d_{i'}/r_{i'} - d_{j'}/r_{j'}$. Then $\tau$ determines a strict ordering of the equivalence classes $[i,j]_{\tau}$. 

\begin{lemma} \thlabel{welldefined} If $(E,\phi)$ is a Higgs bundle of bHN type $\tau = (d_1/r_1,\hdots, d_s/r_s)$, then there is a well-defined smallest equivalence class $[i,j]_{\tau}$ such that $\phi_{kl}$ is well-defined and non-zero for all $(k,l) \in [i,j]_{\tau}$. 
\end{lemma} 

\begin{proof}
It suffices to show that if for a given equivalence class $[i,j]_{\tau}$, the map $\phi_{i'j'}$ is well-defined and equal to $0$ for all $(i',j') < (i,j)$, then $\phi_{kl}$ is well-defined for any $(k,l) \in [i,j]_{\tau}$. We do so by strong induction. For the base case, namely $[i,j] = [s,1]_{\tau} = \{(s,1)\}$, it suffices to show that $\phi_{s1}$ is well-defined. This is true since $\pi_{s-1} \circ \phi|_{E^1}: E^1 \to E/E^{s-1} \otimes L$ trivially induces a map $E_1 \to E_s \otimes L$, given that $E^0=0$ and $E = E^s$. 

As the inductive hypothesis, for a given equivalence class $[i,j]_{\tau}$ suppose that $\phi_{i'j'}$ is well-defined and equal to zero for all $(i',j') < (i,j)$. Then in particular, $\phi_{i+1,j}$ and $\phi_{i,j-1}$ are well-defined and equal to zero. The fact that $\phi_{i+1,j}$ is the zero map implies that $\phi(E^j)$ is contained in $E^i \otimes L$, whilst the fact that $\phi_{i,j-1}$ is the zero map implies that the map $\pi_{i-1} \circ \phi|_{E^j}: E^j \to E/E^{i-1} \otimes L$ restricts to the zero map on $E^{j-1}$. This ensures that $\phi_{ij}: E_j \to E_i \otimes L$ is well-defined. The same argument shows that $\phi_{kl}: E_j \to E_i \otimes L$ is well-defined for any $(k,l) \in [i,j]_{\tau}$. 
\end{proof} 
Note that if this smallest equivalence class $[i,j]_{\tau}$ of \thref{welldefined} satisfies $j \geq i$ (and consequently $l \geq k$ for any $(k,l) \in [i,j]_{\tau}$), then $\phi$ must preserve the HN filtration of $E$ and so the HN and bHN filtrations of $(E,\phi)$ coincide. This leads us to the following definition for the bHN graded of a Higgs bundle.

\begin{defn}[bHN graded of a Higgs bundle] \thlabel{defn:bHNgraded}
Let $(E,\phi)$ denote a Higgs bundle of bHN type $\tau$, and let $\operatorname{gr} E = E_1 \oplus \cdots \oplus E_s$. Let $[i,j]_{\tau}$ denote the equivalence class associated to $(E,\phi)$, as defined in \thref{welldefined}, which we call the \emph{smallest equivalence class associated to $(E,\phi)$}. Let $$\phi_{\operatorname{min}} = \bigoplus_{(k,l) \in [i,j]_{\tau}} \phi_{kl}:  \bigoplus_{(k,l) \in [i,j]_{\tau}} E_l \to \bigoplus_{(k,l) \in [i,j]_{\tau}} E_k \otimes L.$$
 The \emph{bHN graded} of $(E,\phi)$ is the Higgs bundle $$\operatorname{sgr} (E,\phi) : = \begin{cases} \left( \operatorname{gr}(E), \phi_{\operatorname{min}} \right) & \text{ if $i > j$} \\
 \operatorname{gr}(E,\phi) & \text{ otherwise,} 
\end{cases}
$$ 
where the map $\phi_{\operatorname{min}}$ is extended in the natural way to a Higgs field for $\operatorname{gr} (E)$.  
\end{defn}

\begin{example}[bHN graded when $\tau$ has length $2$] Suppose that $\tau$ has length two and let $(E,\phi)$ denote a Higgs bundle of bHN type $\tau$. Let $\operatorname{gr} E = E_1 \oplus E_2$. Then $$\operatorname{sgr}(E,\phi) = \begin{cases} &  \left( E_1 \oplus E_2, \begin{pmatrix} 0 & 0 \\
\phi_{21} & 0 
\end{pmatrix} \right)  \text{ if $\phi_{21} \neq 0$} \\
&  \left( E_1 \oplus E_2, \begin{pmatrix} \phi_{11} & 0 \\
0 & \phi_{22} \end{pmatrix} \right) \text{ if $\phi_{21} = 0$}
\end{cases} $$ 
where $\phi_{ii}: E_i \to E_i \otimes L$ is induced from $\phi$. 
\end{example}

\begin{rk}[Link with the Higgs field scaling action] \thlabel{linkwithcstar}
The isomorphism class of $\operatorname{sgr}(E,\phi)$ is always fixed by the Higgs field scaling action, provided the smallest equivalence class $[i,j]_{\tau}$ associated to $(E,\phi)$ satisfies $i>j$. Indeed in this case $\operatorname{sgr}(E,\phi)$ is a holomorphic chain, and these exactly correspond to Higgs bundles whose isomorphism classes are fixed by the Higgs field scaling action (see \cite{Hausel2003}).
If moreover $(E,\phi)$ is semistable, then an isomorphism class of Higgs bundles fixed by the scaling action can also be obtained by taking the isomorphism class of the semistable limit $\operatorname{lim}_{t \to 0}(E,t\phi)$ under the scaling  action. It is natural to ask how this limit compares to $\operatorname{sgr}(E,\phi)$.

When $r=2$, the two coincide. Indeed, in this case $[i,j]_{\tau} = [2,1]_{\tau}$ and the condition that $\phi_{21} \neq 0$ implies that $(E,\phi)$ is semistable (we will prove this in \thref{rk2strats} in Section \ref{sec:comparison} below), as is $\operatorname{sgr}(E,\phi)$. It follows then from \cite[Prop 3.1]{Gothen2015} that $[\operatorname{sgr}(E,\phi)] = 
[\operatorname{lim}_{t \to 0}(E,t\phi)]$. In higher rank equality may no longer hold. Indeed, in \cite{Gothen2015} Gothen and Zu\~niga-Rojas show that in rank three the underlying bundle of $\operatorname{lim}_{t \to 0} (E,\phi)$ may not be the HN graded of $E$. However by definition the underlying bundle of $\operatorname{sgr}(E,\phi)$ is always the HN graded of $E$. We put forward the following conjecture: if $(E,\phi)$ is semistable and the underlying bundle of $\operatorname{lim}_{t \to 0} (E,\phi)$ is the HN graded of $E$, then the limit coincides with $\operatorname{sgr}(E,\phi)$.  
\end{rk}

\subsubsection{Refinement of the bundle Harder-Narasimhan stratification} \label{subsec:refofbHN}

The notion of bHN graded for a Higgs bundles established in Section \ref{subsec:bHNgraded} above leads to a natural refinement of the bHN stratification of the stack of Higgs bundles. 

\begin{prop} [Refinement of the bHN stratification of $\mathscr{H}_{r,d}(\Sigma,L)$] \thlabel{Higgsrefinement}
Let $\tau  \neq \tau_0$ denote a bHN type of rank $r$ and degree $d$, and of length $s$. For each $[i,j]_{\tau}$, let $$\mathscr{H}^{r,d}_{\tau,[i,j]}(\Sigma,L) \subseteq \mathscr{H}_{\tau}^{r,d}(\Sigma,L)$$ denote the substack of Higgs bundles $(E,\phi)$ such that the smallest equivalence class associated to $(E,\phi)$ is $[i,j]_{\tau}$. 
Moreover, let $$\mathscr{H}_{\tau,[0,0]}^{r,d}(\Sigma,L) : = \bigsqcup_{1 \leq i \leq j \leq s} \mathscr{H}_{\tau,[i,j]}^{r,d}(\Sigma,L).$$  
Then we have: \begin{equation}
\mathscr{H}_{\tau}^{r,d}(\Sigma,L) =  \bigsqcup_{1\leq j < i \leq s}  \mathscr{H}^{r,d}_{\tau,[i,j]}(\Sigma,L) \sqcup \mathscr{H}^{r,d}_{\tau, [0,0]}. \label{eq:Higgsstrat}
\end{equation} 
The resulting refinement of the bHN stratification of $\mathscr{H}_{r,d}(\Sigma,L)$ is a stratification in the sense of \thref{defn:stratification}. The strata are called the \emph{refined bHN strata} of $\mathscr{H}_{r,d}(\Sigma,L)$. 
Moreover, \begin{equation} \mathscr{H}^{r,d}_{\tau,[0,0]}(\Sigma,L)   = \mathscr{H}^{r,d}_{\tau}(\Sigma,L) \cap \mathscr{H}_{r,d}^{\tau}(\Sigma,L). \label{eq:00}
\end{equation}  
\end{prop}

\begin{proof}
The equality \eqref{eq:Higgsstrat} is clear by definition since each $(E,\phi)$ of bHN type $\tau$ lies in $\mathscr{H}_{\tau,[i,j]}^{r,d}(\Sigma,L)$ for some $(i,j)$. 
In order to show that \eqref{eq:Higgsstrat} is a stratification of $\mathscr{H}_{\tau}^{r,d}(\Sigma,L)$ in the sense of \thref{defn:stratification}, we need a suitable partial order on the indexing set. The decomposition given in \eqref{eq:Higgsstrat} is indexed by the union of the set $\{[i,j]_{\tau}  \ | \ i > j \}$ with the singleton set corresponding to collapsing all equivalence classes $[i,j]$ such that $i \leq j$ to one point. We denote its element by $[0,0]_{\tau}$ for simplicity. As seen in Section \ref{subsec:bHNgraded}, a bHN type $\tau = (d_1/r_1,\hdots, d_s/r_s)$ induces an ordering of equivalence classes $[i,j]_{\tau}$ given by $[i,j]_{\tau} < [i',j']_{\tau}$ if $d_i/r_i - d_j/r_j < d_{i'}/r_{i'} - d_{j'}/{r_{j'}}$. This provides a total order on the indexing set for the decomposition of \eqref{eq:Higgsstrat}, with $[i,j]_{\tau} < [0,0]_{\tau}$ whenever $i > j$. 

Next we need to show that given an index $[i,j]_{\tau}$, the union $$\bigsqcup_{[i',j'] > [i,j]} \mathscr{H}_{\tau,[i',j']}^{r,d}(\Sigma,L) $$ is closed in $\mathscr{H}_{\tau}^{r,d}(\Sigma,L)$. To see this, it suffices to observe that $(E,\phi) \in \bigsqcup_{[i',j']_{\tau} > [i,j]_{\tau}} \mathscr{H}_{\tau,[i',j']}^{r,d}(\Sigma,L)$ if and only if $\phi_{kl} = 0$ for all $(k,l) \leq (i,j)$, which is a closed condition. We note that even though $\phi_{kl}$ may not always well-defined, as we have seen in \thref{welldefined} it is well-defined if  $\phi_{k'l'} = 0$ for all $(k',l') < (k,l)$. Since $\phi_{s1}$ is always well-defined, the condition that $\phi_{kl} = 0$ for all $(k,l) \leq (i,j)$ is valid. Hence the union $$\bigsqcup_{[i',j']_{\tau} \leq [i,j]_{\tau}} \mathscr{H}_{\tau,[i',j']}^{r,d}(\Sigma,L)$$ is open in $\mathscr{H}_{\tau}^{r,d}(\Sigma,L)$, and we have that $\mathscr{H}_{\tau,[i,j]}(\Sigma,L)$ is closed inside this open substack.  

It follows from the above that the closure of each stratum $\mathscr{H}_{\tau,[i,j]}^{r,d}(\Sigma,L)$ is contained in the union $$\mathscr{H}_{\tau,[i,j]}^{r,d}(\Sigma,L) \sqcup \bigsqcup_{[i',j']_{\tau} > [i,j]_{\tau}}^{r,d} \mathscr{H}_{\tau,[i',j']},$$ as required for \eqref{eq:Higgsstrat} to satisfy the conditions of \thref{defn:stratification}. 

Combining the total order on the index set associated to a given bHN type $\tau$ with the partial order on the bHN types gives a partial order on the index set for the stratification of $\mathscr{H}_{r,d}(\Sigma,L)$ obtained by combining the stratification given in \eqref{eq:Higgsstrat} for each bHN stratum. 

To prove \eqref{eq:00}, we recall that by \thref{underlying0} \ref{part:equality} the HN and bHN types of $(E,\phi)$ are equal if and only if the corresponding filtrations coincide. But this is possible if and only if the smallest equivalence classs $[i,j]_{\tau}$ associated to $(E,\phi)$ satisfies $i \leq j$, i.e.\ if and only if $(E,\phi) \in \mathscr{H}_{\tau,[0,0]}^{r,d}(\Sigma,L)$. 
 \end{proof}

We conclude this section by endowing each of the strata of the refined bHN stratification (distinct from the semistable stratum) with a surjective algebraic map to a simpler moduli stack, using the map sending a Higgs bundle to its associated bHN graded. Note that restricting to the refined strata is necessary to ensure that the map sending a Higgs bundle to its bHN graded is continuous. 

Given a bHN type $\tau = (d_1/r_1,\hdots, d_s/r_s)$, fix an equivalence class $[i,j]_{\tau}$ such that $j \leq i$. We consider the associated refined bHN stratum $\mathscr{H}_{\tau,[i,j]}^{r,d}(\Sigma,L)$. If $[i,j]_{\tau} = [0,0]_{\tau}$, then $\operatorname{sgr}(E,\phi) = \operatorname{gr}(E,\phi)$ and there is a well-defined map \begin{equation} \operatorname{sgr}:  \mathscr{H}_{\tau,[0,0]}^{r,d}(\Sigma,L) = \mathscr{H}^{r,d}_{\tau}(\Sigma,L) \cap \mathscr{H}_{r,d}^{\tau}(\Sigma,L) \to \prod_{i=1}^{s} \mathscr{H}^{ss}_{r,d}(\Sigma,L) \cap \mathscr{H}^{r,d}_{ss}(\Sigma,L) \label{sgr00}
\end{equation} given by $(E,\phi) \mapsto ((E_1, \phi_1),\hdots, (E_s, \phi_s))$. Thus a Higgs bundle $(E,\phi)$ in this refined bHN stratum is determined by an $s$-tuple of semistable Higgs bundles with underlying semistable bundle together with some extension data indicating how the Higgs bundle can be reconsructed from its graded pieces. 

If $[i,j]_{\tau} \neq [0,0]_{\tau}$ and $(E,\phi) \in \mathscr{H}_{\tau,[i,j]}^{r,d}(\Sigma,L)$, then by definition $\operatorname{sgr}(E,\phi) = (\operatorname{gr} E, \phi_{\operatorname{min}})$ where $\phi_{\operatorname{min}} = \bigoplus_{(k,l) \in [i,j]_{\tau}} \phi_{kl}$. Let $n$ denote the cardinality of the set of pairs $(k,l)  \in [i,j]_{\tau}$. Then $\operatorname{sgr}(E,\phi)$ determines an $n$-tuple of holomorphic chains $((E_l,\phi_{kl}),(E_k \otimes L, 0))$ of rank $r_{kl}:= (r_l,r_k)$ and degree $d_{kl}:= (d_l,d_k+r_k \operatorname{deg} L)$ (see \cite{GarciaPrada2011}). Let $\mathscr{C}_{\tau,[i,j]}(\Sigma,L)$ denote the stack parametrising such objects. Note that it is an algebraic stack as it admits a representable morphism to $\prod_{(k,l) \in [i,j]_{\tau}} \mathscr{V}_{r_l,d_l}(\Sigma) \times \mathscr{V}_{r_k, d_k+r_k \operatorname{deg} L}(\Sigma)$ given by the forgetful map sending $((E_l,\phi_{kl}),(E_k \otimes L, 0))$ to $(E_l,E_k \otimes L )$. We let $\mathscr{C}^{ss}_{\tau,[i,j]}(\Sigma,L)$ denote the open substack given by the preimage of $\prod_{(k,l) \in [i,j]} \mathscr{V}_{r_l,d_l}^{ss}(\Sigma) \times \mathscr{V}_{r_k, d_k+r_k \operatorname{deg} L}^{ss}(\Sigma)$ under this map. Note that $E_k \otimes L$ is semistable if and only if $E_k$ is. Then there is a map  \begin{equation} \operatorname{sgr}: \mathscr{H}_{\tau,[i,j]}^{r,d}(\Sigma,L) \to  \prod_{
\substack{(k,l) \notin [i,j]_{\tau}\\
k > l}} \mathscr{V}_{r_k,d_k}^{ss}  \times \mathscr{C}_{\tau,[i,j]}^{ss}  \label{sgrij}
\end{equation} given by mapping $(E,\phi)$ to the $n$-tuple of holomorphic chains as described above, together with the summands of the HN graded of $E$ which do not appear in any of these holomorphic chains. Note that the definition of this map ensures that the following diagram commutes:
\begin{equation*} 
 \begin{tikzcd}
\mathscr{H}_{\tau,(i,j)}^{r,d}(\Sigma,L) \ar{r}{\operatorname{sgr}} \ar{d}{F} &  \prod_{
\substack{(k,l) \notin [i,j]_{\tau}\\
k > l}} \mathscr{V}_{r_k,d_k}^{ss}  \times \mathscr{C}_{\tau,[i,j]}^{ss}  \ar[d] \\
\mathscr{V}^{\tau}(\Sigma) \ar{r}{\operatorname{gr}} & \prod_{i=1}^s \mathscr{V}_{r_i,d_i}^{ss}(\Sigma),
\end{tikzcd} 
\end{equation*} 
where the vertical map on the right is given by forgetting the morphism $\phi_{ij}$.

The maps defined at \eqref{sgr00} and \eqref{sgrij} are the analogues for the refined bHN strata of the maps $\operatorname{gr}$ and $\operatorname{gr}$ defined for the HN strata of the stack of vector bundles and the HN strata of the stack of Higgs bundles. The existence of these maps imply that a Higgs bundle $(E,\phi) \in \mathscr{H}_{\tau,[i,j]}^{r,d}(\Sigma,L)$ is determined by an associated $s$-tuple, namely its bHN graded, together with some extension data indicating how the Higgs bundle can be reconstructed from its graded pieces.

\begin{rk}[Non-separatedness of the refined bHN strata]
Recall from \thref{HNnotsep} that the HN strata are not in general separated because a Higgs bundle can be deformed to its HN graded to which it need not be isomorphic. By contrast it is not clear whether it is possible to deform a Higgs bundle in a given refined bHN stratum to its associated bHN graded. Nevertheless it is possible to see that the refined bHN strata are not necessarily separated in a different way, by looking at the Higgs field scaling action. 

For example, let $r=3$ and fix a bHN type $\tau$ of length three. Suppose moreover that $[3,2]_{\tau} \leq [2,1]_{\tau}$ (otherwise we just switch the roles of the two equivalence classes). Consider a Higgs bundle \begin{equation} (E,\phi)  = \left( E_1 \oplus E_2 \oplus E_3, \begin{pmatrix} 
0 & 0 & 0 \\
\phi_{21} & 0 & 0 \\
0 & \phi_{32} & 0
\end{pmatrix} \right) \label{higgsbundle} \end{equation} where $\phi_{21}$ and $\phi_{32}$ are both non-zero, and $E_1 \oplus E_2 \oplus E_3 = \operatorname{gr} E$. 
Then $(E,\phi) \in \mathscr{H}_{\tau,[3,2]}^{3,d}(\Sigma,L)$, and moreover $(E,\phi)$ is fixed by the Higgs field scaling action (see \thref{linkwithcstar}). However, it is also possible to obtain the Higgs bundle \eqref{higgsbundle} with $\phi_{21}$ set to zero (or similarly with $\phi_{32}$ set to zero) as a limit of $(E,\phi)$ under the Higgs field scaling action, which is not isomorphic to $(E,\phi)$. This shows that the refined bHN stratum is not separated. 
\end{rk}

\section{Comparison of the stratifications} \label{sec:comparison} 

In this section we compare the three stratifications defined in Section \ref{sec:stratifications} above. In Section \ref{subsec:intersection} we compare the HN and bHN stratifications in arbitrary rank, while in Section \ref{subsec:compoffilt} we go one step further and compare the HN and bHN filtrations in three special cases: when $\operatorname{deg} L = 0$ and when $r=2$ or $3$. In Section \ref{subsec:comparisonof2} we compare the refined bHN stratification and the HN stratifications, specialising to the three cases considered in Section \ref{subsec:comparisonof2}.

\subsection{Comparison of the HN and bHN stratifications} \label{subsec:comparisonof2}

In this section we address the following question: which bHN strata can a given HN stratum can intersect and conversely which HN strata can a given bHN stratum intersect? \thref{underlying0} below establishes bounds on the maximal slope of the HN type of a Higgs bundle in terms of that of its bHN type, and vice-versa. As a consequence of this result we deduce in \thref{intersection} that a HN stratum can intersect only a finite number of bHN strata and vice-versa.

\begin{prop}[Relationship between bHN and HN types] \thlabel{underlying0} Suppose that $(E,\phi)$ is a Higgs bundle of HN type $\mu=(d_1/r_1,\hdots,d_s/r_s)$ and bHN type $\tau=(d_1'/r_1', \hdots, d_t'/r_t')$. Then we have: 
\begin{enumerate}[(i)]
\item $\mu \leq \tau$ with equality if and only if the  HN and bHN filtrations of $(E,\phi)$ coincide; \label{part:equality} 
\item if $r_{\operatorname{max}}$ denotes the maximum of the $r_i$ for $i=1,\hdots, s$, then $$\operatorname{min} \left( \frac{d}{r},\frac{d_1'}{r_1'} - r \operatorname{deg} L \right) \leq \frac{d_1}{r_1} \leq  \frac{d_1'}{r_1'} \leq \frac{d_1}{r_1}  + r_{\operatorname{max}} \operatorname{deg} L.$$ \label{part:inequalities} 
\end{enumerate} 
\end{prop}

It follows from \thref{underlying0} \ref{part:equality} that the HN polygon of a Higgs bundle $(E,\phi)$ lies below its bHN polygon. Moreover, \thref{underlying0} \ref{part:inequalities} provides upper and lower bounds on the possible maximal slopes for the HN (bHN respectively) type of a Higgs bundle with a given bHN (HN respectively) type.

\begin{rk}[Comparison with a result by Nitsure]
In \cite[Prop 3.2]{Nitsure1991}, Nitsure shows that if $(E,\phi)$ is a semistable Higgs bundle and $\tau = (d_1'/ r_1',\hdots, d_t'/ r_t')$ is its bHN type, then \begin{equation} d_1' / r_1' \leq \mu(E) + \frac{(r-1)^2}{r} \operatorname{deg} L. \label{Nsequation} 
\end{equation}  In short, there is a limit to how unstable the underlying bundle of a semistable Higgs bundle can be. \thref{underlying0} \ref{part:inequalities} can be viewed as a generalisation of this result to unstable Higgs bundles. That is, the right-most inequality of \ref{part:inequalities} shows that if $(E,\phi)$ is an unstable Higgs bundle then there is a limit to how unstable the underlying bundle can be, which depends only on the HN type of $(E,\phi)$. 
\end{rk} 
We now prove \thref{underlying0}. 

\begin{proof}[Proof of \thref{underlying0}]
To prove \ref{part:equality}, we first note that if $\mu = \tau$ then the successive quotients associated to both filtrations have the same rank. It follows then from the uniqueness of maximally destabilising subbundles that the filtrations coincide. It is clear that if the filtrations coincide then so do their associated types. The inequality $\mu \leq \tau$ follows from \cite[Thm 2]{Shatz1977} which states that if $F \subseteq E$ is a subbundle, then $F$ lies on or below the HN polygon of $E$.

Suppose now that $(E,\phi)$ has HN type $\mu = (d_1/r_1,\hdots, d_s/r_s)$. Let $0 = E^0 \subset E^1 \subset \cdots \subset E^s = E$ denote its HN filtration (set $E_i = E^{i-1}/E^i$ for each $i$), and let $0={E'}^0 \subseteq {E'}^1 \subseteq \cdots \subseteq {E'}^t=E$ denote its bHN filtration.
While each quotient $(E_i, \phi_{i})$ is semistable as a Higgs bundle, its underlying bundle may have a non-trivial HN filtration $$0 = {E_i^0} \subset {E_i^1} \subset \cdots  \subseteq {E_i^{t_i}} = E^{i}  / E^{i -1}.$$
Let $E^{i ,j}$ denote the preimage of ${E_i^j}$ under the quotient map $E^{i} \rightarrow E^{i} /E^{i-1}$, so that $E_i^j \cong  E^{i, j} / E^{i -1}$. Note that by definition we have the equalities $$E^{i} = E^{i+1, 0} = E^{i, t_{i}}.$$ Since $(E_i,\phi_i)$ is a semistable Higgs bundle for each $i$, we can apply \eqref{Nsequation} to obtain that \begin{equation}  \mu \left(E_i^j / E_i^{j-1} \right) < \mu(E_i^1) \leq \mu(E_i) + \frac{(r_i-1)^2}{r_i} \operatorname{deg}L  \label{previousinequality}
\end{equation} for each $i = 1,\hdots, s$ and $j = 1,\hdots, t_i$.  

The HN filtration of $E$ can be refined to include all of the subbundles $E^{i, j }$, giving a refined filtration $$0 = E^0 \subset E^{1,1} \subset \cdots \subset E^{1, s_1} \subset E^{2,0} \subset E^{2, 1} \subset \cdots \subset   E^{s, t_s} = E$$ of $E$. Although the slopes of these subbundles may no longer satisfy the decreasing condition, each successive quotient $$E^{i,j} / E^{i, j-1} \cong \left( E^{i,j}   / E^{i-1 } \right) / \left(E^{i, j-1} /  E^{i -1} \right) = {E_i^j} / {E_i^{j-1}}$$ is a semistable vector bundle (note that ${E^1}' = E^{1,1}$).

Suppose that $\phi({E'}^1) \subseteq E^{k,l} \otimes L$ but $\phi({E'}^1) \nsubseteq E^{k, l-1} \otimes L$ (if $l=0$ then $E^{k,l}=E^{k-1,t_{k}}$ and we replace $E^{k,l-1}$ by $E^{k-1,t_{k-1}-1}$).  Then the restriction of the Higgs field to ${E'}^1$ induces a non-zero map ${E'}^1 \rightarrow E^{k,l}/ E^{k, l-1} \otimes L$. Since both ${E'}^1$ and $E^{k,l}/ E^{k, l-1} \otimes L$ are semistable vector bundles, it follows that $$\mu({E'}^1) \leq \mu \left( E^{k,l}/ E^{k, l-1} \otimes L \right) = \mu \left( E^{k,l}/ E^{k ,l-1} \right) + \deg L = \mu \left( E^l_k / E_k^{l-1} \right) + \operatorname{deg} L.$$ Combining this inequality with \eqref{previousinequality}, we obtain that $$\mu({E'}^1) \leq  \frac{d_k}{r_k} +  \left( \frac{(r_k - 1)^2}{r_k} + 1 \right)    \deg L .$$ Now let $r_{\operatorname{max}}$ denote the maximum of the ranks $r_{i}$ for $i=1,\hdots, s$. Then $$
\frac{(r_k - 1)^2}{r_k} + 1 = r_k - 1 + \frac{1}{r_k}  \leq r_{\operatorname{max}} $$ since $r_k \geq 1$. 
Using the fact that $d_1/r_1$ is maximal amongst all of the slopes $d_{i}/r_{i}$ for $i=1,\hdots, s$, we obtain that \begin{equation} \mu({E'}^1) = \frac{d_1'}{r_1'} \leq \frac{d_1}{r_1} +  r_{\operatorname{max}} \operatorname{deg} L \label{2ndinequality}.
\end{equation}

Since ${E'}^1$ is the maximally destabilising subbundle of $E$, we have that $\mu({E'}^1) \leq \mu(E^1)$ and so from \eqref{2ndinequality} we obtain: $$\frac{d_1'}{d_2'} - r_{\operatorname{max}}  \operatorname{deg} L \leq \frac{d_1}{r_1} \leq \frac{d_1'}{r_1'} \leq  \frac{d_1}{r_1} +  r_{\operatorname{max}} \operatorname{deg} L.$$ Moreover since $r  \geq r_i$ for all $i$, in particular $r \geq r_{\operatorname{max}}$ and so we see we obtain the desired sequence of inequalities $$\frac{d_1'}{r_1'} - r  \operatorname{deg} L \leq \frac{d_1}{r_1} \leq  \frac{d_1'}{r_1'} \leq \frac{d_1}{r_1}  +  r_{\operatorname{max}}  \operatorname{deg} L.$$  \end{proof}

\begin{cor}[Intersection of the bHN and HN stratifications] \thlabel{intersection}
If $\mu$ is a HN type, then there is only a finite number of possible bHN types $\tau$ for a Higgs bundle $(E,\phi)$ of HN type $\mu$. Similarly, if $\tau$ is a bHN type, then there is only a finite number of HNN types $\mu$ for a Higgs bundles of bHN type $\tau$. 
\end{cor}

\begin{proof}
Let $\mu$ denote a HN type. If $(E,\phi)$ has HN type $\mu$ then by \thref{underlying0} \ref{part:inequalities}, the steepest slope arising in its bHN type $\tau$ is bounded above by $d_1/r_1 + r_{\operatorname{max}} \operatorname{deg} L$, which depends only on $\mu$. But the decreasing condition on the slopes of a HN type ensures that bounding the slope of the maximally destabilising subbundle bounds the number of allowable HN types of $E$. Thus there are only a finite number of possible bHN types (with a fixed rank and degree) for $(E,\phi)$. 

If $\tau$ denotes a bHN type and $(E,\phi)$ has bHN type $\tau$, then by \thref{underlying0} \ref{part:equality} we have that $\mu \leq \tau$ where $\mu$ is its HN type. There are therefore only finitely many possibilities for $\mu$.  
\end{proof} 

\subsection{Comparison of the HN and bHN filtrations} \label{subsec:compoffilt} 

In this section we address the following question: how do the HN and bHN filtration of a Higgs bundle compare? This is a difficult question in general, and we restrict to three special cases in which a complete answer can be given: when $\operatorname{deg}L =0$ (Section \ref{subsubsec:degL=0}), when $r=2$ (Section \ref{subsubsec:r=2}) and when $r=3$ (Section \ref{subsubsec:r=3}). In the case where $r=2$ and $r=3$ we also obtain explicit criteria for determining when a Higgs bundle with unstable underlying bundle is semistable.

\subsubsection{When $\operatorname{deg} L =0$} \label{subsubsec:degL=0}

If $\operatorname{deg} L =0$, which is the case if $\Sigma$ is an elliptic curve and $L$ is the canonical line bundle on $\Sigma$, then it follows from \thref{underlying0} that the maximally destabilising bundles of the HN and bHN filtrations of a given Higgs bundle always coincide. In fact in this case we can prove a stronger result directly, without using \thref{underlying0}. That is, we can prove directly that all of the terms in the HN and bHN filtrations coincide.

\begin{prop}[When $\operatorname{deg} L = 0$]   \thlabel{degL0}
If $\operatorname{deg}L = 0$, and $(E,\phi)$ is an $L$-twisted Higgs bundle on $\Sigma$ then its HN and bHN filtrations coincide, so that $\mathscr{H}_{\tau}^{r,d}(\Sigma,L) = \mathscr{H}^{\tau}_{r,d}(\Sigma,L)$ for any HN type $\tau$. 
\end{prop}

\begin{rk}[Comparison with a result by Franco]
In \cite{Franco2014}, Franco studies Higgs bundles in the classical sense (namely where $L$ is given by the cotangent bundle) on an elliptic curve (hence in this case $\operatorname{deg} L =0$), and proves that $(E,\phi)$ is semistable if and only if $E$ is semistable (see \cite[Prop 4.1]{Franco2014}). \thref{degL0} can therefore be viewed as a generalisation of this result to unstable Higgs bundles. 
\end{rk}

\begin{proof} [Proof of \thref{degL0}]
Suppose that $\operatorname{deg}L =0$ and that $(E,\phi)$ has HN filtration $0 = E^0 \subseteq E^1 \subseteq \cdots \subseteq E^s$ and bHN filtration $0={E'}^0 \subseteq {E'}^1 \subseteq \cdots \subseteq {E'}^t$. We show by induction that $E^i = {E'}^i$ for each $i =1 ,\hdots, s$ (which implies that $s = t$). 

The case where $i = 0$ is trivially true, so suppose that $E^i = {E'}^i$ for all $i \leq k$. We wish to show that $E^{k+1} = {E'}^{k+1}$. Since $E^k= {E'}^k$ is $\phi$-invariant, there is an induced map $\overline{\phi}: E/{E'}^k \to E / {E'}^k \otimes L$. Let $l$ denote the minimal $i \in \{1,\hdots,t\}$ such that $\overline{\phi}({E'}^{k+1}/{E'}^k) \subseteq {E'}^l/{E'}^k \otimes L$. Then $\overline{\phi}$ induces a non-zero map $${E'}^{k+1}/{E'}^k \to ({E'}^l/{E'}^k)/({E'}^{l-1}/{E'}^k) \otimes L = {E'}^l/{E'}^{l-1} \otimes L.$$ Since both vector bundles are semistable, we must have that $$\mu \left({E'}^{k+1}/{E'}^k \right) \leq \mu \left({E'}^l/ {E'}^{l-1} \otimes L \right) = \mu \left({E'}^l/{E'}^{l-1} \right),$$ where the equality follows from the assumption that $\operatorname{deg} L = 0$. Hence it must be that $l=k$, and so $\overline{\phi}$ preserves ${E'}^{k+1}/{E'}^k$. Therefore $\phi$ preserves ${E'}^{k+1}$ and so $E^{k+1} = {E'}^{k+1}$. 
\end{proof} 

\subsubsection{When $r=2$} \label{subsubsec:r=2}

In the rank two case (we no longer assume that $\operatorname{deg} L=0$), a stronger result than \thref{underlying0} is also true. In particular if a Higgs bundle of rank two is unstable then its bHN and HN types coincide.  

\begin{prop}[Relationship between HN and bHN types in rank two] \thlabel{rk2strats} 
Let $(E,\phi)$ be a Higgs bundle of rank $2$.  Then we have: \begin{enumerate}[(i)]
\item if $(E,\phi)$ is unstable, then its HN and bHN filtrations (and therefore types) coincide; \label{firstcond}
\item if $(E,\phi)$ is semistable, then its bHN type $\tau = (d_1', d_2')$ satisfies $d_1' \leq (d+1) \operatorname{deg} L / 2$; \label{thirdcond}
\end{enumerate} 
\end{prop}

Before proving \thref{rk2strats}, we deduce the following \thref{sscritrk2} which provides a criterion for semistability of a Higgs bundle of rank $2$ with unstable underlying bundle.

\begin{cor}[Criterion for semistability of a rank $2$ Higgs bundle with unstable underlying bundle] \thlabel{sscritrk2} 
Let $\tau$ denote a bHN type of rank $2$ and suppose that $(E,\phi)$ has bHN type $\tau$. Then $(E,\phi)$ is semistable if and only if $\phi_{21} \neq 0$, where $\phi_{21}: {E^1} \to E/{E^1} \otimes L$ denotes the composition of the restriction of $\phi$ to the maximally destabilising subbundle ${E^1}$ of $E$ with the quotient map $E \otimes L \to E/{E^1} \otimes L$. 
\end{cor}

\begin{proof}
Suppose that $(E,\phi)$ has bHN type $\tau$ and $\phi_{21} \neq 0$. By \thref{rk2strats} \ref{firstcond}, we know that if $(E,\phi)$ is unstable, then its HN and bHN filtrations coincide. But since $\phi_{21} \neq 0$, the Higgs field does not preserve ${E^1}$, which is a contradiction. So $(E,\phi)$ is semistable. 
Conversely if $(E,\phi)$ has bHN type $\tau$ and is semistable, then $\phi$ cannot preserve ${E^1}$ and so $\phi_{21}$ cannot be zero. 
\end{proof} 

\begin{proof}[Proof of \thref{rk2strats}]
For \ref{firstcond}, suppose that $(E,\phi)$ is unstable of HN type $\mu = (d_1,d_2)$ and bHN type $\tau = (d_1',d_2')$. Let $0 = E^0 \subseteq E^1 \subseteq E^2 = E$ denote its HN filtration and $0 = {E'}^0 \subseteq {E'}^1 \subseteq {E'}^2 = E$ its bHN filtration. Then ${d_1} \leq d_1'$ since otherwise ${E^1}$ would be a destabilising subbundle of $E$ of larger degree than ${E'}^1$, contradicting the bHN type of $(E,\phi)$. 

The composition of the inclusion of ${E'}^1$ into $E$ with the quotient map $E \to E/ {E^1}$ produces a map ${E'}^1 \to E/{E^1}$. Since both ${E'}^1$ and $E^1$ are line bundles, they are stable. If the map is non-zero, then we must have that $d_1' \leq d_2$, which is not possible since $d_2 < d_1 \leq d_1'$. Thus the map is zero and so ${E'}^1 \subseteq E^1$. As both have rank one, they must be equal. 
It follows that the two filtrations coincide. Part \ref{thirdcond} follows immediately from \eqref{Nsequation}. 
\end{proof}

\subsubsection{When $r=3$} \label{subsubsec:r=3}

In rank three it is again possible to give a more explicit comparison of the HN and bHN filtrations. To state the comparison it is useful to define the \emph{splitting type} of a HN type. 

\begin{defn}[Splitting type of a HN type]
Let $\tau=(d_1/r_1,\hdots, d_s/r_s)$ denote a HN type. The \emph{splitting type} of $\tau$ is the $s$-vector $(r_1,r_2,\hdots, r_s)$.  
\end{defn}

\begin{prop}[When $r =3$] \thlabel{rk3strats}
Let $(E,\phi)$ denote a rank three Higgs bundle. We let $E^i$ and ${E'}^i$ denote the terms in the HN and bHN filtrations respectively of $(E,\phi)$. For each $i$ let $\phi_i': {E'}^i \to E/{E'}^i \otimes L$ denote the composition of the restriction of $\phi$ to ${E'}^i$ with the quotient map $\pi_i': E \otimes L \to E/{E'}^i \otimes L$. Given a subsheaf $\mathcal{F} \subseteq G$ of a vector bundle $G$, let $\overline{\mathcal{F}}$ denote the subbundle of $G$ generically generated by $\mathcal{F}$. We define the subbundles $F^i$ and $G^i$ of $E$ and $E \otimes L$ respectively as follows: \begin{itemize}
\item $F^i: = \overline{\operatorname{ker}(\phi_i')} \subseteq E$;
\item $G^i: = \overline{{\pi_i'}^{-1}(\operatorname{im}(\phi_i'))}   \subseteq E \otimes L $.  
\end{itemize} 

Then the HN filtration of $(E,\phi)$ can be recovered from the bHN filtration and from the Higgs field $\phi$, as per the following table.

\begin{center}
\begin{tabular}{c c c  }
Splitting type of bHN type & Splitting type of HN type & HN filtration in terms of bHN filtration  \\
\hline
\multirow{3}{*}{$(1,1,1)$} & $(1,1,1)$ & Both filtrations coincide \\
 & $(2,1)$ & $E^1 \otimes L = G^1 \supseteq {E'}^1 \otimes L$ \\
 & $(1,2)$ & $E^1 = F^2  \subseteq {E'}^2$  \\
\hline
\multirow{3}{*}{$(2,1)$} & $(1,1,1)$ & n/a \\ 
 & $(2,1)$ & Both filtrations coincide\\
 & $(1,2)$ & $E^1 = F^1$ \\
\hline 
\multirow{3}{*}{$(1,2)$} & $(1,1,1)$ & n/a \\
& $(2,1)$ & $E^1 \otimes L = G^1$ \\
 & $(1,2)$ & Both filtrations coincide. 
\end{tabular} 
\end{center} 
The entries `n/a' indicate that the corresponding combination of splitting types for the bHN and HN type of a rank three Higgs bundle cannot occur. 
\end{prop}

Before proving \thref{rk3strats}, we deduce the following \thref{sscriterion} which provides necessary and sufficient conditions for semistability of a Higgs bundle of rank three with unstable underlying bundle. While by definition checking semistability of a Higgs bundle requires checking a slope condition for every Higgs field invariant subbundle, by \thref{sscriterion} below we see that it is sufficient to check the slopes of at most two subbundles, determined in terms of the bHN filtration of the Higgs bundle and its Higgs field.  

\begin{cor}[Criterion for semistability of a rank three Higgs bundle with unstable underlying bundle] \thlabel{sscriterion} 
Let $\tau$ denote a bHN type of rank $3$ and suppose that $(E,\phi)$ has bHN type $\tau$. Let ${E'}^i$ denote the terms in the bHN filtration of $(E,\phi)$.
Using the isomorphism $E \cong E \otimes L \otimes L^{-1}$, we identify $F^2 \otimes L^{-1} \subseteq E \otimes L \otimes L^{-1}$ as a subbundle of $F$. 
Then we obtain the following necessary and sufficient conditions for semistability of $(E,\phi)$: 

\begin{center}
\begin{tabular}{c p{10.3cm}} 
Splitting type of bHN type $\tau$ & $(E,\phi)$ is semistable if and only if \\
\hline
 \multirow{2}{*}{$(1,1,1)$} & $\operatorname{rk} F^2 =1$ and if $F^2$ is $\phi$-invariant then $\operatorname{deg} F^2 \leq 2d/3$, and \\
 & $\operatorname{rk} G^1 = 2$ and  if $G^1 \otimes L^{-1}$ is $\phi$-invariant then $ \operatorname{deg} G^1 \leq 2 d/3 +  \operatorname{deg}L$  \\
 \hline
$(2,1)$ & $\operatorname{rk} F^1 = 1$ and if $F^1$ is $\phi$-invariant then $\operatorname{deg} F^1 \leq d/3$ \\
\hline
$(1,2)$ & $\operatorname{rk} G^1 = 2$ and if $G^1 \otimes L^{-1}$ is $\phi$-invariant then $\operatorname{deg} G^1 \leq 2 d/3 +  \operatorname{deg} L$. 
\end{tabular} 
\end{center}
\end{cor}

\begin{rk}[Link with the Bialynicki-Birula stratification] 
It follows from \thref{sscriterion} above that the slopes of the bundles $F^i$ and $G^i$ defined in \thref{rk3strats} play an important role in detecting the semistability of a Higgs bundle of rank three. In fact these same two bundles and their slopes play an important role in understanding the limit of a semistable Higgs bundle under the Higgs field scaling action, as demonstrated by Gothen and Zu\~niga Rojas in \cite{Gothen2015}. 
\end{rk}

\begin{proof}[Proof of \thref{sscriterion}] 
Suppose that $(E,\phi)$ has bHN type $\tau$ of splitting type $(1,1,1)$. By \thref{rk3strats}, the only possibly destabilising subbundles of $E$ are $F^2$ and $G^1 \otimes L^{-1}$. Suppose that the conditions in the first row of the above table are satisfied. If neither of the two subbundles is $\phi$-invariant, then $(E,\phi)$ is clearly semistable. Suppose then that $F^2$ is $\phi$-invariant.  Then $\mu(F^2) = \operatorname{deg} F^2 \leq 2d/3 = \mu(E)$. Therefore $F^2$ does not destabilise $(E,\phi)$ as $\mu(F^2) = \operatorname{deg} F^2$. Now suppose that $G^1$ is $\phi$-invariant. Then $\mu(G^1 \otimes L^{-1}) = \operatorname{deg}G^1/2 - \operatorname{deg} L \leq 2 d/3$, and therefore $G^1 \otimes L^{-1}$ is not destabilising either. Thus $(E,\phi)$ is semistable. 

Conversely, suppose that $(E,\phi)$ is semistable. Since $F^2$ is the kernel of $\phi_1': {E'}^2 \to E/{E'}^2  \otimes L$, it is a subbundle of the rank two subbundle ${E'}^2 $. If $F^2$ has rank two then it coincides with ${E'}^2 $,and so $F^2$ is $\phi$-invariant and destabilising which contradicts semsitability of $(E,\phi)$. We cannot have that $\operatorname{rk} F^2 = 0$ as $\phi_1'$ is a map from a rank two bundle to a rank one bundle. Thus $\operatorname{rk} F^2 = 1$. Since $(E,\phi)$ is semistable, we must have that $\mu(F^2) = \operatorname{deg} F^2 \leq \mu(E)  = 2d/3$. Now $G^1$ is a subbundle of $E \otimes L$ containing ${E'}^1$. If it has rank one, then ${E'}^1 = G^1$ and so $G^1$ is a destabilising subbundle preserved by $\phi$, which contradicts semistability of $(E,\phi)$. So $G^1$ has rank greater than or equal to two. It cannot have rank three as $\pi_1': E \otimes L \to E/{E'}^1 \otimes L$ is surjective and  $\operatorname{im} \phi_1' \subset E/{E'}^1 \otimes L$ is a strict subbundle. Therefore $\operatorname{rk} G^1 = 2$. If $G^1 \otimes L^{-1}$ is $\phi$-invariant, then since $(E,\phi)$ is semistable we must have that $\mu(G^1 \otimes L^{-1}) = \operatorname{deg} G^1/2 - \operatorname{deg} L \leq 2d/3 = \mu(E)$ which is equivalent to the inequality $\operatorname{deg} G^1 \leq 2d/3 + \operatorname{deg} L$. 

The proofs in the cases where $\tau$ has splitting type $(2,1)$ or $(1,2)$ are completely analogous, with the difference that in both cases there is only one possibility for the potentially destabilising subbundle. 
\end{proof} 

We illustrate \thref{rk3strats} and \thref{sscriterion} below in an example.

\begin{example} \thlabel{p1ex} 
Consider the vector bundle $E:= \mathcal{O}(2) \oplus \mathcal{O}(1) \oplus \mathcal{O}(-1)$ on $\mathbb{P}^1$. Fix a line bundle $\mathcal{O}(t)$ on $\mathbb{P}^1$ of degree $t \geq 3$. It has HN filtration $0 \subseteq \mathcal{O}(2) \subseteq \mathcal{O}(2) \oplus \mathcal{O}(1) \subseteq E$, and therefore HN type $(2,1,-1)$, of splitting type $(1,1,1)$. 
Now let $\phi_{31} \in H^0(\Sigma,\mathcal{O}(-3+t))$ denote a non-zero section, so that $\phi_{31}$ determines a non-zero map $\mathcal{O}(2) \to \mathcal{O}(-1) \otimes \mathcal{O}(t)$. Then $\phi_{31}$ defines a Higgs field $\phi: E \to E \otimes \mathcal{O}(t)$ given by $$\phi : = \begin{pmatrix} 0 & 0 & 0 \\
0 & 0 & 0 \\
\phi_{31} & 0 & 0 
\end{pmatrix};$$ here $\phi$ is written in matrix form with respect to the decomposition of $E \cong \operatorname{gr} E$ given by its HN decomposition. 
Then the map $\phi_2': \mathcal{O}(2) \oplus \mathcal{O}(1) \to \mathcal{O}(-1) \otimes \mathcal{O}(t)$ has kernel $\mathcal{O}(1)$, and so $F^2 = \mathcal{O}(1)$, which is $\phi$-invariant. The map $\phi_1': \mathcal{O}(2) \to (\mathcal{O}(1) \oplus \mathcal{O}(-1)) \otimes \mathcal{O}(t)$ has image $\mathcal{O}(-1) \otimes \mathcal{O}(t)$, which has preimage $\mathcal{O}(-1) \otimes \mathcal{O}(t)$ under the quotient map $\pi_i'$. Therefore $G^1 = (\mathcal{O}(2) \oplus \mathcal{O}(-1)) \otimes \mathcal{O}(t)$. Then $G^1 \otimes L^{-1} = \mathcal{O}(2) \oplus \mathcal{O}(-1)$ is also $\phi$-invariant. Hence by \thref{sscriterion}, we have that $(E,\phi)$ is semistable if and only if $\operatorname{deg} F^2 \leq d/3$ and $\operatorname{deg}G^1 \leq 2d/3 +  \operatorname{deg} L$. While $\operatorname{deg} G^1 = 1+ t  \leq 4/3 + t$, we have that $\operatorname{deg} F^2 = \operatorname{deg} \mathcal{O}(1) = 1 > 2/3$. Therefore $(E,\phi)$ is not semistable. Moreover since $G^1$ is not destabilising, by \thref{rk3strats} we see that $(E,\phi)$ must have HN filtration $0 = E^0 \subseteq F^2 \subseteq E$ and therefore HN type $\mu = (1,1/2,1/2)$ which has splitting type $(1,2)$. 
\end{example}

We now prove \thref{rk3strats} for which we proceed case by case. The two key results used throughout are that if $E \to F$ is a non-zero map of semistable vector bundles then $\mu(E) \leq \mu(F)$, and that all line bundles are semistable.  

\begin{proof}[Proof of \thref{rk3strats}] Let $\mu$ and $\tau$ denote the HN and bHN types of $(E,\phi)$ respectively. We start with the case where $\tau$ has splitting type $(1,1,1)$. To this end, suppose that $\tau = (d_1', d_2', d_3')$. Note that $d_1' + d_2' + d_3' = d$. As our first subcase, we suppose that $\mu$ has splitting type $(1,1,1)$, so that $\mu = (d_1,d_2,d_3)$. Note that we also have $d_1 + d_2 + d_3 = d$. The composition ${E'}^1 \subseteq E \to E/ {E^2}$ is a map of semistable vector bundles (a rank one vector bundle is always semistable). Therefore if the map is non-zero, then we must have that $\mu({E'}^1) = d_1' \leq \mu(E/{E^2})  =d_3$, which implies that $d_i' < d_i$ for each $i$ since $d_1 > d_2 > d_3 $ and $d_1' > d_2' > d_3'$. So then  $$d = d_1' + d_2' + d_3' < d = d_1 + d_2 + d_3$$ which is a contradiction. Therefore the composition is actually the zero map and so ${E'}^1 \subseteq {E^2}$. We now look at the composition ${E'}^1 \hookrightarrow E^2 \to E^2/E^1$ which is again a map of semistable vector bundles. If non-zero, then there is an inequality of slopes $d_1' \leq d_2$, which implies that $d_1' < d_1$ since $d_2 < d_1$. This is a contradiction as ${E'}^1$ is the maximally destabilising subbundle of $E$, so $E^1$ cannot have strictly larger degree. Therefore the composition must be zero and so ${E'}^1 \subseteq E^1$. As both are line bundles and hence of the same rank, they must be equal. We now show that ${E^2} = {E'}^2$. To do so, we consider the composition ${E'}^2/{E'}^1  = {E'}^2/E^1  \hookrightarrow E/E^1 \to E/E^2$, which if non-zero forces the inequality $d_2' \leq d_3$. But then $d_2' + d_3' < d_2 + d_3$ since $d_3' < d_2'$ and $d_3 < d_2$, and so we cannot have $d_1' + d_2' + d_3' = d = d_1 + d_2 + d_3$ (we have already established that $d_1' = d_1$). Thus the composition must be the zero map and so ${E'}^2 \subseteq E^2$. As both have the same rank, they must be equal. 
Thus the HN and bHN filtrations coincide. 

For the second subcase, suppose that $\mu$ has splitting type $(2,1)$, and hence that $\mu= (d_1/2,d_1/2,d_2)$. Note that $d_1 + d_2 = d$. If the composition ${E'}^1 \hookrightarrow E \to E/E^1$ is non-zero, then there is an inequality $d_1' \leq d_2$ as it is a map of semistable bundles. But $d_1 /2 > d_2$ and so $d_1/2 > d_1'$ which contradicts the fact that ${E'}^1$ is the maximally destabilising subbundle of $E$. Thus the composition is zero and so ${E'}^1 \subseteq E^1$. Now since $E^1$ is preserved by $\phi$, we must have that $\phi({E'}^1) \subseteq E^1/{E'}^1 \otimes L$, and therefore $G^1 = \overline{{\pi_1'}^{-1}(\operatorname{im} \phi_1'))} \subseteq E^1 \otimes L$. To show that they are equal, it suffices to show that $G^1$ has rank two. Now $G^1$ has at least rank one since it contains ${E'}^1 \otimes L$, yet it cannot coincide with ${E'}^1 \otimes L$ as otherwise its image under $\pi_1'$ would be trivial which would contradict the fact that ${E'}^1$ is not preserved by $\phi$. Therefore $G^1$ must have rank two.

For the third subcase we suppose that $\mu$ has splitting type $(1,2)$ so that $\mu = (d_1, d_2/2,d_2/2)$. Note that $d_1 + d_2 = 0$. Consider the composition $E^1 \hookrightarrow E \to E/{E'}^2$, which is a map of semistable vector bundles. If non-zero we must have that $d_1 \leq d_3'$. But then since $d_1 + d_2 < d_1 + 2 d_1 = 3 d_1$ (using the fact that $d_1 > d_2/2$), we obtain the following contradiction: $$d = d_1 + d_2 < 3 d_1 \leq 3 d_3' < d_1' + d_2' + d_3' = d.$$ Thus ${E^1} \subseteq {E'}^2$. We wish to show that ${E^1} = F^2 : = \overline{\operatorname{ker} \phi_2'}$. The first step is to note that $\phi_2'$ is not the zero map. Indeed, if it is then ${E'}^2$ is preserved by the Higgs field, which leads to ${E'}^2/E^1$ destabilising $E/E^1$ by a straightforward slope calculation.  
This contradicts the fact that $\mu$ has splitting type $(1,2)$ and so we can conclude that $\phi_2'$ is not the zero map. It follows that $F^2$ has rank one. Now since $E^1 \subseteq {E'}^2$ is preserved by the Higgs field, it must be contained in $\operatorname{ker} \phi_2'$ and therefore in $F^2$. Both are line bundles and therefore they are equal.

We now turn to the case where $\tau$ has splitting type $(2,1)$. Thus we suppose that $\tau = (d_1'/2, d_1'/2,d_2')$. Note that $d_1' + d_2' = d$. For the first subcase, suppose that $\mu$ has splitting type $(1,1,1)$, so that $\mu = (d_1,d_2,d_3)$. Note that $d_1 + d_2 + d_3 = d$. If the composition $E^1 \hookrightarrow E \to E/{E'}^1$ is non-zero, then $d_1 \leq d_2'$ as the composition is a map of semistable bundles. But then $d_1 + d_2 + d_3 = d < d = d_1' + d_2'$ which is a contradiction. Hence $E^1 \subseteq {E'}^1$. Since ${E'}^1$ is semistable, we must have that $d_1 \leq d_1'/2$. Now if we assume that the composition ${E'}^1 \hookrightarrow E \to E/ E^2$ is non-zero, we must have that $d_1'/2 \leq d_3$. But then $d_1'/2 < d_1$ since $d_3  < d_1$ and so we have a contradiction. Thus ${E'}^1 \subseteq E^2$. As both have the same rank, they must be equal.  
But then $E^1 \subseteq E^2 = {E'}^1$ is a destabilising subbundle, which contradicts the fact that ${E'}^1$ is semistable. Therefore $\mu$ cannot have splitting type $(1,1,1)$. 

Fo the second subcase, suppose that $\mu$ has splitting type $(2,1)$, so that $\mu = (d_1/2,d_1/2,d_2)$. If the composition ${E'}^1 \hookrightarrow E \to E/E^1$ is non-zero, then there is an inequality $d_1'/2 \leq d_2$ as the composition is a map of semistable bundles. But then $d_1/2 > d_2 > d_1'/2$ which contradicts the fact that ${E'}^1$ is maximally destabilising for $E$. Hence the composition is zero and so ${E'}^1 \subseteq E^1$. 
Since both have the same rank, they are equal. 

For the third subcase, suppose that $\mu$ has splitting type $(1,2)$ so that $\mu = (d_1, d_2/2, d_2/2)$. Note that $d_1 + d_2 = d$. If the composition $E^1 \hookrightarrow E \to E/{E'}^1$ is non-zero, then $d_1 \leq d_2'$ as the composition is a map of semistable bundles. But then we obtain the contradiction $$ d= d_1 + d_2 < 3 d_1 \leq 3 d_2' < d_1' + d_2' = d,$$ where the second and third inequalities follow from the fact that $d_1 > d_2/2$ and that $d_1'/2 > d_2'$ respectively. Therefore the composition is zero and so ${E^1} \subseteq {E'}^1$. Since $E^1$ must be preserved by the Higgs field, we must in particular have that $\phi(E^1) \subseteq {E'}^1 \otimes L$ so that $E^1 \subseteq  \overline{\operatorname{ker} \phi_1'} =: F^1$. But since we are assuming that $\mu$ has splitting type $(1,2)$, we know that $\phi$ does not preserve ${E^1}'$ and therefore $\operatorname{ker} \phi_1'$ is non-trivial. Since $E/{E^1}' \otimes L$ has rank one and ${E'}^1$ has rank two, it follows that $F^1$ must have rank one. Therefore the inclusion ${E^1} \subseteq  F^1$ is an equality.    

Finally we consider the case where $\tau$ has splitting type $(1,2)$ and hence we suppose that $\tau = (d_1', d_2'/2, d_2'/2)$. Note that $d_1' + d_2' = d$. For the first subcase we suppose that $\mu$ has splitting type $(1,1,1)$ so that $\mu = (d_1, d_2, d_3)$. Note that $d_1 + d_2 + d_3=d$. Our aim is to show that $E^1 = {E'}^1$, and then to conclude that $E^2/E^1$ destabilises $E/{E'}^1$, which will give us a contradiction. To do so, we consider the composition ${E'}^1 \hookrightarrow E \to E/E^2$, which is a map of semistable bundles. If non-zero, then we have that $d_1' \leq d_3$, which implies that $d_1' < d_1$. This contradicts the fact that ${E'}^1$ is maximally destabilising, so the composition must be zero. Hence ${E'}^1 \subseteq E^2$. If the composition ${E'}^1 \hookrightarrow E^2 \to E^2 / E^1$ is non-zero, then we have that $d_1' \leq d_2$ as the composition is a map of semistable bundles. This leads to the same contradiction as before, namely that $d_1' < d_1$. Hence we must have that ${E'}^1 \subseteq E^1$. As both are line bundles, they must be equal. 
But then $E^2/E^1 \subseteq E/E^1 = E/{E'}^1$ is a destabilising subbundle, which contradicts the semistability of $E/{E'}^1$. 

For the second subcase, suppose that $\mu$ has splitting type $(2,1)$ so that $\mu = (d_1/2, d_1/2,d_2)$. If the composition ${E'}^1  \hookrightarrow E \to E/E^1$ is non-zero, then we have that $d_1' \leq d_2$. But $d_2 < d_1/2$ so we obtain that $d_1' < d_1/2$ which contradicts the fact that ${E'}^1$ is maximally destabilising. Therefore the composition is zero and so ${E'}^1 \subseteq E^1$. Now since we have assumed that $\mu$ has splitting type $(2,1)$, we know that ${E'}^1$ is not preserved by $\phi$, while $E^1$ is. Therefore we must have that $\operatorname{im} \phi_1'$ is contained in $ E^1 / {E'}^1 \otimes L$. Thus ${E'}^1 \otimes L   \subseteq { \pi_1'}^{-1}(\operatorname{im} \phi_1') \subseteq E^1 \otimes L $. Hence $E^1 \otimes L $ contains $G^1 : = \overline{{ \pi_1'}^{-1}(\operatorname{im} \phi_1')}$ as a subbundle. The latter has rank two (since ${E'}^1 \otimes L$ can only be strictly contained in it given that $\phi_1'$ is non-zero), and therefore the inclusion must be an equality. 

For the third and final subcase we suppose that $\mu$ has splitting type $(1,2)$, so that $\mu  = (d_1, d_2/2, d_2/2)$. If the composition $E^1 \subseteq E \to E/{E'}^1$ is non-zero, then we have that $d_1 \leq d_2'/2$ since the composition is a map of semistable bundles. The fact that $d_1' > d_2'/2$ yields the contradiction $$d = d_1' + d_2' > d_2'/2 + d_2' = 3 (d_2'/2) \geq 3 d_1 > 3d/3 = d,$$ where the last inequality follows from observing that $E^1$ is destabilising so that $d_1 > d/3$. Therefore the composition is zero and so $E^1 \subseteq {E'}^1$. As both are line bundles, they must be equal. \end{proof}

\begin{rk}[Relationship between the HN and bHN filtrations and type in higher rank] 
The fact that in rank two and three the HN filtration can be computed from the bHN filtration and the Higgs field is quite special, and we do not expect it to be true in general. Indeed in our proofs of \thref{rk2strats} and \thref{rk3strats} we repeatedly appealed to the fact that line bundle are automatically semistable, but line bundles no longer appear in higher rank and the argument no longer carries over. 
\end{rk}

\subsection{Comparison of the HN and refined bHN stratifications}

 \label{subsec:intersection} 
In this section we comment on the relationship between the refined bHN stratification and the HN stratification, guided by the following question: does the refined bHN stratum in which a Higgs bundle lies determine the HN type of a Higgs bundle? We do so in the three special cases considered in Section \ref{subsec:compoffilt}, namely when $\operatorname{deg}L =0$, when $r=2$, and when $r=3$. This is because these are the only cases in which we can completely describe the relationship between the bHN and HN stratifications, a necessary first step towards describing the relationship between the refined bHN and HN stratifications. 

\subsubsection{When $\operatorname{deg} L = 0$} 

\begin{prop}[The three stratifications when $\operatorname{deg} L =0$] 
If $\operatorname{deg} L =0$, then the HN, bHN and refined bHN stratifications of $\mathscr{H}_{r,d}(\Sigma,L)$ all coincide. 
\end{prop} 

\begin{proof} 
By \thref{degL0}, when $\operatorname{deg} L =0$ the bHN and HN filtrations of $(E,\phi)$ coincide. Therefore given a bHN type $\tau$, the refined bHN strata $\mathscr{H}_{\tau,[i,j]}(\Sigma,L)$ are all empty except for $\mathscr{H}_{\tau,[0,0]}(\Sigma,L)$. Thus the refinement of the bHN stratification is actually trivial, and so all three stratifications coincide.  
\end{proof}

\subsubsection{When $r=2$} 

\begin{prop}[Refined bHN stratification when $r=2$]
The refined bHN stratification of $\mathscr{H}_{2,d}(\Sigma,L)$ is the intersection of the bHN stratification with the HN stratification. 
\end{prop} 

\begin{proof} 
 The refined bHN stratification of $\mathscr{H}_{\tau}^{2,d}(\Sigma,L)$ is given by $$\mathscr{H}_{\tau}^{2,d}(\Sigma,L) = \mathscr{H}_{\tau,[2,1]}^{2,d}(\Sigma,L) \sqcup \mathscr{H}_{\tau,[0,0]}^{2,d}(\Sigma,L).$$ 
By definition, the Higgs bundle $(E,\phi)$ lies in $\mathscr{H}_{\tau,[2,1]}^{2,d}(\Sigma,L)$ if and only if the map $\phi_{21}: E'_1 \to E'_2 \otimes L$ is non-zero, where $E_1' \oplus E_2'$ is the HN graded of $E$. Thus $ {E'}^1 = E_1'$ is not preserved by the Higgs field. By \thref{sscritrk2} it follows that $(E,\phi)$ is semistable. Conversely if $(E,\phi)$ is semistable, then the Higgs field cannot preserve ${E'}^1 $ and so $\phi_{21}$ is non-zero. Therefore $$\mathscr{H}_{\tau,[2,1]}^{2,d}(\Sigma,L)  = \mathscr{H}_{\tau}^{2,d}(\Sigma,L) \cap \mathscr{H}_{2,d}^{ss}(\Sigma,L).$$ 
To finish the proof it suffices to appeal to \thref{Higgsrefinement}, by which we have that $$ \mathscr{H}_{\tau,[0,0]}^{2,d}(\Sigma,L) =  \mathscr{H}_{\tau}^{2,d}(\Sigma,L) \cap \mathscr{H}_{2,d}^{\tau}(\Sigma,L).$$ 
\end{proof}

\subsubsection{When $r \geq 3$} 

In the above two cases, knowledge of the refined bHN stratum in which a Higgs bundle $(E,\phi)$ lies uniquely determines its HN type. However these represent very special cases, and in general there is no reason to expect that knowing the refined bHN type of a Higgs bundle $(E,\phi)$ (i.e.\ knowing `how much' of the bHN filtration is preserved by the Higgs field) determines the HN type of $(E,\phi)$. We can see this already when $r=3$, thanks to \thref{rk3strats}. That is, let $\tau$ denote a bHN type and suppose that $(E,\phi) \in \mathscr{H}_{\tau,[3,1]}^{3,d}(\Sigma,L)$. Then $\phi_{31} \neq 0$ and so none of the terms in the bHN filtration of $(E,\phi)$ are preserved. This does not mean that $(E,\phi)$ is semistable, since $\phi$ may still preserve a destabilising subbundle not appearing in the bHN filtration. A specific instance of this phenonemon has already been given in \thref{p1ex}. Thus in general the refined bHN stratum in which a Higgs bundle lies does not determine its HN type (unless it is the highest refined bHN stratum).

\section{Filtrations of the stack of Higgs bundles by quotient stacks} \label{sec:idasqs}

In this section we use the HN and bHN stratifications to study the geometry of $\mathscr{H}_{r,d}(\Sigma,L)$: we use them to filter $\mathscr{H}_{r,d}(\Sigma,L)$ by an increasing union of quotients stacks in two different ways. These filtrations will be used in Section \ref{sec:HNGIT} to relate the HN and bHN stratifications to GIT instability stratifications, and the refined bHN stratification to a Bialynicki-Birula stratification.

The fact that the semistable stratum $\mathscr{H}_{r,d}^{ss}(\Sigma,L)$ is a quotient stack follows from Nitsure's construction of the moduli space $\mathcal{M}_{r,d}^{ss}(\Sigma,L)$ of semistable Higgs bundles (see \cite {Nitsure1991}). We use the set-up of \cite{Nitsure1991} to generalise the identification of the semistable stratum as a quotient stack to a compatible identification of the remaining unstable strata as quotient stacks.  
To this end, in Section \ref{subsec:paramspace} we describe the parameter space for Higgs bundles constructed by Nitsure in \cite{Nitsure1991} and which enables the identification of the semistable stratum as a quotient stack. 
In Section \ref{subsec:idasqsbHN} we show how the bHN stratification of this parameter space can be used to filter the stack of Higgs bundles by quotient stacks (see \thref{qsforHN}). In Section \ref{subsec:idasqbHN} we prove the analogous result for the HN strata (see \thref{HNstratsetup}).

\subsection{Nitsure's parameter space for Higgs bundles} \label{subsec:paramspace}

In this section we review Nitsure's construction of a parameter space for Higgs bundles. The parameter space parametrises Higgs bundles equipped with some extra structure, and this extra structure is captured by the action of an algebraic group so that two Higgs bundles in the parameter space are isomorphic if and only if they lie in the same orbit. Nitsure's construction builds on the construction of a parameter space for vector bundles presented by Newstead in \cite{Newstead1978}, and which first appeared in work of Mumford \cite{Mumford1963} and Seshadri \cite{Seshadri1967}.

\label{subsec:Nsconstruction}

\subsubsection{Parameter space for vector bundles} 
Let $m = d+ r (1-g)$ and let $Q_{r,d}$ denote the Quot scheme parametrising quotient sheaves $q: \mathcal{O}_{{\Sigma}}^{\oplus m} \rightarrow \mathcal{E}$ of rank $r$ and degree $d$. We let $\mathcal{O}_{{\Sigma} \times Q_{r,d}}^{\oplus m} \rightarrow \mathcal{U}_{r,d}$ denote the universal quotient sheaf on ${\Sigma} \times Q_{r,d}$. There is a natural $\operatorname{GL}(m)$ action on the scheme $Q_{r,d}$ given by $$A \cdot q  = q \circ A$$ for any $A \in \GL(m)$. Define $R_{r,d} \subseteq Q_{r,d}$ to be the subset of all $ q: \mathcal{O}_{\Sigma}^{\oplus m} \to \mathcal{E}  \in Q_{r,d}$ satisfying the property that: 

\begin{enumerate}[(i)]
\item $\mathcal{E}$ is locally free (hence a vector bundle since $\Sigma$ is a curve); 
\item the canonical map $H^0(q): H^0({\Sigma}, \mathcal{O}_{\Sigma}^{\oplus m}) \rightarrow H^0({\Sigma},\mathcal{E})$ is an isomorphism. 
\end{enumerate}

The set $R_{r,d}$ is a $\GL(m)$-invariant open and reduced subset of $\mathcal{Q}_{r,d}$ (see \cite[Thm 5.3]{Newstead1978}). 
Moreover, by \cite[Lem 5.2 \& Thm 5.3]{Newstead1978} again, the following conditions are satisfied, provided $d>r(2g-1)$:  
\begin{enumerate}[(i)]
\item the family $\mathcal{U}_{r,d}|_{\Sigma \times R_{r,d}}$ of vector bundles over ${\Sigma}$ parametrised by $R_{r,d}$ has the local universal property\footnote{Given a moduli problem, a family $\mathcal{U}$ parametrised by a scheme $B$ has the local universal property if given any family $\mathcal{F}$ parametrised by a scheme $B'$, and any point $b \in B'$, there exists an open neighbourhood $U$ of $b'$ and a map $U \to B$ such that the restriction of $\mathcal{F}$ to $U$ coincides with the pull-back of $\mathcal{U}$ along the map $U \to B$.} for families of vector bundles $E$ of rank $r$ and degree $d$ satisfying: \begin{enumerate}[(a),nosep]
\item $E$  is generated by its sections; \label{cond1}
\item $H^1({\Sigma}, E) =0$; \label{cond2} 
\end{enumerate}  
\item two vector bundles in $R_{r,d}$ are isomorphic if and only if they lie in the same orbit under the $\GL(m)$ action on $R_{r,d}$. \label{2'}
\end{enumerate} 

\subsubsection{Extension to a parameter space for Higgs bundles} 

In \cite{Nitsure1991}, Nitsure extends this construction to obtain a parameter space for Higgs bundles. Let $\operatorname{HF}_{r,d}$ denote the functor from the category of $R_{r,d}$-schemes to the category of groups defined by $$(f: S \rightarrow R_{r,d}) \mapsto H^0(S, {\pi_S}_{\ast} (\operatorname{id}_{\Sigma} \times f)^{\ast} (\pi_{\Sigma}^{\ast} L \otimes \End (\mathcal{U}_{r,d}|_{\Sigma \times R_{r,d}})),$$ where $\pi_{\Sigma}: {\Sigma} \times R \rightarrow {\Sigma}$ denotes the projection onto ${\Sigma}$\footnote{The notation $\operatorname{HF}_{r,d}$ stands for `Higgs field'. See \thref{understandingH}.}.

\begin{rk}[Understanding the functor $\operatorname{HF}_{r,d}$] \thlabel{understandingH}
By the universal property of $R_{r,d}$ for vector bundles, an $R_{r,d}$-scheme $f: S \rightarrow R_{r,d}$ corresponds to a family $E_S$ of vector bundles on $\Sigma$ parametrised by $S$. Then $\operatorname{HF}_{r,d}(S \rightarrow R_{r,d})$ can be interpreted as the group of all Higgs fields $\phi_S$ such that $(E_S, \phi_S)$ is a family of Higgs bundles parametrised by $S$. In particular, if $S$ is a point, then $S \to R_{r,d}$ is a vector bundle $E$ on $\Sigma$ and its image under $\operatorname{HF}_{r,d}$ is the group $H^0(\Sigma,L \otimes \operatorname{End} E)$ of Higgs fields for $E$.  
\end{rk} 

The functor $\operatorname{HF}_{r,d}$ is representable by a scheme $f: F_{r,d} \rightarrow R_{r,d}$ which satisfies the following properties (see \cite[\S 3]{Nitsure1991}):
\begin{enumerate}[(i)]
\item the family $(E_{F_{r,d}},\phi_{F_{r,d}})$ of Higgs bundles over $\Sigma$ parametrised by ${F_{r,d}}$, determined by the identity map ${F_{r,d}} \to {F_{r,d}}$, has the local universal property for families of Higgs bundles $(E,\phi)$ of rank $r$ and degree $d$ such that $E$ is generated by its global sections and $H^1(\Sigma,E) = \{0\}$;
\item the $\GL(m)$-action on $R_{r,d}$ lifts to an action on ${F_{r,d}}$, given by $A \cdot (q,\phi) = (A \cdot q, \phi);$
\item two Higgs bundles in ${F_{r,d}}$ are isomorphic if and only if they lie in the same orbit under the $\operatorname{GL}(m)$ action on $F_{r,d}$; \label{universalfamily2}
\item given $(q: \mathcal{O}_{\Sigma}^{\oplus m} \to E,\phi) \in F_{r,d}$, there is an isomorphism $\operatorname{Aut}(E,\phi) \cong \operatorname{Stab}_{GL(m)} (q,\phi)$. \label{automorphism} 
\end{enumerate}

\begin{notation} \thlabel{notationforF}
If $\mu$ is a HN type of rank $r$ and degree $d$, then we let $F_{r,d}^{\mu} \subseteq F_{r,d}$ and $F^{r,d}_{\mu} \subseteq F_{r,d}$ denote the locally closed subscheme consisting of those Higgs bundles of HN type $\mu$ and of bHN type $\mu$ respectively. Similarly we let $F_{r,d}^{\leq \mu} \subseteq F_{r,d}$  and $F_{\leq \mu}^{r,d} \subseteq F_{r,d}$ denote the open subscheme consisting of those Higgs bundles of HN type smaller than or equal to $\mu$ and of bHN type smaller than or equal to $\mu$ respectively. 
\end{notation}

\subsection{Filtration by quotient stacks via the bHN stratification} \label{subsec:idasqsbHN}

The bHN stratification of the stack of Higgs bundles induces a filtration by open substacks, given by $$\mathscr{H}_{r,d}(\Sigma,L) = \bigcup_{\tau} \mathscr{H}_{\leq \tau}^{r,d}(\Sigma,L).$$ The aim of this section is to prove that each $\mathscr{H}^{r,d}_{\leq \tau}(\Sigma,L)$ is isomorphic to a quotient stack, and that these isomorphisms are compatible with the natural inclusions $\mathscr{H}^{r,d}_{\tau'} (\Sigma,L) \subseteq \mathscr{H}^{r,d}_{\leq \tau}(\Sigma,L)$ for $\tau' \leq \tau$.

\begin{theorem}[Filtration of the stack of Higgs bundles by quotient stacks via the bHN stratification] \thlabel{qsforHN} Let $\tau$ be a bHN type of rank $r$ and degree $d$ and fix some $e \in \mathbb{Z}$. Then for $e$ sufficiently large,  there is an isomorphism $$\mathscr{H}^{r,d}_{\leq \tau}(\Sigma,L) \cong \left[F^{r,d+re}_{\leq \tau+re} / \operatorname{GL}(m+re) \right]$$ which restricts to an isomorphism $$\mathscr{H}^{r,d}_{\tau'} (\Sigma,L) \cong \left[F^{r,d+re}_{\tau'+re} / \operatorname{GL}(m+re) \right]$$ for each $\tau' \leq \tau$. 
\end{theorem} 

Note that if we fix a bHN type $\tau$ of rank $r$ and degree $d$ and if $L_e$ is a line bundle on $\Sigma$ of degree $e$, then tensoring the underlying bundle of a Higgs bundle by $L_e$ gives an isomorphism of stacks $\mathscr{H}^{r,d}_{\leq \tau}(\Sigma,L) \cong \mathscr{H}^{r,d+re}_{\leq \tau+e}(\Sigma,L)$. Now for a given rank $r$ and degree $d$, the scheme $F_{r,d}$ parametrises Higgs bundles $(E,\phi)$ of rank $r$ and degree $d$ satisfying the property that $E$ is generated by its sections and $H^1(\Sigma,E) = 0$. The key to proving \thref{qsforHN} is therefore to show that given a bHN type $\tau$ of rank $r$ and degree $d$, then provided the degree is sufficiently large (depending on $\tau$), any Higgs bundle of bHN type $\tau$ satisfies the properties needed to be realisable as a point in $F_{r,d}$. We do so in \thref{dsufflargevb} below.

\begin{lemma} \thlabel{dsufflargevb}
 Let $\tau=(d_1/r_1,\hdots,d_s/r_s)$ denote a HN type of rank $r$ and degree $d$. 
 Then if $d > d_1/r_1(r-1) +  2 g -1$, any vector bundle $E$ of HN type $\tau' \leq \tau$ satisfies the property that: 
\begin{enumerate}[(a)]
\item $E$  is generated by its sections, and \label{cond1}
\item $H^1({\Sigma}, E) =0$. \label{cond2} 
\end{enumerate}   \label{vbs}
\end{lemma} 

\begin{rk}[Comparison with semistable vector bundles]
\thref{dsufflargevb} is a generalisation to unstable HN types of the corresponding result for the trivial HN type $\tau_0 = (d/r,\hdots, d/r)$, as stated in \cite[Lem 5.2]{Newstead1978}. That is if $d > r(2g-1)$, then any semistable vector bundle $E$ of rank $r$ and degree $d$ satisfies \ref{cond1} and \ref{cond2}. When considering unstable HN types, the bound on the degree depends not just on the rank $r$ but also on the maximal slope of the HN type. In the course of the proof of \thref{HNstratsetup} we will see that by tensoring any vector bundle of a HN type $\tau' \leq \tau$ by a line bundle on $\Sigma$ of sufficiently large degree (depending on $\tau$), we can ensure that the degree of the resulting bundle satisfies the desired inequality. 
\end{rk} 

\begin{proof}[Proof of \thref{dsufflargevb}]
We adapt Newnstead's proof of \cite[Lem 5.2]{Newstead1978} to the unstable case. We first show \ref{cond2} is satisfied if $d > d_1/r_1(r-1) +  2 g -1$. Fix a HN type $\tau$ and suppose that $E$ has HN type $\tau' \leq \tau$. If $H^1(\Sigma,E) \neq 0$, then using the isomorphism $H^1(\Sigma,E) \cong H^0(\Sigma, E^{\vee} \otimes T^{\ast} \Sigma)^{\vee}$ we can suppose that there is a non-zero global section of $E^{\vee} \otimes T^{\ast} \Sigma$, or equivalently a non-zero homomomrphism $h: E \to T^{\ast} \Sigma$. Let $H \subseteq E$ denote the subbundle generically generated by $\operatorname{ker} h$. Then $H$ has rank $r-1$, and we have: $$\operatorname{deg} H \geq \operatorname{ker} h \geq \operatorname{deg} E - \operatorname{deg} T^{\ast} \Sigma = d - (2g-2).$$  Since $E$ has HN type $\tau' \leq \tau$, the slope of any subbundle of $E$ is bounded above by the maximal slope $d_1/r_1$ of $\tau$. Thus $$d - (2g-2) \leq \operatorname{deg} H \leq \frac{d_1}{r_1}(r-1),$$ which gives the contradiction $$d \leq \frac{d_1}{r_1}(r-1) + 2 g-2.$$ 

To prove that \ref{cond1} is satisfied, let $x \in \Sigma$  and consider the exact sequence $$0 \to \mathcal{I}_x \otimes E \to E \to E_{x} \to 0,$$ where $\mathcal{I}_x$ is the sheaf of ideals corresponding to the point $x$, and $E_x$ is torsion sheaf with support $\{x\}$ corresponding to the fibre of $E$ at $x$. To show that $E$ is generated by its sections, it is enough to show that the induced map of global sections $$H^0(\Sigma,E) \to H^0(\Sigma,E_x) \cong E_x$$ is surjective. We do so by showing that $H^1(\Sigma,\mathcal{I}_x \otimes E) = 0$. Since $\mathcal{I}_x$ is torsion-free and therefore locally free as $\Sigma$ is a curve, it corresponds to a line bundle $L_x$. Moreover, since $E$ has HN type $\tau'$, the tensor product $E \otimes L_x$ has HN type $\tau' + \operatorname{deg}(L_x)$, where $\tau' + \operatorname{deg} (L_x)$ is the HN type obtained by adding the fixed constant $\operatorname{deg} (L_x)$ to each entry of $\tau'$. The exact sequence $$0 \to L_x \to \mathcal{O}_{\Sigma} \to k_x \to 0$$ where $k_x$ denotes the skyscraper sheaf at $x$ shows that the line bundle $L_x$ has degree $-1$ since $k_x$ has degree $1$. It follows that $$\operatorname{deg} (E \otimes L_x) = \operatorname{deg}(\mathcal{I}_x \otimes E) = d -r > d_1/r_1(r-1) + 2 g-1 - r = (d_1/r_1 -1)(r-1)+2 g.$$ Therefore by the argument from the previous paragraph, applied to a vector bundle of HN type $\tau'  - 1$, we must have that $H^1(\Sigma,\mathcal{I}_x \otimes E)=0$.  
\end{proof} 

We can now prove \thref{qsforHN}. 

\begin{proof}[Proof of \thref{qsforHN}]

Let $(E,\phi)$ be a Higgs bundle of bHN type $\tau' \leq \tau$. By \thref{dsufflargevb}, we know that if $d$ is sufficiently large (depending on the maximal slope appearing in $\tau$) then $E$ is generated by its sections and $H^1(\Sigma,E) = \{0\}$. Therefore $(E,\phi)$ appears as a point in $F^{r,d}_{\leq \tau}$. It follows that the restriction of the universal family $(E_{F_{r,d}},\phi_{F_{r,d}})$ to $F^{r,d}_{\leq \tau} \times \Sigma \subseteq F_{r,d} \times \Sigma$ has the local universal property for families of Higgs bundles of bHN type $\tau' \leq \tau$. Note that the same is true if we restrict to a single bHN stratum for a given $\tau' \leq \tau$. Using the fact that the stabiliser in $\GL(m)$ of a point in $F_{r,d}$ is isomorphic to the automorphism group of the corresponding Higgs bundles, we can conclude that there is an isomorphism of stacks $$\mathscr{H}^{r,d}_{\leq \tau}(\Sigma,L) \cong \left[ F^{r,d}_{\leq \tau} / G_d \right].$$ This isomorphic restricts to an isomorphism for each stratum indexed by some $\tau' \leq \tau$.  

To conclude the proof, we must deal with the case where $d \leq {d_1}/{r_1}(r-1)+ 2g-1$. This can be done by observing that if we tensor the underlying bundle of a Higgs bundle by a line bundle $L_e$ of degree $e \in \mathbb{Z}$, then as long as $e$ is sufficiently large we can ensure that the resulting Higgs bundle has degree $d+re > ({d_1}/{r_1} + e)(r-1)+ 2g-1$ which is the bound needed for \thref{dsufflargevb} applied to the HN type $\tau+re$. Thus we obtain for $e$ sufficiently large the following isomorphisms: $$\mathscr{H}^{r,d}_{\leq \tau}(\Sigma,L) \cong \mathscr{H}^{r,d+re}_{\leq \tau + e} (\Sigma,L) \cong \left[ F^{r,d+re}_{\leq \tau +e} / \GL(m+re) \right].$$ 
\end{proof}

\subsection{Filtration by quotient stacks via the HN stratification} \label{subsec:idasqbHN}

In this section we prove the analogue of \thref{qsforHN} for the HN stratification of $\mathscr{H}_{r,d}(\Sigma,L)$. That is, the HN stratification of $\mathscr{H}_{r,d}(\Sigma,L)$ induces a filtration by open substacks $$\mathscr{H}_{r,d}(\Sigma,L) = \bigcup_{\mu} \mathscr{H}_{r,d}^{\leq \mu}(\Sigma,L),$$ and we prove the following

\begin{theorem}[Filtration of $\mathscr{H}_{r,d}(\Sigma,L)$ by quotient stacks via the HN stratification]  \thlabel{HNstratsetup}
Let $\mu$ be a HN type of rank $r$ and degree $d$ and fix some $e \in \mathbb{Z}$. Let $m =  d+r(1-g)$. Then for $e$ sufficiently large, there is an isomorphism of stacks $$\mathscr{H}_{r,d}^{\leq \mu}(\Sigma,L) \cong \left[F_{r,d+re}^{\leq \mu+e} / \operatorname{GL}(m+re) \right]$$ which restricts to an isomorphism $$\mathscr{H}_{r,d}^{\mu'} (\Sigma,L) \cong \left[F_{r,d+re}^{\mu'+e} / \operatorname{GL}(m+re)) \right]$$ for each $\mu' \leq \mu$.  
\end{theorem}

The key to proving \thref{HNstratsetup} is to show that given a HN type $\mu$ of rank $r$ and degree $d$, then provided the degree is sufficiently large (depending on $\mu$), any Higgs bundle of HN type $\mu' \leq \mu$ satisfies the properties needed to be realisable as a point in $F_{r,d}$. 
This is the analogue of the result for vector bundles given in \thref{dsufflargevb}.

\begin{lemma} \thlabel{dsufflargehbs}
Let $\mu = (d_1/r_1,\hdots, d_s/r_s)$ denote a HN type of rank $r$ and degree $d$. Suppose that  $d > \left(d_1/r_1 + \operatorname{deg} L (r_{\operatorname{max}} - 1 \right)) (r-1) + 2g-1$ where $r_{\operatorname{max}} = \operatorname{max}\{r_i  \ | \ i \in \{1,\hdots, s\}\}$. Then the underlying bundle of any Higgs bundle $(E,\phi)$ of HN type $\mu' \leq \mu$ satisfies \begin{enumerate}[(a)]
\item $E$  is generated by its sections, and
\item $H^1({\Sigma}, E) =0$. 
\end{enumerate}   
\end{lemma}

\begin{rk}[Comparison with semistable Higgs bundles]
\thref{dsufflargehbs} is a generalisation to unstable HN types of the corresponding result for the trivial HN type $\mu_0 = (d/r,\hdots, d/r)$, as stated in \cite[Cor 3.4]{Nitsure1991}. That is, if $d >  (r-1)^2 \operatorname{deg} L$, then any semistable Higgs bundles $(E,\phi)$ of rank $r$ and degree $d$ satisfies \ref{cond1} and \ref{cond2}. 
\end{rk}

\begin{proof}[Proof of \thref{dsufflargehbs}]
Suppose that $d > \left(d_1/r_1 + \operatorname{deg} L (r_{\operatorname{max}} - 1 \right)) (r-1) + 2g-1$. By applying \thref{dsufflargevb} to the HN type $ \mu + \operatorname{deg} L (r_{\operatorname{max}} -1)$ obtained by adding the fixed constant $\operatorname{deg} L (r_{\operatorname{max}} -1)$ to each entry of $\mu$, we know that any vector bundle $E$ of HN type $\tau \leq \mu + \operatorname{deg} L (r_{\operatorname{max}} -1)$ satisfies \ref{cond1} and \ref{cond2}. 
If $(E,\phi)$ has HN type $\mu$, then by \thref{underlying0} \ref{part:inequalities} we know that its bHN type $\tau=(d_1'/r_1',\hdots, d_t'/r_t')$ satisfies $d_1'/r_1' \leq d_1/r_1 + \operatorname{deg} L (r_{\operatorname{max}} - 1)$. Note that the result remains true  if we assume that $(E,\phi)$ has HN type $\mu' \leq \mu$. The inequality for $d_1'/r_1'$  implies that $\tau \leq \mu + \operatorname{deg} L (r_{\operatorname{max}} -1)$, and therefore by the above paragraph we can conclude that $E$ satisfies \ref{cond1} and \ref{cond2}. 
\end{proof}

We can now prove \thref{HNstratsetup}.

\begin{proof}[Proof of \thref{HNstratsetup}]
Let $\mu= (d_1/r_1,\hdots, d_s/r_s)$ denote a HN type of rank $r$ and degree $d$ and suppose that $d > \left(d_1/r_1 + \operatorname{deg} L (r_{\operatorname{max}} - 1 \right)) (r-1) + 2g-1$. Then by \thref{dsufflargehbs}, any Higgs bundle of HN type $\mu' \leq \mu$ appears in the parameter space $F_{r,d}^{\leq \mu}$. Thus the restriction of the universal family $(E_{F_{r,d}},\phi_{F_{r,d}})$ to $F_{r,d}^{\leq \mu} \times \Sigma \subseteq F_{r,d} \times \Sigma$ has the local universal property for families of Higgs bundles of HN type $\mu' \leq \mu$. Using the fact that the stabiliser in $\GL(m)$ of a point in $F_{r,d}$ is isomorphic to the automorphism group of the corresponding Higgs bundle, we can conclude that there is an isomorphism of stacks $$\mathscr{H}_{r,d}^{\leq \mu}(\Sigma,L) \cong \left[ F_{r,d}^{\leq \mu} / \GL(m) \right].$$ 

To conclude the proof, we must deal with the case where $d \leq \left(d_1/r_1 + \operatorname{deg} L (r_{\operatorname{max}} - 1 )\right) (r-1) + 2g-1$. But this case can be dealt with exactly as in the proof of \thref{qsforHN}. 
\end{proof}

\section{The three stratifications as GIT instability and Bialynicki-Birula stratifications}  \label{sec:HNGIT}

In this section we prove that each of the three instability stratifications defined in Section \ref{sec:stratifications} can be reconstructed from algebraic stratifications arising from the GIT constructions of the moduli space of semistable vector bundles and Higgs bundles. More precisely, we show that the HN and bHN stratifications can each be recovered from a GIT instability stratification, and that the refined bHN stratification can be recovered from a Bialynicki-Birula stratification.

The GIT construction of the moduli space of semistable vector and Higgs bundles builds on the identification of $\mathscr{V}_{r,d}^{\leq \tau}(\Sigma)$ and $\mathscr{H}_{r,d}^{\leq \mu}(\Sigma,L)$  as quotient stacks: the group actions are linearised with respect to an ample line bundle on a projective completion of the parameter spaces, thus enabling the application of GIT. Sections \ref{subsec:auxGN} and \ref{subsec:auxGhatN} introduce the two projective varieties $G(r,m)^N$ and $\widehat{G}(r,m)^N$ used to obtain these projective completions. In Section \ref{subsec:HNGIT} we show that a correspondence can be established between HN types and GIT instability types such that the GIT instability stratification of the projective completions of the parameter spaces pull back to the bHN and HN stratifications respectively. In Section \ref{subsec:refbHNasBB} we show that the refined bHN stratification can be recovered from a Bialynicki-Birula stratification.

\subsection{The auxiliary variety $G(r,m)^N$} \label{subsec:auxGN}

We first define the variety and show that it admits a linear $\GL(m)$-action. We then establish a link between the variety $R_{r,d}$ parametrising vector bundles on $\Sigma$ and $G(r,m)^N$. Finally we describe the associated GIT instability stratification. 

\subsubsection{Definition and linear $\GL(m)$-action} \label{subsubsec:defofGN}  Let $G(r,m)$ denote the Grassmannian of $r$-dimensional quotients of $k^m$. There is a natural action of $\GL(m)$ on $G(r,m)$ given by multiplication on the right. 
By composing the Pl\"ucker embedding of $G(r,m)$ into a projective space with the Segre embedding of the product of the projective spaces into a projective space, we obtain an embedding $
G(r,m)^N  \hookrightarrow \mathbb{P}(V)$ where $$V: = \bigotimes_{1 \leq k \leq N} \bigwedge^{m-r} k^m.$$ The action of $\GL(m)$ on $G(r,m)^N$ can be viewed as the restriction of action of $\GL(m)$ on $\mathbb{P}(V)$ determined by the natural representation  of $\GL(m)$ on the vector space $V$. This is a natural linearisation of the $\GL(m)$ action on $G(r,m)^N$.

\subsubsection{Link between $G(r,m)^N$ and the parameter space $R_{r,d}$ for vector bundles} \label{linkwithR} Any $x \in {\Sigma}$ determines a morphism $\iota_x: R_{r,d} \rightarrow G(r,m)$ given by mapping a vector bundle $E \in R_{r,d}$ to the fibre $E_x$. By choosing $N$ points $x_1,\hdots, x_N \in {\Sigma}$, we obtain a $\GL(m)$-equivariant map $\iota_{x_1,\hdots,x_N}: R_{r,d} \rightarrow G(r,m)^N$. By \cite[Thm 5.6]{Newstead1978}, if $d$ is sufficiently large, then there exists a sequence of points in $x_1,\hdots, x_N \in \Sigma$ such  the map $\iota_{x_1,\hdots,x_N}: R_{r,d} \rightarrow G(r,m)^N$ is injective.

\subsubsection{GIT instability stratification of $G(r,m)^N$} 
The restriction of the $\GL(m)$-representation on $V$ to the maximal torus $T \subseteq \GL(m)$ consisting of diagonal matrices gives a decomposition of $V$ into weight spaces. We can choose a basis of $V$ with respect to which the action of $T$ is diagonal. The basis we choose is indexed by $N$-tuples $I=(I_k)_{k=1}^N$ of subsets $I_k$ of $\{1,\hdots,m\}$ of cardinality $r$, whose elements are listed in increasing order. For each $I_k$,  let $I_k(j)$ denote its $j$-th element. Fixing a basis $\{e_1,\hdots, e_m\}$ of $k^m$ yields a basis of $V$ consisting of vectors $$E_I : = \bigotimes_{1 \leq  k \leq  N} e_{I_k(1)} \wedge \cdots \wedge e_{I_k(r)}.$$  

We can assume that each $e_i$ is a weight vector for the action of $T$ on $k^m$. The corresponding torus weight can be identified with a character $\chi_i$ of the Lie algebra $\mathfrak{t}$ of $T$, defined by $\chi_i(a_1,\hdots,a_m) = a_i$. Then each $E_I$ is a weight vector for the action of $T$, and the associated weight is given  by $$\alpha_I : = \sum_{1 \leq k  \leq N} \sum_{1 \leq i \leq  r} \chi_{I_k(i)}.$$ We fix a $\GL(m)$-invariant inner product on $\mathfrak{t}$ of $T$ and use it to systematically identify $\mathfrak{t}$ with $\mathfrak{t}^{\vee}$.  
 
The inclusion $G(r,m)^N \hookrightarrow \mathbb{P}(V)$ maps an element $y = \langle y_1,\hdots, y_N \rangle \in G(r,m)^N$ to the point in $\mathbb{P}(V)$ with homogeneous coordinates $$y_I = \prod_{k=1}^N \operatorname{det} y_{I_k}$$ where $y_{I_k}$ is the submatrix of the $r \times m$ matrix $y_k$ obtained by picking out the $r$  columns $\{I_k(1),\hdots,I_k(r)\}$. Setting $X: = \widehat{G}(r,m)^N$, applying the definitions of Section \ref{subsubsec:GITinstab} we have that $$Z_{\beta} : = \{ y \in G(r,m)^N \ | \ y_I = 0 \text{ if } \alpha_I \cdot \beta \neq || \beta||^2\}$$ and 
$$Y_{\beta} : = \{ y \in G(r,m)^N \ | \ \text{if } \alpha_I \cdot \beta < || \beta ||^2 \text{ then } y_I = 0, \text{ and } y_I \neq 0 \text{ for some } I \text{ such that } \alpha_I \cdot \beta = || \beta ||^2\},$$ with a  retraction $p_{\beta}: Y_{\beta} \to Z_{\beta}$. 
As seen in Section \ref{subsubsec:GITinstab}, the closed subvariety $Z_{\beta}$ is invariant under the action of the subgroup $\operatorname{Stab} \beta \subseteq \GL(m)$, which has a linearisation induced by that of $\GL(m)$ acting on $G(r,m)^N$. We consider the semistable locus $Z_{\beta}^{ss}$ determined by the linearisation obtained after twisting by the character $- \chi_{\beta}$. Then $Y_{\beta}^{ss} = p_{\beta}^{-1}(Z_{\beta}^{ss})$ and we have $S_{\beta} := \GL(m) Y_{\beta}^{ss},$ giving a stratification 
$$G(r,m)^N = \bigsqcup_{\beta \in \mathcal{B}} S_{\beta}.$$

\subsection{The auxiliary variety $\widehat{G}(r,m)^N$} \label{subsec:auxGhatN}

As in the case of $G(r,m)^N$, we first define the variety and show that it admits a linear $\GL(m)$-action. We then establish a link between the variety $F_{r,d}$ parametrising Higgs bundles on $\Sigma$ and $\widehat{G}(r,m)^N$. Finally we describe the associated GIT instability stratification.

\subsubsection{Definition and linear $\GL(m)$-action} \label{subsubsec:GLactiononG-hat} 
Let $U_r^m$ denote the universal bundle on $G(r,m)$. 
The variety $\widehat{G}(r,m)$ is the projectivised bundle on $G(r,m)$ defined by $$\widehat{G}(r,m): = \mathbb{P}(\mathcal{O}_{G(r,m)} \oplus \End U_r^m) \cong \mathbb{P}((\mathcal{O}_{G} \oplus \End U_r^m) \otimes (\operatorname{det} U_r^m)^{-1}).$$ Points in $\widehat{G}(r,m)$ are equivalence classes $\langle y, \pcoor{c:\phi} \rangle$ where $y$ is an $r \times m$ matrix, $c \in k$ and  $\phi$ is an $r \times r$ matrix with $c$ and $\phi$ not both zero. The equivalence is given by $\langle y , \pcoor{c:\phi} \rangle = \langle y, \pcoor{\beta c: \beta \phi} \rangle$ for any $\beta \in k^{\ast}$ and $\langle y , \pcoor{c:\phi} \rangle = \langle \alpha y, \pcoor{(\det \alpha)^{-1} c : (\det \alpha)^{-1} \alpha \phi \alpha^{-1}} \rangle$ for any $\alpha \in \GL(r)$.

The action of $\GL(m)$ on $G(r,m)$ induces an action on $\widehat{G}(r,m)$. By \cite[Prop 2.2]{Nitsure1991}, the line bundle $\mathcal{O}_{\widehat{G}(r,m)}(1)$ on $\widehat{G}(r,m)$ is very ample. Therefore the same is true for the line bundle $\bigotimes_{k=1}^N \widehat{\pi}_k^{\ast} (\mathcal{O}_{\widehat{G}}(1))$ on $\widehat{G}(r,m)^N$ where $\widehat{\pi}_k: \widehat{G}(r,m)^N \rightarrow \widehat{G}(r,m)$ denotes the projection onto the $k^{\text{th}}$ factor. Hence there is an inclusion $$ \widehat{G}(r,m)^N \hookrightarrow \mathbb{P}(\widehat{V}^{\vee})$$ where $$\widehat{V}:= H^0 \left(\widehat{G}(r,m)^N, \bigotimes_{k=1}^N \widehat{\pi}_k^{\ast} \left( \mathcal{O}_{\widehat{G}(r,m)} (1) \right) \right).$$ By \cite[\S 2]{Nitsure1991} there is an identification $H^0( \widehat{G},\mathcal{O}_{\widehat{G}}(1)) \cong H^0 \left( G, \operatorname{det} U_r^m \oplus \operatorname{det} U_r^m \otimes \End U_r^m \right).$
It follows that $$\widehat{V} \cong \bigotimes_{k=1}^N  (\widehat{V}_1^k \oplus \widehat{V}_2^k) $$ where $$\widehat{V}_1^k = H^0(G(r,m), \det U_r^m) \cong  H^0 \left( G(r,m)^N, \pi_k^{\ast} \left( \operatorname{det} U_r^m \right) \right)$$ and $$ \widehat{V}_2^k =  H^0(G(r,m), \det U_r^m \otimes \End U_r^m) \cong H^0 \left( G(r,m)^N, \pi_k^{\ast} \left( \operatorname{det} U_r^m \otimes \End U_r^m\right) \right).$$ 

Since the action of $\GL(m)$ on $\widehat{G}(r,m)^N$ is induced by the natural representation of $\GL(m)$ on each $\widehat{V}_1^k$ and $\widehat{V}_2^k$, we obtain a linear action of $\GL(m)$ on $\widehat{G}(r,m)^N$. 

\subsubsection{Link between $\widehat{G}(r,m)^N$ and the parameter space $F_{r,d}$ for Higgs bundles}  \label{linkwithF}

The fibre in $F_{r,d}$ over a vector bundle $E \in R_{r,d}$ can be identified with $H^0({\Sigma},  \operatorname{End} (E \otimes L )$, as noted in \thref{understandingH}. Thus a point $\phi$ in this fibre determines a morphism $\phi:E \rightarrow E \otimes L$. Fixing a basis for $L_x$ induces an endomorphism $\phi_x: E_x \rightarrow E_x$, that is, a point in $\operatorname{End} E_x$. Thus we obtain a morphism $\widehat{\iota}_x: F_{r,d} \rightarrow \widehat{G}(r,m)$ given by $$(E,\phi) \mapsto \langle E_x, [\phi_x:1 ] \rangle$$ and this morphism
 lies over the morphism $\iota_x:R_{r,d} \to G(r,m)$. 
By choosing $N$ points $x_1,\hdots, x_N \in {\Sigma}$, we obtain $\GL(m)$-equivariant maps $\iota$ and $\hat{\iota}$ making the following diagram commute:  
\begin{equation}  
\begin{tikzcd}
   F_{r,d} \arrow{r}{\hat{\iota}} \arrow["f",d] &  \widehat{G}(r,m)^N \arrow["\pi",d] \\ 
  R_{r,d} \arrow{r}{\iota} & G(r,m)^N, 
\end{tikzcd} \label{comm1} 
\end{equation} 
where $\pi: \widehat{G}(r,m)^N \to G(r,m)^N$ denotes the natural projection. To simplify notation we now omit $r, d$ and $m$ from the notation. 

While the map $\widehat{\iota}: F \to \widehat{G}^N$ is not necessarily injective even for $d$ and $N$ sufficiently large, in contrast to $\iota$, it is injective when restricted to certain open subsets of $F_{r,d}$. The open subset is defined as follows. Given $A \geq 0$, let $T_{A}$ denote the set of HN types of those vector bundles $E$ of rank $r$ and degree $d$ on $\Sigma$ which satisfy the inequality $\mu(E') \leq \mu(E) + A$ for any non-zero subbundle $E'  \subseteq E$. Let $F_{A} \subseteq F$ denote the open subvariety of Higgs bundles for which the underlying bundle has HN type in $T_{A}$.   
Then by \cite[Prop 5.3]{Nitsure1991}, given any $A \geq 0$, for $d$ (depending on $A,r,g$) and $N$ (depending on $d$) sufficiently large, there exists a sequence of $N$ points $x_1,\hdots, x_N \in {\Sigma}$ such that the restriction $\widehat{\iota}: F_A \rightarrow \widehat{G}^N$ is injective.

\subsubsection{GIT instability stratification for $\widehat{G}(r,m)^N$} \label{subsubsec:gitforghat} 
The restriction of the representation of $\GL(m)$ on $\widehat{V}$ to the maximal torus $T \subseteq \GL(m)$ consisting of diagonal matrices gives a decomposition of $\widehat{V}$ into weight spaces. We can choose a basis of $\widehat{V} = \bigotimes_{k=1}^N (\widehat{V}_1^k \oplus \widehat{V}_2^k)$ with respect to which the action of $T$ is diagonal.  To define such a basis, we first consider the vector space $\widehat{V}_1^k =H^0(G(r,m), \operatorname{det} U_r^m)$. Let $I \subseteq \{1,\hdots, m\}$ be a subset of cardinality $r$, and let $I(l)$ denote the $l^{\text{th}}$ entry of $I$ for $l \in \{ 1,\hdots, r\}$, when the elements are listed in increasing order. We define a global section $S_I \in H^0(G(r,m), \det U_r^m)$ by $S_I(  y ) = \det y_J$ where $y_I$ is the $r \times r$ submatrix of $y$ obtained by taking the columns $I(1), \hdots, I(r)$. The set of all sections $S_I$ defined in this way forms a basis for $H^0(G(r,m), \det U_r^m)$.

To give a basis for $\widehat{V}_{2}^k = H^0(G(r,m),\det U_r^m \otimes \End U_r^m)$, note that given a subset $I \subseteq \{1,\hdots, m\}$ as above and elements $i,j \in \{1,\hdots, m\}$, we can define a section $\sigma_I^{ij}$ of $\End U_r^m$ on the open subset of $G(r,m)$ consisting of all $y$ such that $S_I(y) \neq 0$. This section is defined by $\sigma_I^{ij}(y) = y_I \sigma_{ij} y_I^{-1}$ where $\sigma_{ij}$ is the $r \times r$ matrix with zeroes everywhere except for a one in the $(ij)$-th position. This allows us to define a global section $S_I^{ij} \in H^0(G(r,m), \det U_r^m \otimes \End U_r^m)$ by $$S_I^{ij} (y) = S_I(y) \cdot \sigma_I^{ij}(y).$$ Sections of the form $S_I^{ij}$ form a 
basis for the vector space $H^0(G(r,m), \det U_r^m \otimes \End U_r^m)$.

More generally, we let $I$ and $J$ denote $N$-tuples of subsets $I_1,\hdots, I_N, J_1,\hdots, J_N \subseteq \{1,\hdots,m\}$ of cardinality $r$ respectively, and we let $i$ and $j$ correspond to $N$ choices $i_1,\hdots, i_N$ and $j_1, \hdots, j_N$ respectively of elements in $\{1,\hdots, r\}.$ The vector space $\bigotimes_{k=1}^N  \left( V_{1k} \oplus V_{2k} \right)$ then admits a basis consisting of elements $$S_{I}^{ij} := \bigotimes_{k=1}^N S_{I_k}^{i_k j_k}  \text{ and } S_{J} := \bigotimes_{k=1}^N S_{J_k}.$$ We let $\check{S}_J$ and $\check{S}_I^{ij}$  denote the associated dual basis vectors for $\widehat{V}_1^{\vee} : = \bigotimes_{k=1}^N V_{1k}^{\vee}$ and $\widehat{V}_2^{\vee}: = \bigotimes_{k=1}^N V_{2k}^{\vee}$ respectively.

With the above set-up, we can compute the weights for the diagonalised action of the maximal torus $T \subseteq \operatorname{GL}(m)$ consisting of diagonal matrices on $\widehat{V}^{\vee}$.  The basis vectors $\check{S}_I^{ij}$ and $\check{S}_J$ were chosen precisely because they are weight vectors for the $T$-action. Again let $\mathfrak{t}$ denote the Lie algebra of $T$. For each $l = 1, \hdots, m$ let $\chi_l$ be the character of $T$ given by $\operatorname{diag}(a_1,\hdots, a_m) \mapsto a_l$, which upon differentiation at the identity produces an element of $\mathfrak{t}^{\ast} \cong \mathfrak{t}$. We also call this element $\chi_l$ for ease of notation. The torus $T$ acts on $\check{S}_I$ and $\check{S}_I^{ij}$ by multiplication by the following respective characters:\begin{equation} \widehat{\alpha}_I := - \sum_{k=1}^N \sum_{l \in I_k} \chi_l \text{  and  } \widehat{\alpha}_{I}^{ij} : = \sum_{k=1}^N \left( \chi_{I_k(j_k)} - \chi_{I_k(i_k)}  - \sum_{l \in I_k} \chi_l \right). \label{characters}
\end{equation}

The embedding $$\widehat{G}(r,m)^N \hookrightarrow \mathbb{P}(\widehat{V}^{\vee}) $$ maps a point $\widehat{y}=(\langle y_1, [c_1:\phi_1] \rangle, \hdots, y_N,[ c_N:\phi_N]\rangle)$ to the point with homogeneous coordinates \begin{equation} \widehat{y}_{J} = \prod_{k=1}^N c_k \operatorname{det} {y_{k}}_{J_k} \text{ and } \widehat{y}_{I}^{ij} = \prod_{k=1}^N  \operatorname{det} {y_{k}}_{I_k} \operatorname{tr}\left( \phi_k {y_k}_{I_k} \sigma_{I_k}^{i_kj_k} {y_k}_{I_k}^{-1}  \right). \label{coords}
\end{equation}

Setting $\widehat{X}:= \widehat{G}(r,m)^N$ and letting $\widehat{Z}_{\beta},$ $\widehat{Y}_{\beta}$ and $\widehat{S}_{\beta}$ denote the analogues of $Z_{\beta}$, $Y_{\beta}$ and $S_{\beta}$ for $X$, the above computations can be used to explicitly describe the GIT unstable strata $\widehat{S}_{\beta} = \GL(m) \widehat{Y}_{\beta}^{ss}$ in the GIT instability stratification $$\widehat{X} : = \widehat{G}(r,m)^N = \bigsqcup_{\beta \in \mathcal{B}} \widehat{S}_{\beta}$$ associated to the linear action of $\GL(m)$ on $\widehat{X}$.

\subsection{The bHN and HN stratifications as GIT instability stratifications} \label{subsec:HNGIT}

In Sections \ref{linkwithR} and \ref{linkwithF} above we defined for any choice of $N$ points on $\Sigma$ maps $\iota:R_{r,d} \to G(r,m)^N$ and $\widehat{\iota}: F_{r,d} \to \widehat{G}(r,m)^N$ making the diagram \eqref{comm1} commute. 
In this section we show that the bHN and HN stratifications of $F_{r,d}$ (more precisely of $F_{\leq \tau}^{r,d}$ and $F_{r,d}^{\leq \mu}$ for a fixed bHN $\tau$ and HN type $\mu$) coincide with the pull-backs of the GIT stability stratifications of $G(r,m)^N$ and $\widehat{G}(r,m)^N$ under the maps $\pi \circ \widehat{\iota}$ and $\widehat{\iota}$ respectively, provided the degrees appearing in $\tau$ and $\mu$ and $N$ are sufficiently large. In Section \ref{subsubsec:bHNtoGIT} we consider the case where we fix the bHN type; the case where we fix the HN type instead is treated in Section \ref{subsubsec:HNtoGIT}. 

\subsubsection{Relating bHN types to GIT instability types}  \label{subsubsec:bHNtoGIT}

In \cite[\S 11]{Kirwan1984} a correspondence $\tau \mapsto \beta(\tau)$ is established between HN types $\tau$ for vector bundles and GIT instability types $\beta$ associated to the linear action of $\GL(m)$ on $G(r,m)^N$. Note that the latter are elements of the $m$-dimensional vector space $\mathfrak{t} \cong \mathfrak{t}^{\vee}$. The corresondence is the following. 

\begin{defn}[$\tau \mapsto \beta(\tau)$ correspondence] \thlabel{mubetacorresp}
 Let $\tau = (d_1/ r_1, \hdots, d_1/ r_1, \hdots, d_s/r_s)$ denote a HN type. Let $k_i = - N r_i$ and $m_i = d_i+ r_i(1-g)$.  Set $k = \sum_{i = 1}^s k_{i}$. We define an associated vector $$\beta(\tau) := \left( \frac{k_1}{m_1}, \hdots, \frac{k_1}{m_1}, \frac{k_2}{m_2}, \hdots, \frac{k_s}{m_s} \right)$$ where each $k_i/ m_i$ appears $m_i$ times.  
\end{defn}
The sequence of strict inequalities $d_1/r_1 > \cdots > d_s/r_s$ implies a sequence of strict inequalities $k_1/m_1 > \cdots > k_s/m_s$ (see \cite[Prop 16.9]{Kirwan1984}). It is shown in \cite[Cor 11.5]{Kirw1985} that given a HN type $\tau = (d_1/r_1,\hdots, d_s/r_s)$, if the degrees $d_i$ are sufficiently large, there exists an $N \in \mathbb{N}_{>0}$ and $N$ points $x_1,\hdots, x_N \in \Sigma$ such that a vector bundle $E \in R_{r,d}$ has HN type $\tau' \leq \tau$ if and only if $\iota(E)$ lies in $S_{\beta(\tau')}$, where $\iota: R_{r,d} \to G(r,m)^N$ is the map determined by $x_1,\hdots,x_N$. 
\thref{bHNtoGIT} below extends this result to bHN types for Higgs bundles. 

\begin{prop}[Relating bHN types to GIT instability types] \thlabel{bHNtoGIT}
Let $\tau = (d_1/r_1,\hdots, d_s/r_s)$ denote a bHN type of rank $r$ and degree $d$. Then if the degrees $d_i$ are sufficiently large, there exists an $N \in \mathbb{N}_{>0}$ and $N$ points $x_1,\hdots, x_N \in \Sigma$ such that a Higgs bundle $(E,\phi) \in F_{r,d}$ has bHN type $\tau' \leq \tau$ if and only if $\widehat{\iota}(E,\phi)$ lies in $\pi^{-1}(S_{\beta(\tau')})$, where $\widehat{\iota}: F_{r,d} \to \widehat{G}(r,m)^N$ is the map determined by $x_1, \hdots, x_N$. 
\end{prop} 

\begin{proof} 
By \cite[Cor 11.5]{Kirw1985}, if the degrees $d_i$ are sufficiently large then for $N$ sufficiently large we can choose $N$ points on $\Sigma$ such that a vector bundle $E \in R_{r,d}$ has HN type $\tau' \leq \tau$ if and only if $\iota(E)$ lies in $S_{\beta(\tau')}$. The same $N$ points give a map $\widehat{\iota}: F_{r,d} \to \widehat{G}(r,d)^N$ and we have the following commutative diagram:\begin{equation}  
\begin{tikzcd}
f^{-1}(R_{r,d}^{\tau'}) \arrow["\widehat{\iota}",r]  \arrow[hookrightarrow]{d}& \pi^{-1}(S_{\beta(\tau')}) \arrow[hookrightarrow]{d} \\ 
   F_{r,d} \arrow{r}{\widehat{\iota}} \arrow["f",d] &  \widehat{G}(r,m)^N \arrow["\pi",d] \\ 
  R_{r,d} \arrow{r}{\iota} & G(r,m)^N \\
  R_{r,d}^{\tau'} \arrow[hookrightarrow]{u} \ar["\iota",r] & S_{\beta(\tau')}, \arrow[hookrightarrow]{u}
\end{tikzcd} \label{comm2} 
\end{equation} 
where $\tau'$ is a HN type with $\tau' \leq \tau$. Indeed, if $(E,\phi) \in F_{r,d}$ has bHN type $\tau' \leq \tau$, then we know that $\iota \circ f (E,\phi) = \iota(E)$ lies in $S_{\beta(\tau')}$ where $f: F_{r,d} \to R_{r,d}$ is the forgetful map. Therefore $\widehat{\iota}(E,\phi)$ must lie in $\pi^{-1}(S_{\beta(\tau')})$, where $\pi: \widehat{G}(r,m)^N \to G(r,m)^N$ is the natural projection. Conversely, if $\widehat{\iota}(E,\phi)$ lies in $S_{\beta(\tau')}$ for some $\tau' \leq \tau$, then $\pi(\widehat{\iota}(E,\phi)) = \iota(f(E,\phi)) = \iota(E)$ lies in $S_{\beta(\tau')}$. Therefore $E \in R_{r,d}$ has HN type $\tau'$ and so $(E,\phi)$ has bHN type $\tau'$. 
\end{proof} 

We conclude this section by noting its reformulation in terms of the stack of Higgs bundles.

\begin{cor}[bHN stratification in terms of a GIT instability stratification] \thlabel{corbHN}
Fix a bHN type $\tau$ of rank $r$ and degree $d$. 
Then under the isomorphism $$\mathscr{H}^{r,d}_{\leq \tau}(\Sigma,L) \cong \left[F^{r,d+re}_{\leq \tau+e} / \operatorname{GL}(m+re) \right]$$ established in \thref{qsforHN}, for $e$ sufficiently large the bHN stratification of  $\mathscr{H}^{r,d}_{\leq \tau}(\Sigma,L)$ coincides with the intersection of the pull-back to $\left[\widehat{G}(r,m+re)^N / \GL(m+re) \right]$ of the GIT instability stratification of $\left[G(r,m+re)^N / \GL(m+re) \right]$ with its substack $\left[F^{r,d+re}_{\leq \tau+e} / \GL(m+re) \right]$. 
\end{cor} 

\begin{proof}
The result follows immediately from \thref{bHNtoGIT}, by observing that if $(E,\phi)$ has bHN type $\tau = (d_1/r_1,\hdots, d_s/r_s)$ and $L_e$ is a line bundle on $\Sigma$ of degree $d$, then $(E \otimes L_e, \phi \otimes \operatorname{id}_{L_e})$ has bHN type $\tau + e = ((d_1+r_1e)/r_1,\hdots, (d_s+ r_s e)/r_s)$. Thus by choosing $e$ sufficiently large we can ensure that the degrees appearing in the bHN type $\tau + e$ are large enough so that \thref{bHNtoGIT} can be applied.  
\end{proof}

\subsubsection{Relating HN types to GIT instability types} \label{subsubsec:HNtoGIT}

In this section we prove the analogue of \thref{bHNtoGIT} for HN types, namely

\begin{prop} \thlabel{bigprop} 
Let $\mu = (d_1/r_1,\hdots, d_s/r_s)$ denote a HN type of rank $r$ and degree $d$. Then if the degrees $d_i$ are sufficiently large, for $N$ sufficiently large there  exists $N$ points $x_1,\hdots, x_N \in {\Sigma}$ such that a Higgs bundle $(E,\phi) \in F_{r,d}$ has HN type $\mu' \leq \mu$ if and only if $\widehat{\iota}(E,\phi)$ lies in $\widehat{S}_{\beta(\mu')}$, where $\widehat{\iota}:F_{r,d} \to \widehat{G}(r,m)^N$ is the map determined by $x_1,\hdots, x_N$.  
\end{prop}

\begin{rk}[Comparison with Kirwan's result for vector bundles] 
We note that the above \thref{bigprop} is the generalisation to Higgs bundles of the analogous result for vector bundles, established in \cite[Cor 11.5]{Kirw1985} and which we have already used to prove \thref{bHNtoGIT}.  
\end{rk}

Before proving \thref{bigprop}, we note its reformulation in terms of the stack of Higgs bundles. 

\begin{cor}[HN stratification as a GIT instability stratification] \thlabel{HNasGIT}
Fix a HN type $\mu$ of rank $r$ and degree $d$. Then under the isomorphism $$\mathscr{H}^{r,d}_{\leq \mu}(\Sigma,L) \cong \left[F_{r,d+re}^{\leq \mu+e} / \operatorname{GL}(m+re) \right]$$ established in \thref{qsforHN}, for $e$ sufficiently large the HN stratification of  $\mathscr{H}_{r,d}^{\leq \mu}(\Sigma,L)$ coincides with the intersection of the GIT instability stratification of $\left[\widehat{G}(r,m+re)^N / \GL(m+re) \right]$ with its substack $\left[F_{r,d+re}^{\leq \mu+e} / \GL(m+re) \right]$. 
\end{cor} 

\begin{proof} 
The result follows from \thref{bigprop}, just as \thref{corbHN} follows from \thref{corbHN}. 
\end{proof}

We now turn to proving \thref{bigprop}. Given a GIT instability stratification $X = \bigsqcup_{\beta \in \mathcal{B}} S_{\beta}$ for a linear $G$-action on a projective variety $X$, by definition $S_{\beta} = G Y_{\beta}^{ss}$ (see Section \ref{subsubsec:GITinstab}). Hence to show that $x \in S_{\beta}$, it suffices to show that $g \cdot x \in Y_{\beta}^{ss}$ for some $g \in G$. This is the strategy we use to prove the forward implication of \thref{bigprop}, and we do so in two steps: 
\begin{enumerate}[(i)]
\item firstly we show that $g \cdot x \in Y_{\beta}$ for some $g \in G$ (see \thref{lemma1});
\item secondly we show that $p_{\beta}(g \cdot x) \in Z_{\beta}^{ss}$ (see \thref{lemma2}). 
\end{enumerate}

\begin{lemma} \thlabel{lemma1}
Let $\mu=(d_1/r_1,\hdots, d_s/r_s)$ denote a HN type of rank $r$ and degree $d$. Then if the degrees $d_i$ and $N$ are sufficiently large, there are $N$ points $x_1, \hdots, x_N \in \Sigma$ such that if $(E,\phi) \in F_{r,d}$ has HN type $\mu$, then $g \cdot \widehat{\iota}(E,\phi) \in \widehat{Y}_{\beta(\mu)}$ for some $g \in \GL(m)$, where $\widehat{\iota}: F_{r,d} \to \widehat{G}(r,m)^N$ is the map associated to $x_1,\hdots, x_N$. 
\end{lemma}

\begin{proof} [Proof of \thref{lemma1}] 
We wish to show that for $d_i$ and $N$ sufficiently large, the image $\widehat{\iota}(g \cdot (E,\phi))$ lies in $\widehat{Y}_{\beta(\mu)}$. We do so by using the set-up described in Section \ref{subsubsec:gitforghat}, which can be summarised as follows. There is a $\GL(m)$-equivariant embedding $$\widehat{G}(r,m)^N \hookrightarrow \mathbb{P}(\widehat{V}_1^{\vee} \oplus \widehat{V}_2^{\vee})$$ with $\GL(m)$ acting via a representation on $V_1$ and $V_2$. Moreover, there is a basis of torus weight vectors $\{\check{S}_I\}$ for $\widehat{V}_1^{\vee}$, indexed by $N$-tuples $I$ of subsets $I_k \subseteq \{1,\hdots,m\}$ of cardinality $r$. Similarly there is a basis of torus weight vectors $\{\check{S}_I^{i,j}\}$ for $\widehat{V}_2$, indexed by triples $(I,i,j)$ where $I$ is an $N$-tuple  as above together and $i = (i_1,\hdots,i_N), j = (j_1,\hdots,j_N) \in \{1,\hdots, r\}^N$. 
The central torus in $\GL(m)$ acts with weight $\widehat{\alpha}_I \in \mathfrak{t}^{\vee}$ on $\widehat{E}_I$ and with weight $\widehat{\alpha}_I^{ij}$ on $\widehat{E}_I^{ij}$ (both weights are described explicitly at \eqref{characters}. Given a point $\widehat{y} \in \widehat{G}(r,m)^N$, its image in $\mathbb{P}(V_1^{\vee} \oplus V_2^{\vee})$ has homogeneous coordinates $\widehat{y}_I$ and $\widehat{y}_I^{ij}$, which can also be written down explicitly -- see \eqref{coords}. 

Consider the HN filtration $ 0 = E^0 \subset E^1 \subset \cdots \subset E^s = E$ of $(E,\phi)$. By taking global sections we obtain a filtration $0 = H^0(\Sigma,E^0) \subset H^0(\Sigma,E^1) \subset \cdots \subset H^0(\Sigma,E^s) = H^0(\Sigma,E)$ of $H^0(\Sigma,E)$. Since $d > r(g-1)$ we have that $h^1(\Sigma,E) = 0$ and so the Riemann-Roch formula implies that $h^0(\Sigma,E) = d- r(1-g) = m$. We can thus fix an isomorphism $H^0(\Sigma,E) \cong k^m$, and from here on we identify the two vector spaces. For ease of notation we set $M^{\gamma} = H^0(\Sigma,E^{\gamma})$, so that we have a filtration $0 = M^0 \subset M^1 \subset \cdots \subset M^s = k^m$ of $k^m$. We note that $\operatorname{dim} M^{\gamma} / M^{\gamma-1} = m_{\gamma}$ for each $\gamma = 1,\hdots, s$. Let $e_1, \hdots, e_m$ denote the standard basis of $k^m$. 
By replacing $\widehat{y}$ by $g \cdot \widehat{y}$ for an appropriate $g \in \operatorname{GL}(m)$, we may assume that $M^{\gamma}$ is spanned by the basis vectors $\{e_l  \ | \ l \leq \operatorname{dim} M^{\gamma} \}$ for $\gamma = 1 ,\hdots, s$. For ease of notation we relabel this new $g \cdot \widehat{y}$ as $\widehat{y}$.

Showing that $ \widehat{y} \in \widehat{Y}_{\beta(\mu)}$ for sufficiently large $d$ and $N$ requires showing that the following three statements are true for $d$ and $N$ sufficiently large: \begin{enumerate}
\item if $\widehat{\alpha}_I \cdot \beta(\mu) < || \beta(\mu) ||^2$ then $\widehat{y}_I = 0$; \label{pfstep1}
\item if $\widehat{\alpha}_{I}^{ij}  \cdot \beta(\mu') < || \beta(\mu) ||^2$ then $\widehat{y}_I^{ij} = 0$; \label{pfstep2}
\item there is a non-zero coordinate, either $\widehat{y}_I$ or $\widehat{y}_I^{ij}$, such that $\widehat{\alpha}_I \cdot \beta(\mu) = || \beta(\mu) ||^2$ or $\widehat{\alpha}_I^{ij} \cdot \beta(\mu) = || \beta (\mu)||^2$.  \label{pfstep3}
\end{enumerate}

The proof of the above three statements are as follows. \\

\paragraph*{Proof of \eqref{pfstep1}} 
Suppose that $\widehat{y}_I\neq 0$ for some $I$. We wish to show that $\widehat{\alpha}_I \cdot \beta(\mu) \geq ||\beta(\mu)||^2$. Recall that $\beta(\mu') = (k_1/m_1,\hdots, k_s/m_s)$ with each $k_i/m_i$ repeated $i$ times, where $k_i = -N r_i$ and $m_i = d_i + r_i (1-g)$.  Using the expression for $\widehat{\alpha}_I$ given in \eqref{characters}, we have: \begin{align} 
\widehat{\alpha}_I \cdot \beta(\mu) & = - \sum_{k=1}^N \left(\sum_{l \in J_k} \chi_l \right) \cdot \beta \nonumber \\  
& =  - \sum_{k=1}^N \sum_{\gamma = 1}^s \frac{k_{\gamma}}{m_{\gamma}} \# \{l \in \{1,\hdots,r\} \ | \ m_1 + \cdots + m_{\gamma -1} < I_k(l) \leq m_1 + \cdots + m_{\gamma} \}  \nonumber  \\
& =  - \sum_{k=1}^N \sum_{\gamma = 1}^s \frac{k_{\gamma}}{m_{\gamma}} \# \{l \in \{1,\hdots,r\} \ | \ \operatorname{dim} M^{\gamma-1} < I_k(l) \leq \operatorname{dim} M^{\gamma}\} \nonumber  \\ 
& =  - \sum_{\gamma=1}^N \sum_{\gamma=1}^{s-1} \left( \frac{k_{\gamma}}{m_{\gamma}} - \frac{k_{\gamma+1}}{m_{\gamma+1}} \right) \# \{ l  \in \{1,\hdots, r\}    \ | \ I_k(l) \leq \operatorname{dim}M^{\gamma} \} - N r \frac{k_s}{m_s}. \label{lasteq}
\end{align} 

Since $\widehat{y}_I = \prod_{k=1}^N \operatorname{det} (y_k)_{I_k}$ by \eqref{coords}, and since by assumption $\widehat{y}_I \neq 0,$ it follows that $\operatorname{det} (y_k)_{I_k} \neq 0$ for all $k = 1,\hdots, N$. Hence the columns $I_k(1),\hdots, I_k(r)$ of the $r \times m$ matrix $y_k$ are linearly independent vectors in $k^r$ for each $k$. In particular, the columns $I_k(l)$ such that $I_k(l) \leq \operatorname{dim} M^{\gamma}$ are linearly independent vectors in $k^r$, and these vectors are contained in $y_k(M^{\gamma})$. Thus the number of integers $l \in \{1,\hdots,r\}$ such that $I_k(l) \leq \operatorname{dim} M^{\gamma}$ is bounded above by the dimension of the image $y_k(M^{\gamma}) \subseteq k^r$. But $\operatorname{dim} y_k(M^{\gamma}) = \operatorname{rk} E^{\gamma} = r_1 + \hdots + r_{\gamma},$ provided the degrees $d_{\gamma}$ are sufficiently large. Hence following on from \eqref{lasteq}, we obtain that \begin{align*} 
\widehat{\alpha}_I \cdot \beta(\mu) & \geq  - \sum_{k=1}^N \sum_{\gamma = 1}^{s-1} \left( \frac{k_{\gamma}}{m_{\gamma}} - \frac{k_{\gamma+1}}{m_{\gamma+1}} \right) \left( r_1 + \cdots + r_{\gamma} \right) - Nr \frac{k_s}{m_s} \\
& = -  \sum_{\gamma = 1}^{s-1} \left( \frac{k_{\gamma}}{m_{\gamma}} - \frac{k_{\gamma+1}}{m_{\gamma+1}} \right)N  \left( r_1 + \cdots + r_{\gamma} \right) - N r \frac{k_s}{m_s}\\
& = \sum_{\gamma = 1}^{s-1} \left( \frac{k_{\gamma}}{m_{\gamma}} - \frac{k_{\gamma+1}}{m_{\gamma+1}} \right) \left( k_1 + \cdots + k_{\gamma} \right) + k \frac{k_s}{m_s}  = \sum_{\gamma=1}^{s} \frac{k_{\gamma}^2}{m_{\gamma}} = || \beta(\mu) ||^2. 
\end{align*}

\paragraph*{Proof of \eqref{pfstep2}} Now suppose that $\widehat{y}_I^{ij} \neq 0$. By \eqref{characters} we have that $\widehat{\alpha}_{I}^{ij} = \widehat{\alpha}_I + \sum_{k=1}^N (\chi_{I_k(j)} - \chi_{I_k(i)})$. It follows from \eqref{pfstep1} that $\widehat{\alpha}_I \cdot \beta(\mu) \geq || \beta(\mu) ||^2$ since the assumption that $\widehat{y}_I^{ij} \neq 0$ implies in particular that $\operatorname{det} (y_k)_{I_k} \neq 0$ for each $k$, by the expression for $\widehat{y}_I^{ij}$ given in \eqref{coords}. Thus it suffices to prove that \begin{equation} \sum_{k=1}^N ( \chi_{I_k(j_k)} - \chi_{I_k(i_k)} ) \cdot \beta(\mu')  \geq 0. \label{dotproduct} 
\end{equation} 
We start by calculating $\sum_{k=1}^N \chi_{I_k(j_k)} \cdot \beta(\mu)$. First note that $\chi_{I_k(j_k)} \cdot \beta(\mu) = k_{\gamma} / m_{\gamma}$ where $\gamma \in \{1,\hdots, s\}$ is such that $\operatorname{dim} M^{\gamma-1} < I_k(j_k) \leq \operatorname{dim} M^{\gamma}$. Therefore \begin{align*} \sum_{k=1}^N \chi_{I_k(j_k)} \cdot \beta(\mu) & = \sum_{\gamma=1}^{s} \frac{k_{\gamma}}{m_{\gamma}} \# \{k \in \{1,\hdots,N\}  \ | \ \operatorname{dim} M^{\gamma-1} < I_k(j_k) \leq \operatorname{dim} M^{\gamma} \} \\
& = \sum_{\gamma=1}^{s-1} \left( \frac{k_{\gamma}}{m_{\gamma}} - \frac{k_{\gamma+1}}{m_{\gamma+1}}\right) \# \{ k \in \{1,\hdots, N \}  \ | \ I_k(j_k) \leq \operatorname{dim} M^{\gamma} \} +N \frac{k_{\gamma}}{m_{\gamma}}. 
\end{align*} 
The same equality is valid with $j_k$ replaced by $i_k$. It follows that \begin{align*}  \sum_{k=1}^N ( \chi_{I_k(j_k)} - \chi_{I_k(i_k)} ) \cdot \beta(\mu) & =  \sum_{\gamma = 1}^{s-1} \left( \frac{k_{\gamma}}{m_{\gamma}} - \frac{k_{\gamma+1}}{m_{\gamma+1}} \right) ( \# \{ k \in \{1,\hdots,N\}  \ | \ I_k(j_k) \leq \operatorname{dim} M^{\gamma}  \}  \\
& \hspace{3.5cm} - \# \{ k  \in \{1,\hdots, N\} \ | \   I_k(i_k) \leq \operatorname{dim} M^{\gamma} \} ).
\end{align*} 
Since $k_{\gamma}/m_{\gamma} > k_{\gamma+1}/{m_{\gamma+1}}$ for every $\gamma = 1,\hdots, s-1$, we need only show that  \begin{equation} \# \{ k \in \{1,\hdots,N\}  \ | \ I_k(j_k) \leq \operatorname{dim} M^{\gamma}  \}  \geq  \# \{ k  \in \{1,\hdots, N\} \ | \   I_k(i_k) \leq \operatorname{dim} M^{\gamma} \} 
\label{inequalitytoshow} 
\end{equation}  for all $\gamma=1,\hdots,s$. 

By assumption $\widehat{y}_I^{ij} \neq 0$, therefore using the expression for $\widehat{y}_I^{ij}$ given in \eqref{coords} we obtain that $$\operatorname{tr}(\phi_k {y_k}_{I_k} \sigma_{i_k j_k} {y_k}_{I_k}^{-1} ) \neq 0 $$ for every $k  =1,\hdots, N$. By the cyclic property of the trace, the equality $$\operatorname{tr}(\phi_k {y_k}_{I_k} \sigma_{i_k j_k} {y_k}_{I_k}^{-1} ) = \operatorname{tr} ((y_k)_{I_k}^{-1}  \phi_k (y_k)_{I_k} \sigma_{i_k j_k} )$$ holds. Moreover, since by definition $\sigma_{i_kj_k}$ is the $r \times r$ matrix with zeroes everywhere except for a one in the $(i_kj_k)$-th entry, it follows that for any $r \times r$ matrix $A$ there is an equality $$ \operatorname{tr} (A \sigma_{i_k j_k}) = a_{j_k i_k}.$$ Setting $A = (y_k)_{I_k}^{-1}  \phi_k (y_k)_{I_k}$ and using the fact that $\widehat{y}_{I_k}^{i_k j_k} \neq 0$, we can deduce that the $(j_k i_k)$-th entry of the matrix $(y_k)_{I_k}^{-1}  \phi_k (y_k)_{I_k}$ is non-zero.  

The maps $\phi_k$ are induced by the Higgs field $\phi$, and each $E^{\gamma}$ in the HN filtration of $(E,\phi)$ is preserved by $\phi$. It follows that each $y_k(M^{\gamma})$ is preserved by $\phi_k:k^r \to k^r$. This property of $\phi_k$ implies that the matrix $(y_k)_{I_k}^{-1}  \phi_k (y_k)_{I_k}$ is block upper-triangular.  
As a result, if $i_k \leq \operatorname{dim} y_k(M^{\gamma})$, then the same must be true for $j_k$ since the $(j_k i_k)$-th entry of $(y_k)_{I_k}^{-1}  \phi_k (y_k)_{I_k}$ cannot be zero. 

Suppose that $I_k(i_k) \leq \operatorname{dim} M^{\gamma}$. Note that given any $l \in \{1,\hdots, r\}$, the inequality $l \leq \operatorname{dim} y_k(M^{\gamma})$ holds if and only if $I_k(l) \leq \operatorname{dim} M^{\gamma}$. This follows from the fact that we have chosen a basis $\{e_1, \hdots ,e_m\}$ for $k^m$ such that $M^{\gamma}$ is spanned by the basis vectors $\{e_l  \ | \ l \leq \operatorname{dim} M^{\gamma} \}$. Therefore $i_k \leq  \operatorname{dim} y_k(M^{\gamma})$. The result from the above paragraph implies that $j_k \leq  \operatorname{dim} y_k(M^{\gamma})$, which is equivalent to the inequality $I_k(j_k) \leq \operatorname{dim} M^{\gamma}$. Hence 
the inequality given in \eqref{inequalitytoshow} is true, and the required ineuqality of \eqref{dotproduct} follows. \\

\paragraph*{Proof of \eqref{pfstep3}} 
It remains only to show that $\widehat{y}$ has a non-zero projective coordinate such that the dot product of the corresponding weight with $\beta(\mu)$ is equal to $||\beta(\mu)||^2$. It follows from the argument given in \eqref{pfstep1} that 
$\widehat{\alpha}_I \cdot \beta(\mu) = || \beta(\mu) ||^2$ if and only if $$k_{\gamma} = \# \{ (l,k) \in \{1,\hdots,s\} \times \{1,\hdots,N\}  \ | \ \operatorname{dim} M^{\gamma-1} <  I_k(l) \leq \operatorname{dim} M^{\gamma}   \}$$ for every $\gamma = 1,\hdots, s$ and every $k=1,\hdots, N$. Similarly, it follows from the argument given in \eqref{pfstep2} that  $\widehat{\alpha}_I^{ij}  \cdot \beta(\mu) = || \beta||^2$ if and only if $$k_{\gamma} = \# \{ (l,k) \in \{1,\hdots,s\} \times \{1,\hdots,N\} \ | \ \operatorname{dim} M^{\gamma-1} <  J_k(l) \leq \operatorname{dim} M^{\gamma}   \}$$ for every $\gamma = 1,\hdots, s$ and every $k=1,\hdots, N$ and $$\operatorname{dim} M^{\gamma-1} < I_k(i_k) \leq \operatorname{dim} M^{\gamma} \Leftrightarrow  \operatorname{dim} M^{\gamma-1} < I_k(j_k) \leq \operatorname{dim} M^{\gamma}$$ for every $k = 1,\hdots,N$. 

To identify such a non-zero coordinate, consider the HN graded $\operatorname{gr}(E,\phi) = (E_1, \phi) \oplus \cdots \oplus (E_s, \phi_s)$. Each $(E_i,\phi_i)$ is a semistable Higgs bundles of rank $r_i$ and degree $d_i$ and can be identified as a point in $F_{r_i,d_i}$ for $d_i$ sufficiently large. Just as the map $\widehat{\iota}: F_{r,d} \to \widehat{G}(r,m)^N$ was defined given a choice of $N$ points on $\Sigma$, a map $\widehat{\iota_i}: F_{r_i,d_i} \to \widehat{G}(r_i,m_i)^N$ can be defined analogously for each $i = 1, \hdots ,s$. Let $\widehat{Y}_i= ( \langle  y_{i1}, [c_{i1}:\phi_{i1}] \rangle, \hdots \langle y_{iN}, [c_{iN}:\phi_{iN}] \rangle ) $ denote the image of $(E_i,\phi_i)$ under $\widehat{\iota}_i$ for each $i$, so that there is an associated element 
\begin{equation} \widehat{Y} : = (\langle y_{11}, [c_{11}:\phi_{11}] \rangle, \hdots,  \langle y_{sN}, [c_{sN}:\phi_{sN}]\rangle ) \in \prod_{i=1}^s \widehat{G}(r_i,m_i)^N. \label{bigyhat} 
\end{equation} 

Each $\widehat{G}(r_i,m_i)^N$ can be embedded into a large projective space, analogously to $\widehat{G}(r,m)^N$. The product of the resulting projective spaces for each $i$ can itself be embedded in a larger projective space via the Segre embedding. The key point is that the projective coordinates of $\widehat{Y}$ inside this big projective space exactly correspond to the non-zero coordinates $\widehat{y}_I$ and $\widehat{y}_{I}^{ij}$ satisfying $\widehat{\alpha}_I \cdot \beta(\mu') = ||\beta(\mu)||^2$ and $\widehat{\alpha}_{I}^{ij} \cdot \beta(\mu) = ||\beta (\mu)||^2$. Since the coordinates of $\widehat{Y}$ are projective, at least one is non-zero, and it corresponds therefore to a non-zero coordinate $\widehat{y}_{I}^{ij}$ or $\widehat{y}_I$ satisfying $\widehat{\alpha}_{I}^{ij} \cdot \beta (\mu) = ||\beta(\mu)||^2$ if it is $\widehat{y}_{I}^{ij}$ and $\widehat{\alpha}_I \cdot \beta (\mu) = ||\beta(\mu)||^2$ if it is $\widehat{y}_I$. This shows that $g \cdot \widehat{y} \in \widehat{Y}_{\beta(\mu)}$ as required. 
\end{proof}

\begin{lemma} \thlabel{lemma2}
Let $\mu=(d_1/r_1,\hdots, d_s/r_s)$ denote a HN type of rank $r$ and degree $d$. Then if the degrees $d_i$ and $N$ are sufficiently large, there exists $N$ points $x_1, \hdots, x_N \in \Sigma$ such that if $(E,\phi) \in F_{r,d}$ has HN type $ \mu$, then $g \cdot \widehat{\iota}(E,\phi) \in \widehat{Y}_{\beta(\mu)}^{ss}$ for some $g \in \GL(m)$, where $\widehat{\iota}: F_{r,d} \to \widehat{G}(r,m)^N$ is the map associated to $x_1,\hdots, x_N$. 
\end{lemma} 

\begin{proof} 
Suppose that $(E,\phi) \in F_{r,d}$ has HN type $\mu$. By \thref{lemma1}, for $d$ and $N$ large enough there is some $g \in \GL(m)$ such that  $ g \cdot \widehat{\iota}((E,\phi)) \in \widehat{Y}_{\beta(\mu)}$. For simplicity, let $g$ be the identity and write $\widehat{y}: = \widehat{\iota}(E,\phi)$. By definition of $\widehat{Y}_{\beta(\mu)}$, the point $\widehat{y}$ lies in $\widehat{Y}_{\beta(\mu)}$ if and only if $p_{\beta(\mu)}(\widehat{y})$ lies in $\widehat{Z}_{\beta(\mu)}^{ss}$. The image $p_{\beta(\mu)}(\widehat{y})$ coincides with the image $\widehat{\iota}(\operatorname{gr}(E,\phi))$, and the component of $Z_{\beta}$ containing $p_{\beta}(\widehat{y})$ can be identified with the product $\prod_{i=1}^s  \widehat{G}(r_{i}, m_{i})^N$. Note that under this identification $p_{\beta(\mu)}(\widehat{y})$ corresponds to the element $\widehat{Y}$ defined above at \eqref{bigyhat} in the proof of \thref{lemma1}. 
By the Hilbert-Mumford criterion, $p_{\beta}(\widehat{y}) \in \widehat{Z}_{\beta(\mu)}$ lies in $\widehat{Z}_{\beta(\mu)}^{ss}$ if for every one-parameter subgroup $\lambda$ of $\operatorname{Stab} \beta(\mu)$ we have that $\mu(\widehat{Y}_i, \lambda) \geq \lambda \cdot \beta(\mu)$, where $\mu(\widehat{Y}_i, \lambda)$ denotes the negative of the weight with which the one-parameter subgroup $\lambda$ acts on the fibre of the line bundle over $\operatorname{lim}_{t \mapsto 0} \lambda(t) \cdot \widehat{Y}_i$ (see \cite[Def 12.20]{Kirwan1984}).\footnote{In \cite[Def 12.20]{Kirwan1984} the inequality is the other way around. This is because the definition of the Hilbert-Mumford weight $\mu(-, \lambda)$ given there is the negative of the definition which we use in this paper, following \cite[p77]{Nitsure1991}.}

Since by definition $\widehat{Y}_i = \widehat{\iota}_i(E_i,\phi_i)$ and each $(E_i,\phi)$ is semistable, we can apply \cite[Prop 5.5]{Nitsure1991} to conclude that provided $d_i$ and $N$ are sufficiently large then $\widehat{Y}_i$ lies in the GIT-semistable locus for the action of $\SL(m_i)$ on $\widehat{G}(r_i,m_i)^N$. Therefore $\mu(\widehat{Y}_i, \lambda) \geq 0$ for every one-parameter subgroup $\lambda$ of $\operatorname{SL}(m_i)$. 

Now any one-parameter subgroup $\lambda$ of $\GL(m_i)$ can be written as the product of a one-parameter subgroup of $\SL(m_i)$ with a central one-parameter subgroup of $\GL(m_i)$ given by $t \mapsto (t^a, \hdots, t^a)$ for some $ a \in k$, where $(t^a,\hdots, t^a)$ represents the matrix with $t^a$ on each diagonal entry. By abuse of notation let $\lambda$ denote the corresponding element of the Lie algebra $\mathfrak{t}_i$ of the central torus $T_i$ in $\GL(m_i)$, so that $\lambda = (a,\hdots, a)$.  Then $(t^a,\hdots, t^a)$ acts on $\widehat{Y}_i$ with weight $Nr_ia = - k_i \operatorname{tr} \lambda/m_i.$ Therefore $\mu(\widehat{Y}_i,\lambda) = k_i \operatorname{tr} \lambda /m_i.$ 
It follows that any one-parameter subgroup $\lambda_i$ of $\GL(m_i)$, which we identify as an element of $\mathfrak{t}_i \cong k^{m_i}$, satisfies $\mu(\widehat{Y}_i,\lambda_i) =k_i \operatorname{tr} \lambda_i / m_i.$

Thus for a one-parameter subgroup $\lambda$ of $\prod_{i=1}^s \GL(m_i)$, where we label by $\lambda_1, \hdots, \lambda_s$ the blocks on the diagonal with each $\lambda_i$ determining an element of $\mathfrak{t}_i$, we have that $$\mu(\widehat{Y},\lambda) \geq  \sum_{i=1}^{s} \operatorname{tr} \lambda_i k_i =  \lambda \cdot \beta(\mu).$$ Hence $\widehat{Y} = p_{\beta}(\widehat{y}) \in Z_{\beta(\mu)}^{ss}$ as required. 
\end{proof}

We can now prove \thref{bigprop}, which follows easily from \thref{lemma1} and \thref{lemma2}. 

\begin{proof}[Proof of \thref{bigprop}] 
If $\mu=\mu_0$ then the result is true by \cite[Prop 5.5]{Nitsure1991}. Thus we can assume that $\mu$ is not the trivial HN type $\mu_0$. Given a HN type $\mu$ of rank $r$ and degree $d$, there is only a finite number of HN types $\mu'$ of rank $r$ and degree $d$ such that $\mu' \leq \mu$. Thus we can choose $d$ and $N$ large enough such that the results of \thref{lemma1} and \thref{lemma2} are valid for any Higgs bundle of HN type $\mu' \leq \mu$, rather than just those of HN type $\mu$. This yields the desired result. 
\end{proof} 

\subsection{The refined bHN stratification as a Bialynicki-Birula stratification} \label{subsec:refbHNasBB}

In this section we show that the refined bHN stratification of a given bHN stratum $\mathscr{H}_{\tau}^{r,d}(\Sigma,L)$ can be obtained from a Bialynicki-Birula (BB) stratification associated to the action of a one-parameter subgroup of $\GL(m)$ on a subvariety of $\widehat{G}(r,m)^N$. 

By \thref{bHNtoGIT}, we know that given a bHN stratum $\mathscr{H}_{\tau}^{r,d}(\Sigma,L)$  there is for $e$ and $N$ sufficiently large an isomorphism and inclusion of stacks $$\mathscr{H}_{\tau}^{r,d}(\Sigma,L) \cong \left[F_{\tau+e}^{r,d+re}/\GL(m+re) \right] \hookrightarrow \left[ \pi^{-1}(S_{\beta(\tau+e)}) / \GL(m+re) \right] \subseteq \left[ \widehat{G}(r,m+re)^N / \GL(m+re)\right].$$ As reviewed in Section \ref{subsubsec:GITinstab}, we have that $$S_{\beta(\tau+e)}= \GL(m+re) Y_{\beta(\tau)}^{ss} \cong \GL(m+re) \times_{P_{\beta(\tau+e)}} Y_{\beta(\tau)}^{ss}$$ where $P_{\beta(\tau+e)}$ is the parabolic subgroup of $\GL(m+re)$ associated to the maximally destabilising one-parameter subgroup $\lambda_{\beta(\tau+re)}(\GG_m)$ of $\GL(m+re)$. The action of $\lambda_{\beta(\tau+re)}$ on $\pi^{-1}(\overline{Y_{\beta(\tau+re)}^{ss}})$ yields a BB stratification of the space. \thref{higgstrat} below shows that taking the $\GL(m)$-sweep of this stratification and pulling it back to $\mathscr{H}_{\tau}^{r,d}(\Sigma,L)$ recovers the refined bHN stratification.

\begin{prop}\thlabel{higgstrat}
Let $\tau = (d_1/r_1,\hdots, d_s/r_s)$ denote a bHN type of rank $r$ and degree $d$. Then if the degrees $d_i$ and $N$ are sufficiently large, the refined bHN stratification of $\mathscr{H}_{\tau}^{r,d}(\Sigma,L)$ can be recovered from the stratification of $$\left[ \pi^{-1}(S_{\beta(\tau)}) / \GL(m) \right] = \left[ \GL(m)\pi^{-1}(Y_{\beta(\tau)}^{ss}) /\GL(m) \right]$$ induced by the Bialyinicki-Birula stratification for the action of $\lambda_{\beta(\tau)}(\GG_m)$ on $\pi^{-1}(\overline{Y_{\beta(\tau+re)}^{ss}}) \subseteq \widehat{G}(r,m)^N$, under the inclusion and isomorphism $$\mathscr{H}_{\tau}^{r,d}(\Sigma,L) \cong \left[F_{\tau}^{r,d}/\GL(m) \right] \hookrightarrow \left[ \pi^{-1}(S_{\beta(\tau)}) / \GL(m) \right].$$  
\end{prop}

\begin{proof}
By definition, the refined bHN stratification of $\mathscr{H}_{\tau}^{r,d}(\Sigma,L)$ is given by \begin{equation*}\mathscr{H}_{\tau}^{r,d}(\Sigma,L) = \bigsqcup_{1 \leq j < i  \leq s} \mathscr{H}_{\tau,[i,j]}^{r,d}(\Sigma,L) \sqcup \mathscr{H}_{\tau,[0,0]}^{r,d}(\Sigma,L). \label{HiggsrefforHNstrat}
\end{equation*}   
 
Suppose that the degrees $d_i$ appearing in $\tau$ are sufficiently large so that there is an isomorphism and inclusion $$\mathscr{H}_{\tau}^{r,d}(\Sigma,L) \cong \left[F_{\tau}^{r,d}/\GL(m) \right] \hookrightarrow \left[ \pi^{-1}(S_{\beta(\tau)}) / \GL(m) \right].$$ This is possible by \thref{qsforHN} and \thref{bHNtoGIT}. 
Let $X: =\pi^{-1}( \overline{Y_{\beta(\tau)}^{ss}})$. Then $$ \pi^{-1}(S_{\beta(\tau)}) = \pi^{-1}(S_{\beta(\tau)}) \cap \GL(m) X.$$ 

We start by showing that the open stratum of the Bialynicki-Birula stratification associated to the action of $\lambda_{\beta(\tau)}(\GG_m)$ on $X$ pulls back to give the open stratum $\mathscr{H}_{\tau,[s,1]}^{r,d}(\Sigma,L)$ of the refined bHN stratification of $\mathscr{H}_{\tau}^{r,d}(\Sigma,L)$. 

If $\GG_m$ acts linearly on a projective variety $X \subseteq \mathbb{P}(V)$ with respect to a representation $\GG_m \to \GL(V)$ (this is the case for the action of $\lambda_{\beta(\tau)}(\GG_m)$ on $X$), then the open stratum can be described as follows. The representation $V$ decomposes into a sum of $\GG_m$-weight spaces $\bigoplus_{i \in \mathbb{Z}} V_i$. If we let $V_{\operatorname{min}}$ denote the non-zero weight space with minimal index and $\zmin := X \cap \mathbb{P}(V_{\operatorname{min}})$, then the open Bialynicki-Birula stratum is $$X_{\operatorname{min}}: =  \{x \in X  \ | \ \operatorname{lim}_{t \mapsto 0} t \cdot x \in \zmin \}.$$ In the case we are considering, the action of $\lambda_{\beta(\tau)}(\GG_m)$ on $\pi^{-1}(\overline{Y_{\beta(\tau)}^{ss}}) \subseteq \widehat{G}(r,m)^N \subseteq \mathbb{P}(\widehat{V}^{\vee})$ is linearised via the restriction of the representation $\GL(m) \to \widehat{V}^{\vee}$ defined in Section \ref{subsubsec:GLactiononG-hat}.  

For a point $\widehat{y}=(\langle y_1, [c_1{:}\phi_1] \rangle, \hdots, y_N,[ c_N{:}\phi_N]\rangle) \in \pi^{-1}(\overline{Y_{\beta(\tau)}^{ss}})$ to be in the minimal weight space for $\pi^{-1}(\overline{Y_{\beta(\tau)}^{ss}})$, it must be fixed by $\lambda_{\beta(\tau)}(\GG_m)$. Thus in particular its projection to $G(r,m)^N$ under $\pi$ must be fixed by $\lambda_{\beta(\tau)}(\GG_m)$, since $\pi$ is $\GL(m)$-equivariant, and moreover it must lie in the minimal weight space for $\overline{Y_{\beta(\tau}}^{ss})$. But by \cite[Prop 4.2.0.2]{Jackson2018}, this minimal weight space is precisely $Z_{\beta(\tau)}^{ss}$.  
Thus to determine the minimal weight space in $\pi^{-1}(\overline{Y_{\beta(\tau)}^{ss}})$, it suffices to determine the minimal weight space for the action of $\lambda_{\beta(\tau)}(\GG_m)$ on the fibre $\pi^{-1}(z)$ for any $z \in Z_{\beta(\tau)}^{ss}$. We fix such a point $z = (\langle z_1 \rangle, \hdots, \langle z_N \rangle)$. Then any point $\widehat{z} \in \pi^{-1}(z)$ is of the form $$\widehat{z} = ( \langle  z_1, [c_1{:} \phi_1] \rangle, \cdots, \langle z_N, [c_N{:} \phi_N] \rangle ).$$ By definition of the action of $\GL(m)$ on $\widehat{G}(r,m)^N$, the one-parameter subgroup $\lambda_{\beta(\tau)
}(\mathbb{G}_m)$ acts trivially on the coordinates $c_i$, and via conjugation on the matrices $\phi_i$. The weights of this action are therefore given by $$\{0\} \cup \left\{ \left. \frac{k_j}{m_j} - \frac{k_i}{m_i}  \ \right| \ i,j \in \{1,\hdots, s\}, i \neq j \right\}.$$ Since $k_1/m_1 > k_2/m_2 > \cdots > k_s/m_s$, it follows that the minimal weight is $k_s/m_s - k_1/m_1$. 

Given $\phi \in \operatorname{Mat}_{r \times r}(k)$, we can consider for each $i,j \in \{1,\hdots, s\}$ the block $\phi_{ij} \in  \operatorname{Mat}_{r_i \times r_j}(k)$. Matrices $\phi$ which are zero everywhere except for the bottom left entry $\phi_{s1}$ are weight vectors for the conjugation action of $\lambda_{\beta(\tau)}(\mathbb{G}_m)$ on $\operatorname{Mat}_{r \times r}(k)$ with weight $k_s/{m_s}- k_1/m_1$.  Thus the minimal weight space for the action of $\lambda_{\beta(\tau)}(\GG_m)$ on $\pi^{-1}(\overline{Y_{\beta(\tau)}^{ss}})$ consists of points $(\langle z_1, [0{:} \phi_1] \rangle, \cdots, \langle z_N , [0{:} \phi_N] \rangle) \in \pi^{-1}(\overline{Y_{\beta}^{ss}})$ such that $( \langle z_1 \rangle, \hdots, \langle z_N \rangle) \in Z_{\beta(\tau)}^{ss}$ and such that for each $k = 1 ,\hdots, N$ the equality $(\phi_k)_{ij} = 0$ holds for $i,j \in \{1,\hdots, s\}$ except when $i = s$ and $j=1$ (note that this forces $(\phi_k)_{s1} \neq 0$). 

It follows that the open stratum of the BB stratification corresponding to this minimal weight space consists of points $(\langle y_1, [0{:} \phi_1] \rangle, \cdots, \langle y_N , [0{:} \phi_N] \rangle) \in \pi^{-1}(\overline{Y_{\beta}^{ss}})$ such that $( \langle y_1 \rangle, \hdots, \langle y_N \rangle) \in Y_{\beta}^{ss}$ and such that $(\phi_k)_{s1} \neq 0$ for each $k = 1, \hdots, N$. The intersection of this open stratum with the image of $F_{\tau}^{r,d}$ under $\widehat{\iota}$ can be interpreted in the following moduli-theoretic way: it consists of Higgs bundles $(E,\phi)$ with Harder-Narasimhan type $\tau$ such that $\phi_{s1} \neq 0.$ Thus the resulting open stratum of $\mathscr{H}_{\tau}^{r,d}(\Sigma,L)$ coincides with the open stratum $\mathscr{H}_{\tau,[s,1]}^{r,d}(\Sigma,L)$ of the refined stratification of $\mathscr{H}_{\tau}^{r,d}(\Sigma,L)$.

The proof that the BB strata for higher weights $k_j/m_j - k_i/m_i < 0$ coincide with the strata $\mathscr{H}_{\tau,[i,j]_{\tau}}(\Sigma,L)$ is completely analogous.

It remains therefore only to relate the highest stratum $\mathscr{H}_{\tau,[0,0]}^{r,d}(\Sigma,L) = \mathscr{H}_{\tau}^{r,d}(\Sigma,L) \cap \mathscr{H}^{\tau}_{r,d}(\Sigma,L)$ of the Higgs stratification to a BB stratum. We show that it is the stratum corresponding to the zero weight space. Indeed, in the  zero weight space for the action of $\l_{\beta(\tau)}(\GG_m)$ on $\pi^{-1}(\overline{Y_{\beta(\tau)}^{ss}})$, points $\widehat{y} = (\langle y_1, [c_1{:} \phi_1] \rangle, \cdots, \langle y_N , [c_N{:} \phi_N] \rangle) \in \pi^{-1}(\overline{Y_{\beta(\tau)}^{ss}})$ in the corresponding BB stratum satisfy the property that each $\phi_i$ is a block upper triangular matrix (to ensure that there is no non-zero coordinate with smaller weight). 

Thus if a Higgs bundle lies in this Bialynicki-Birula stratum, then its Higgs field must preserve the HN filtration (this follows from the fact that each $\phi_i$ is block upper triangular), so that its HN and bHN types coincide. Conversely, if a Higgs bundle satisfies this property, then its image in $\pi^{-1}(\overline{Y_{\beta(\tau)}^{ss}})$ will have zero weight simply based on the fact that it lies in $F_{r,d}$. Indeed, the inclusion $F_{\tau}^{r,d} \hookrightarrow \widehat{G}(r,m)^N$ maps $F_{\tau}^{r,d}$ into the set of points $(\langle y_1, [c_1{:} \phi_1] \rangle, \cdots, \langle y_N , [c_N{:} \phi_N] \rangle) \in \widehat{G}(r,m)^N$ such that all coordinates $c_i$ are non-zero, and these coordinates have weight zero for the action of $\l_{\beta(\tau)}(\GG_m)$. Thus the intersection with $F_{\tau}^{r,d}$ of the BB stratum associated to the zero weight space for the action of $\l_{\b(\tau)}$ on $\pi^{-1}(\overline{Y_{\beta(\tau)}^{ss}})$  corresponds inside the HN stratum $\mathscr{H}_{\tau}^{r,d}(\Sigma,L)$ to the intersection $\mathscr{H}_{\tau}^{r,d}(\Sigma,L) \cap \mathscr{H}^{\tau}_{r,d}(\Sigma,L)$.
\end{proof}


\begin{thebibliography}{10}

\bibitem{Atiyah1983}
{\sc Atiyah, M., and Bott, R.}
\newblock The {Y}ang-{M}ills equations over {R}iemann surfaces.
\newblock {\em Philosophical Transactions of the Royal Society 308\/} (1983),
  523--615.

\bibitem{Berczi2018}
{\sc B{\'e}rczi, G., Hoskins, V., and Kirwan, F.}
\newblock Stratifying quotient stacks and moduli stacks.
\newblock In {\em Geometry of Moduli\/} (Cham, 2018), J.~A. Christophersen and
  K.~Ranestad, Eds., Springer International Publishing, pp.~1--33.

\bibitem{Casalaina-Martin2018}
{\sc Casalaina-Martin, S., and Wise, J.}
\newblock {\em {The Geometry, Topology and Physics of Moduli Spaces of Higgs
  Bundles}}, vol.~36 of {\em Lecture Note Series, Institute for Mathematical
  Sciences, National University of Singapore}.
\newblock World Scientific Publishing Company, 09 2018, ch.~{An Introduction to
  Moduli Stacks, with a View towards Higgs Bundles on Algebraic Curves},
  pp.~199--399.

\bibitem{Collier2019}
{\sc Collier, B., and Wentworth, R.}
\newblock {Conformal limits and the Bialynicki-Birula stratification of the
  space of $\lambda$-connections}.
\newblock {\em Advances in Mathematics 350\/} (2019), 1193--1225.

\bibitem{Fan2020}
{\sc Fan, Y.}
\newblock {Construction of the Moduli Space of Higgs Bundles using Analytic
  Methods}.
\newblock arXiv:2004.07182 [math.DG], 2020.

\bibitem{Fan2022}
{\sc Fan, Y.}
\newblock An analytic approach to the quasi-projectivity of the moduli space of
  higgs bundles.
\newblock {\em Advances in Mathematics 406\/} (2022), 108506.

\bibitem{Franco2014}
{\sc Franco, E., Garcia-Prada, O., and Newstead, P.~E.}
\newblock Higgs bundles over elliptic curves.
\newblock {\em Illinois J. Math. 58}, 1 (2014), 43--96.

\bibitem{Garcia-Prada2008}
{\sc Garcia-Prada, O.}
\newblock {Moduli spaces and Geometric Structures, Appendix to third edition of
  Differential Analysis on Complex Manifolds by Raymond O. Wells}.
\newblock {\em Graduate Texts in mathematics\/} (2008).

\bibitem{GarciaPrada2011}
{\sc Garcia-Prada, O., Heinloth, J., and Schmitt, A.}
\newblock On the motives of moduli of chains and {H}iggs bundles.
\newblock {\em Journal of the European Mathematical Society 16\/} (04 2011).

\bibitem{Gothen2015}
{\sc Gothen, P., and Zuniga-Rojas, R.}
\newblock {Stratifications on the Moduli Space of Higgs Bundles}.
\newblock {\em Portugaliae Mathematica 74\/} (11 2015).

\bibitem{Gurjar2014}
{\sc Gurjar, S., and Nitsure, N.}
\newblock {Schematic Harder-Narasimhan stratification for families of principal
  bundles and $\Lambda$ -modules}.
\newblock {\em Proc. Indian Acad. Sci. (Math. Sci.) 124}, 3 (August 2014),
  315--332.

\bibitem{Halpern-Leistner2015}
{\sc Halpern-Leistner, D.}
\newblock $\theta$-stratifications, $\theta$-reductive stacks, and
  applications.
\newblock In {\em Proceedings of the 2015 AMS Summer Institute in Salt Lake
  City\/} (2015).
\newblock arXiv:1411.0627.

\bibitem{Halpern-Leistner2016}
{\sc Halpern-Leistner, D.}
\newblock {The equivariant Verlinde formula on the moduli of Higgs bundles}.
\newblock arXiv: 1608.01754 [math.AG], 8 2016.

\bibitem{Halpern-Leistner2018}
{\sc Halpern-Leistner, D.}
\newblock On the structure of instability in moduli theory.
\newblock arXiv:1411.0627v4 [math.AG], July 2018.

\bibitem{Hausel1998b}
{\sc Hausel, T.}
\newblock {\em {Geometry of the moduli space of Higgs bundles}}.
\newblock PhD thesis, University of Cambridge, 1998.

\bibitem{Hausel2003}
{\sc Hausel, T., and Thaddeus, M.}
\newblock Mirror symmetry, {L}anglands duality, and the {H}itchin system.
\newblock {\em Inventiones mathematicae 153}, 1 (July 2003), 197--229.

\bibitem{Hausel2004}
{\sc Hausel, T., and Thaddeus, M.}
\newblock Generators for the cohomology ring of the moduli space of rank 2
  higgs bundles.
\newblock {\em Proceedings of the London Mathematical Society 88}, 3 (2004),
  632–658.

\bibitem{Hesselink1978}
{\sc Hesselink, W.~H.}
\newblock Uniform instability in reductive groups.
\newblock {\em Journal für die reine und angewandte Mathematik 0303\_0304\/}
  (1978), 74--96.

\bibitem{Hitchin1987}
{\sc Hitchin, N.}
\newblock {The self-duality equations on a {R}iemann surface}.
\newblock {\em Proceedings of the London Mathematical Society 3}, 55 (1987),
  59--126.

\bibitem{Hoskins2014}
{\sc Hoskins, V.}
\newblock Stratifications for moduli of sheaves and moduli of quiver
  representations.
\newblock arXiv:1407.4057 [math.AG], July 2014.

\bibitem{Hoskins2012}
{\sc Hoskins, V., and Kirwan, F.}
\newblock Quotients of unstable subvarieties and moduli spaces of sheaves of
  fixed {H}arder-{N}arasimhan type.
\newblock {\em Proceedings of the London Mathematical Society 3}, 105 (2012),
  852--890.

\bibitem{Huang2020}
{\sc Huang, P.}
\newblock {\em {Th\'eorie de Hodge non-abélienne et des spécialisations}}.
\newblock PhD thesis, {Laboratoire de Math\'ematiques J. A. Dieudonn\'e}, 06
  2020.

\bibitem{Jackson2018}
{\sc Jackson, J.}
\newblock {\em Moduli {S}paces of {U}nstable {C}urves and {S}heaves via
  {N}on-{R}eductive {G}eometric {I}nvariant {T}heory}.
\newblock PhD thesis, University of Oxford, 2018.

\bibitem{Jackson2021}
{\sc Jackson, J.}
\newblock {Moduli of unstable objects: sheaves of Harder-Narasimhan length
  two}.
\newblock arXiv:2111.07428 [math.AG], 2021.

\bibitem{Hoskins2021}
{\sc Jackson, J., and Hoskins, V.}
\newblock Quotients by parabolic groups and moduli of unstable objects.
\newblock arXiv:2111.07429 [math.AG].

\bibitem{Kempf1978}
{\sc Kempf, G.~R.}
\newblock Instability in invariant theory.
\newblock {\em Annals of Mathematics 108}, 2 (1978), 299--316.

\bibitem{Kirwan1984}
{\sc Kirwan, F.}
\newblock {\em Cohomology of {Q}uotients in {S}ymplectic and {A}lgebraic
  {G}eometry}.
\newblock Mathematical Notes, Princeton University Press, 1984.

\bibitem{Kirw1985}
{\sc Kirwan, F.}
\newblock On spaces of maps from {R}iemann surfaces to {G}rassmannians and
  applications to the cohomology of moduli of vector bundles.
\newblock {\em Arkiv f\"or Matematik 24}, 1-2 (1985), 221--275.

\bibitem{Mumford1963}
{\sc Mumford, D.}
\newblock {\em {Projective Invariants of Projective Structures and
  Applications}}.
\newblock 01 1963, pp.~526--530.

\bibitem{Mumford1994}
{\sc Mumford, D., Fogarty, J., and Kirwan, F.}
\newblock {\em {Geometric Invariant Theory}}.
\newblock Springer, 1994.

\bibitem{Ness1984}
{\sc Ness, L., and Mumford, D.}
\newblock A stratification of the null cone via the moment map.
\newblock {\em American Journal of Mathematics 106\/} (12 1984).

\bibitem{Newstead1978}
{\sc Newstead, P.}
\newblock {\em Introduction to {M}oduli {P}roblems and {O}rbit {S}paces}.
\newblock Narosa Publishing House, 1978 (reprint 2012).

\bibitem{Nitsure1991}
{\sc Nitsure, N.}
\newblock Moduli space of semistable pairs on a curve.
\newblock {\em Proceedings of the London Mathematical Society s3}, 62 (March
  1991), 275--300.

\bibitem{Nitsure2011}
{\sc Nitsure, N.}
\newblock {Schematic Harder-Narasimhan Stratification}.
\newblock {\em International Journal of Mathematics 22}, 10 (2011), 1365--1373.

\bibitem{Ramanan1989}
{\sc Ramanan, S., Narasimhan, M., and Beauville, A.}
\newblock Spectral curves and the generalised {T}heta divisor.
\newblock {\em Journal für die reine und angewandte Mathematik 398\/} (1989),
  169--179.

\bibitem{Sancho2022}
{\sc Sancho, {\'A}.~A.}
\newblock {Shatz and Bialynicki-Birula stratifications of the moduli space of
  Higgs bundles}.
\newblock {\em Hokkaido Mathematical Journal 51}, 1 (2022), 25 -- 56.

\bibitem{Schaposnik2020}
{\sc Schaposnik, L.}
\newblock {Higgs Bundles -- Recent Applications}.
\newblock {\em Notices of the American Mathematical Society 67}, 5 (May 2020),
  625--634.

\bibitem{Seshadri1967}
{\sc Seshadri, C.~S.}
\newblock {Space of Unitary Vector Bundles on a Compact Riemann Surface}.
\newblock {\em Annals of Mathematics 85}, 2 (1967), 303--336.

\bibitem{Shatz1977}
{\sc Shatz, S.}
\newblock The decomposition and specialization of algebraic families of vector
  bundles.
\newblock {\em Compositio Mathematica 35}, 2 (1977), 163--187.

\bibitem{Simpson1994a}
{\sc Simpson, C.}
\newblock Moduli of representations of the fundamental group of a smooth
  projective variety {II}.
\newblock {\em Publications Math\'ematiques de l'IH\'ES 80\/} (1994), 5--79.

\bibitem{Simpson1994}
{\sc Simpson, C.~T.}
\newblock Moduli of representations of the fundamental group of a smooth
  projective variety {I}.
\newblock {\em Publications math\'ematiques de l'IH\'ES 79\/} (1994).

\bibitem{Swoboda2021}
{\sc Swoboda, J.}
\newblock {Moduli Spaces of Higgs Bundles – Old and New}.
\newblock {\em Jahresbericht der Deutschen Mathematiker-Vereinigung 123\/} (03
  2021), 1--66.

\bibitem{Wentworth2016}
{\sc Wentworth, R.}
\newblock {\em Geometry and {Q}uantization of {M}oduli {S}paces}.
\newblock Birkh\"auser, Cham, 2016, ch.~Higgs {B}undles and {L}ocal {S}ystems
  on {R}iemann {S}urfaces.

\bibitem{Wilkin2006}
{\sc Wilkin, G.}
\newblock {Morse Theory for the Space of Higgs Bundles}.
\newblock {\em Communications in Analysis and Geometry 16\/} (12 2006).

\end{thebibliography}

\end{document}